\documentclass[10pt,reqno]{amsart}

\usepackage{amsmath}%
\usepackage{amsthm}
\usepackage{amsfonts}%
\usepackage{amssymb,latexsym}%
\usepackage{graphicx}
\usepackage[english]{babel}
\usepackage{syntonly}
\usepackage{color}  
\usepackage[all]{xy}
\usepackage{cancel}
\usepackage[only,llbracket,rrbracket]{stmaryrd}
\usepackage{mathrsfs}  
\usepackage{tipa}  
\usepackage{mathtools}  

\usepackage{xspace}
\usepackage{textcase}

\usepackage{multicol} 

\usepackage{lscape}
\usepackage{rotating} 

\DeclareMathAlphabet\mathbfcal{OMS}{cmsy}{b}{n} 

\numberwithin{equation}{section}

\usepackage{subfiles} 

\makeatletter                                      
\newcommand{\labitem}[2]{                     
\def\@itemlabel{#1}                     
\item					                      
\def\@currentlabel{#1}\label{#2}}          
\makeatother			                     

\usepackage{tikz} 
\usetikzlibrary{arrows,shapes,patterns,fit,backgrounds}
\tikzstyle{every picture}+=[remember picture]

\newtheorem{thm}{Theorem}[section]
\newtheorem{lemma}[thm]{Lemma}
\newtheorem{prop}[thm]{Proposition}

\newtheorem{cor}[thm]{Corollary}

\theoremstyle{definition}
\newtheorem{defn}[thm]{Definition}

\theoremstyle{remark}
\newtheorem{rem}[thm]{Remark}

\newcommand{\myhat}[1]{\widehat{#1}}

\newcommand{\ie}{i.e., }

\newcommand{\cf}{cf.\ }

\newcommand{\clas}[1]{\mathbb{#1}} 
\newcommand{\alg}[1]{\mathbf{#1}}
\newcommand{\spa}[1]{\mathfrak{#1}}
\newcommand{\ralg}[1]{\mathcal{#1}}
\newcommand{\Fig}[1]{{\llbracket}#1{\rangle\kern-.6ex\rangle}} 
\newcommand{\Idfg}[1]{{\langle\kern-.6ex\langle}#1{\rrbracket}} 
\newcommand{\Fimg}[1]{\Fig{#1}} 
\newcommand{\Fimgg}{\Fimg{\phantom{.}}} 
 
\newcommand{\Idfgg}{\Idfg{\phantom{.}}}

\newcommand{\clopu}{\mathrm{ClUp}}  
\newcommand{\aclopu}{\mathbf{ClUP}}

\newcommand{\taua}{\tau_{\Al}}

\newcommand{\eps}{\varepsilon} 
\newcommand{\nvarphi}{\Phi}

\newcommand{\Xb}{X_{B}}

\newcommand{\Xbone}{X_{B_{1}}}
\newcommand{\Xbtwo}{X_{B_{2}}}

\newcommand{\xis}{\xi}

\newcommand{\xisa}{T}

\newcommand{\nxisa}{\Xi}

\newcommand{\subsw}{\subseteq^{\nat}}
\newcommand{\pws}{\mathcal{P}}

\newcommand{\ups}{\mathcal{P}^{\uparrow}}

\newcommand{\funid}{\mathrm{id}}
\newcommand{\inc}{\bot}

\newcommand{\Sm}{\mathcal{S}}
\newcommand{\vdashs}{\vdash_{\mathcal{S}}}

\newcommand{\Fm}{\mathbf{Fm}}
\newcommand{\Hom}{\mathrm{Hom}}   
\newcommand{\nat}{\omega}

\newcommand{\Fi}{\mathrm{Fi}}
\newcommand{\lFi}{\alg{Fi}}

\newcommand{\Irr}{\mathrm{Irr}}

\newcommand{\Opt}{\mathrm{Op}}
\newcommand{\sOpt}{\mathfrak{O}\mathrm{p}}

\newcommand{\Id}{\mathrm{Id}}

\newcommand{\FF}{\mathcal{F}}
\newcommand{\II}{\mathcal{I}}
\newcommand{\varphif}{\varphi_{\FF}}

\newcommand{\varphifm}{\myhat{\varphi}_{\FF}}
\newcommand{\varphimh}{\myhat{\varphi}}

\newcommand{\varphim}{\myhat{\varphi}}  
\newcommand{\varphifA}{\varphif[\Al]}

\newcommand{\varphib}{\varphi_{\BB}}

\newcommand{\MefA}{\mathrm{M}_{\FF}(\Al)}

\newcommand{\Mee}{\mathrm{M}}
\newcommand{\Meea}{\Mee(A)}
\newcommand{\MeeA}{\Mee(\Al)}

\newcommand{\MeeB}{\Mee(\BB)}

\newcommand{\PP}{\mathbf{P}}  

\newcommand{\IdF}{\Id_{F}}

\newcommand{\LL}{\mathbf{L}}  

\newcommand{\MM}{\mathbf{M}}   
\newcommand{\Fim}{\Fi}

\newcommand{\Irrm}{\Irr_{\wedge}}

\newcommand{\Optm}{\Opt}

\newcommand{\sOptm}{\sOpt_{\wedge}}

\newcommand{\Al}{\mathbf{A}}  
\newcommand{\BB}{\mathbf{B}}   
\newcommand{\Alg}{\mathsf{Alg}}

\newcommand{\FiSA}{\Fi_{\Sm}(\Al)}
\newcommand{\FiSB}{\Fi_{\Sm}(\BB)}
\newcommand{\FiS}{\Fi_{\Sm}}

\newcommand{\lFiSA}{\lFi_{\Sm}(\Al)}

\newcommand{\IrrSA}{\Irr_{\Sm}(\Al)}
\newcommand{\IrrSB}{\Irr_{\Sm}(\BB)}

\newcommand{\OptSA}{\Opt_{\Sm}(\Al)}
\newcommand{\OptSB}{\Opt_{\Sm}(\BB)}
\newcommand{\OptS}{\Opt_{\Sm}}%
\newcommand{\sOptS}{\mathfrak{O}\mathrm{p}_{\Sm}}

\newcommand{\IdSA}{\Id_{\Sm}(\Al)}

\newcommand{\IdsS}{\Id_{s\Sm}}%
\newcommand{\IdsSA}{\Id_{s\Sm}(\Al)}

\newcommand{\cAlg}{\mathsf{Alg}}
\newcommand{\cPr}{\mathsf{Pr}}

\newcommand{\Cons}{\mathrm{Fg}_{\Sm}}
\newcommand{\Consa}{\mathrm{Fg}_{\Sm}^{\Al}}
\newcommand{\Conspa}{\mathrm{Fg}_{\Sm}^{\varphifA}}

\newcommand{\equivsa}{\equiv_{\Sm}^{\Al}}
\newcommand{\equivsb}{\equiv_{\Sm}^{\BB}}
\newcommand{\equivsfm}{\equiv_{\Sm}^{\Fm}}
\newcommand{\leqsa}{\leq_{\Sm}^{\Al}}

\newcommand{\Consb}{\mathrm{Fg}_{\Sm}^{\BB}}

\newcommand{\leqsb}{\leq_{\Sm}^{\BB}}
\newcommand{\geqsa}{\geq_{\Sm}^{\Al}}

\newcommand{\Lee}{\mathrm{L}}

\newcommand{\DSp}{generalized Priestley\xspace}

\newcommand{\Sp}{$\Sm$-Priestley\xspace}

\newcommand{\AAAND}{\,\,\,\,\,\,\,\,\,\,\,\,\,\,\,\,\,\,}

\newcommand{\IFF}{\,\,\, \text{ iff }\,\,\,}

\newcommand{\pseudoddt}{E4}

\usepackage[pdftex,bookmarks,colorlinks]{hyperref}

\title{Priestley-style duality for filter-distributive congruential logics}

\author{Mar\'ia Esteban}
\address{Departament de Filosofia, Universitat de Barcelona (UB). Montalegre 6, 08001 Barcelona, Spain }\email{mariaegmp@gmail.com}
\thanks{}

\author{Ramon Jansana}
\address{Departament de Filosofia, Universitat de Barcelona (UB). Montalegre 6, 08001 Barcelona, Spain}
\email{jansana@ub.edu}
\thanks{}
\begin{document}

\date{\today}

\begin{abstract} 
 We first present a Priestley-style dualitiy for the classes of algebras that are the algebraic counterpart  of some  congruential, finitary and  filter-distributive logic with theorems. Then we analyze which properties of the dual spaces correspond to properties that the logic might enjoy, like the deduction theorem or the existence of a disjunction.
\end{abstract}

\maketitle

\section{Introduction}

The classes of algebras that correspond to many well-known  logics have a distributive lattice reduct. Among them we find  Boolean algebras, Heyting algebras, modal algebras, positive modal algebras, De Morgan algebras and MV-algebras. These classes of algebras are also the algebraic counterpart of some congruential logic\footnote{We use `congruential logic'  in accordance with \cite{GeJaPa10} instead of the more common `fully selfextensional logic' we find in abstract algebraic logic.}  equal or closely related to the  logic from which they originally arise. This fact can be taken to explain from a logical perspective why topological Priestley dualities exist for many of them because the prime filters of the algebras are in fact the irreducible \textit{logical filters} of the congruential logic.  A crucial property of these congruential  logics is  that in any of their algebras the lattice of logical filters is distributive. In abstract algebraic logic the logics with this property are known  as filter-distributive.

Besides the  logics whose algebras have a distributive lattice reduct,  there are logics which are congruential and filter-distributive  but whose algebras have only a meet-semilattice  or a join-semilattice reduct or even no semilattice reduct at all;  for example Hilbert algebras, which  are the algebras that constitute   the algebraic counterpart of the implication fragment of intuitionistic logic, a fragment which  is itself a congruential and filter-distributive logic. 

 Our aim in this paper  is to  develop a framework to obtain Priestley-style dualities for the classes of algebras that are the algebraic counterparts  of the congruential, finitary and  filter-distributive logics. The  point of view we take is that of logic. The starting point is any logic $\Sm$ with the mentioned properties, its algebraic counterpart, namely the canonical class   $\Alg\Sm$ of algebras associated with it in abstract algebraic logic,  and the lattices of the logical filters of the algebras in $\Alg\Sm$. In  these lattices  we have its irreducible elements, which are also prime because the lattice is distributive. In general, these irreducible logical filters are not enough to be the points of a dual space if we are interested in a Priestley-style duality.    We need a less restrictive notion encompassing the irreducible logical filters. To introduce it, we consider the notion of strong logical  ideal, that comes from a generalization of the notion of strong Frink ideal for Hilbert algebras introduced in \cite{CeJa12}, and define the optimal logical filters as the logical filters whose complement is a strong logical ideal.  The optimal filters will be the points of the dual space.
 
A way to understand the optimal logical filters of an algebra $\Al \in \Alg\Sm$ is to consider the $\Sm$-semilattice of $\Al$,  a notion  introduced in \cite{GeJaPa10}. It is  the dual of the join-semilattice of the finitely generated logical filters of $\Al$. The lattice of the filters of this meet-semilattice turns out to be isomorphic to the lattice of the logical filters of $\Al$ and therefore inherits the distributivity from this later one.  Hence,  the $\Sm$-semilattice of $\Al$ is a distributive meet-semilattice, in the sense of  \cite[Sec.\ II.5]{Gr78}. Therefore  we can use and take inspiration from  the duality   for distributive meet-semilattices developed in \cite{BeJa11} to obtain the dualities we are after.

For every finitary, congruential and filter distributive logic $\Sm$ with theorems we present  a duality between the category of the algebras in $\Alg\Sm$ together with the algebraic homomorphisms between them and a category of  Priestly-style spaces augmented with an algebra of clopen up-sets; such structures will be called \Sp spaces. Then we characterize  properties of the the category of \Sp spaces that correspond to basic logical properties that the logic $\Sm$ might enjoy. We do it for the property of having a binary formula that behaves like a conjunction, the property of having a set of binary formulas that collectively behaves as a disjunction, the property of having a binary formula that behaves like an implication that satisfies the \textit{modus ponens} rule and the deduction theorem, and finally,  the property of having an inconsistent formula, \ie a formula that implies every formula. 

The paper is structured as follows. In Section 2 we present the preliminaries we need on posets, distributive meet-semilattices, and congruential logics. We also review the duality given in \cite{BeJa11}. Section 3 is devoted to the representation theorems for $\Sm$-algebras that we obtain using the optimal filters. In Section 4 we  study the $\Sm$-semilattice of an $\Sm$-algebra. In Section 5 we introduce the dual objects  of the $\Sm$-algebras and in Section 6  the duals of the homomorphisms between  $\Sm$-algebras. Section 7 shows the categorical duality. Finally,  in Section 8 we present the results on the dual properties of the logical properties we mentioned above.   

\section{Preliminaries}

\subsection{Notation}
Let $X, Y$ be  sets. For any $B\subseteq X$, we denote by  $B^{c}$ the relative complement of $B$ w.r.t.\ $X$  when no confusion can arise, i.e., $B^{c} =\{a\in X:a\notin B\}$. We denote the powerset of a set $X$ by $\mathcal{P}(X)$. We indicate that $B$ is a finite subset of $X$ by $B \subsw X$.

For any binary relation $R\subseteq X\times Y$ and every $x \in X$ we let $R(x) := \{y \in Y: xRy\}$, and we define  the function $\Box_{R}:\pws(Y)\to \pws(X)$  by setting for every $U \subseteq Y$ 
$$\Box_{R}(U)\coloneqq\{x\in X: R(x) \subseteq U\}.$$

For algebras $\Al$ and $\BB$ of the same type,  we denote  by $\Hom(\Al,\BB)$ the set of all homomorphisms from $\Al$ to $\BB$. 

If $\langle X, \tau, \leq\rangle$ is an ordered topological space, $\clopu(X)$ denotes the set of all its clopen up-sets. 

\subsection{Posets and distributive meet-semilattices}

Let $\PP = \langle P,\leq\rangle$ be a partial ordered set (a poset for short). 
For every $a\in P$, we let ${\uparrow} a \coloneqq\{b\in P:a\leq b\}$ and  ${\downarrow} a\coloneqq\{b\in P:b\leq a\}$. 
For every $U\subseteq P$, we define ${\uparrow}U\coloneqq \bigcup\{{\uparrow}a:a\in U\}$ and ${\downarrow} U\coloneqq\{{\downarrow}a:a\in U\}$. 
Moreover, we say that $U\subseteq P$ is an \emph{up-set} of $\PP$ (resp.\ a \emph{down-set}) when ${\uparrow}U=U$ (resp.\ ${\downarrow}U=U$). 
By $\ups(P)$ we denote the collection of all up-sets of $\PP$, and for 
 $U\subseteq P$, we denote by $\max(U)$  the collection of maximal elements of $U$. 
 A set $U\subseteq P$ is \emph{up-directed} (resp.\ \emph{down-directed}) when for every $a,b\in U$ there exists $c\in U$ such that $a,b\leq c$  (resp.\ $c\leq a,b$). 
 
 An \emph{order filter} of $\PP$ is any non-empty up-set that is down-directed. Dually,  an \emph{order ideal} is any non-empty down-set that is up-directed. By 
$\Id(\PP)$ we denote the collection of the order ideals  of $\PP$ and by $\Fi(\PP)$ the collection of its order filters. These collections may not be closure sytsems, i.e., they may not be closed under intersections of arbitrary subfamilies. 
 A \emph{Frink ideal}  
 of $\PP$, a concept introduced in  \cite{Fr54}, is a set $I\subseteq M$ such that for every finite $I'\subseteq I$ and every $b\in M$, $\bigcap\{{\uparrow}a:a\in I'\}\subseteq {\uparrow}b$ implies $b\in I$; in other words: if every lower bound of the set of upper bounds of $I'$ belongs to $I$.
We denote by $\IdF(\PP)$ the collection   of all Frink-ideals of $\PP$, which  is  a closure system, and  by $\Idfgg$ its associated  closure operator, 
\ie for every $B\subseteq M$, $\Idfg{B}$ is the least Frink ideal containing $B$ that can be  described as follows
\begin{equation}
a \in\Idfg{B} \hspace{0.1cm}  \text{ iff }  \hspace{0.1cm} \text{there exists a finite } B' \subseteq B \text{ s.t.}  \bigcap_{b \in B'}{\uparrow}b\subseteq {\uparrow}a.
\end{equation}
Notice that according to the definition the emptyset  may be  a Frink ideal, 
but this happens if and only if there is no bottom element in $\PP$.


An algebra $\MM=\langle M,\wedge,1\rangle$ of type $(2,0)$ is a \emph{meet-semilattice with top}  when the binary operation $\wedge$ is idempotent, commutative, associative, and $a\wedge1=1$ for all $a\in M$. 
The (meet) partial order  of  $\MM$ is the relation $\leq$ such that for every $a,b\in M$, $a\leq b$ if and only if  $a\wedge b=a$. 
Then for every $a, b \in M$,  $ a \wedge b$ is the  meet of $a, b$ w.r.t.\ $\leq$ and $1$ is its greatest upper bound. 
The meet-semilattices with  top  that also have  a least lower bound  are called \emph{bounded meet-semilattices};  the least lower bound  is then denoted by $0$. 

Let $\MM$ be a meet-semilattice with top. As a poset, we have the order ideals and the order filters  of $\MM$. These   ones turn out  to be the non-empty up-sets closed under $\wedge$ and they are usually called  \emph{meet filters}, or simply \emph{filters}. Now  
 the collection $\Fi(\MM)$ 
is a closure system, and therefore it is a complete lattice where the infimum of a subset of $\Fi(\MM)$ is given by the intersection. We denote by $\Fimgg$ the  associated  closure operator. It  assigns to each $B\subseteq M$ the  \emph{filter generated by} $B$, i.e.\ the least  filter containing $B$, which can be characterized as follows. 
 For every   $a\in M$:
$$a\in \Fimg{B} \IFF a=1 \text{ or } (\exists n\in\omega)(\exists b_{0},\dots,b_{n}\in B)\,b_{0}\wedge \dots\wedge b_{n}\leq a.$$
Since $\Fi(\MM)$ is a lattice, we have its meet-irreducible elements. We denote the set of these filters by  $\Irrm(\MM)$ and call them the \emph{irreducible}  filters of $\MM$.  

A meet-semilattice $\MM$ with top   is \emph{distributive} (\cf   \cite[Sec.\ II.5]{Gr78}) when for every $a,b_{1},b_{2}\in M$ with $b_{1}\wedge b_{2}\leq a$, there exist $c_{1},c_{2}\in M$ such that $b_{1}\leq c_{1}$, $b_{2}\leq c_{2}$ and $a=c_{1}\wedge c_{2}$. 
It is well known that $\MM$ is distributive if and only if the lattice of the  filters of $\MM$ is distributive. Another characterization of distributive meet-semilattices is given in  \cite[Thm.\ 10]{Ce03b}: A meet-semilattice with top 
$\MM$ is distributive if and only if for all $F\in\Fim(\MM)$, $F\in\Irrm(\MM)$ if and only if $F^{c}\in\Id(\MM)$.

 In \cite{BeJa11} a Priestley-style duality  is presented  for two categories with objects the bounded distributive meet-semilattices and one category with morphisms the  usual algebraic homomorphisms and the other with morphisms  
 the algebraic homomorphisms that in addition preserve the existing finite joins (including the lower bounds when they exist); these morphisms  are called  there \emph{sup-homomor\-phisms}. The duality in \cite{BeJa11}
 can be slightly modified to obtain a duality for distributive meet-semilattices with top  (but not necessarily a least lower bound) as explained in a sketchy manner in \cite[Sec.\ 9]{BeJa11}. 
It should be noticed that  the duality  sketched there works only if we  modify the definitions  of Frink ideal and optimal filter of \cite{BeJa11} so that they may include the  empty set and the total  set respectively, when a lower bound does not exist. A careful presentation of this modified duality can be found in  \cite{Es13}. We expound  it briefly  here  since we will make use of it later in the paper. A useful way to define sup-homomorphism is as follows. Let $\MM_1$ and $\MM_2$ be  meet-semilattices with top. A homomorphism $h:  \MM_1\to \MM_2$ is a sup-homomorphism if for every $b \in M_1$ and every finite $A \subseteq M_1$
\begin{equation}
\text{ if } \bigcap_{a \in A}{\uparrow}a\subseteq {\uparrow}b, \text{ then } \bigcap_{a \in A}{\uparrow}h(a) \subseteq {\uparrow}h(b).
\end{equation}
Note that taking $A$ empty, $\bigcap_{a \in A}{\uparrow}a\subseteq {\uparrow}b$ holds if and only if $b$ is a lower bound of $\MM_1$; henceforth, if $h:  \MM_1\to \MM_2$ is a sup-homomorphism and $\MM_1$ has a lower bound, then $\MM_2$ should have one and $h$ should preserve the lower bounds. 

Regarding objects, the strategy followed in \cite{BeJa11}  to obtain the dual of  a  distributive meet-semilattice with top can be described as follows. 
First every  distributive meet-semilattice $\MM$ with top is embedded into a  distributive lattice $\Lee(\MM)$ with top by a meet-semilattice embedding $e$ that is also a sup-homomorphism and such that $e[M]$ is join dense in $\Lee(\MM)$. The pair $(\Lee(\MM), e)$ is called the \emph{distributive envelope of\ $\MM$}.  This lattice is in category-theoretic terms  the \emph{free distributive lattice extension of $\MM$} w.r.t.\ the forgetful functor from the category of  distributive lattices with top together with their algebraic homomorphisms to the category of distributive meet-semilattices with top and the  sup-homomorphisms. Thus  the distributive envelope of $\MM$ is (up to isomorphism) the only distributive lattice $\LL$ with top  such that there exists a meet-semilattice embedding $e: \MM \to \LL$ which is a sup-homomorphism and $e[M]$ is join-dense in $\LL$. It holds that $\MM$ has a lower bound if and only if  $\Lee(\MM)$ has a lower bound.

 The Priestley dual of $\MM$ is essentially  the
Priestley dual space of $\Lee(\MM)$ together with a dense set that allows to recover  $\MM$ inside the lattice of the clopen up-sets. 
One can take as  points of the dual space  the inverse images of the prime filters of $\Lee(\MM)$ by the embedding $e: \MM \to \Lee(\MM)$, together with $M$ when $\MM$ has no lower bound. Under the identification of $M$ with $e[M]$, the points of the space, which are  called \emph{optimal filters},  are then the intersection with $M$ of the prime filters of $\Lee(\MM)$, with the addition of $M$ if $\MM$ has no lower bound.  

The notion of optimal filter can be defined for  every meet-semilattice $\MM$   using   the notion of Frink ideal.  Moreover,  the optimal filters  can be used to give one of the  particular constructions of the distributive envelope of a distributive meet-semilattice. 
Let $\MM$ be a meet-semilattice with top. A filter $F\in\Fim(\MM)$ is  \emph{optimal} when there is $I\in\IdF(\MM)$ such that 
 $F$ is a maximal element of $\{G\in \Fim(\MM):G\cap I=\emptyset\}$ and 
$I$ is a maximal element of $\{J\in\IdF(\MM):F\cap J=\emptyset\}$.
We denote by $\Optm(\MM)$ the set of all optimal  filters of $\MM$. 
It is easy to check that the irreducible  filters are optimal. 
Moreover,  $M$ is an optimal  filter if and only if there is no bottom element in $\MM$. If $\MM$ is distributive, then  
$F\in\Optm(\MM)$ if and only if $F^{c}\in\IdF(\MM)$. This shows that the definition of optimal filter given here and the definition in \cite{BeJa11}  are coextensive. For distributive meet-semilattices with top the irreducible filters and the  optimal filters  are characterized in the next proposition using the concept of $\wedge$-prime set.  A set $X \subseteq M$ 
 of a meet-semilattice $\MM$ with top  is said to be \emph{$\wedge$-prime}  (or simply \emph{prime}) if it is   proper  and  for all non-empty finite $U\subseteq  M$, if 
$\bigwedge U \in I$, then $U \cap I \not = \emptyset$. We apply the concept $\wedge$-prime to order ideals and Frink ideals. Note that if $\emptyset$ is a Frink ideal, then it is $\wedge$-prime.
We address the reader to sections 2.3, 3.2.2 and Appendix A in \cite{Es13} for the proof of the next proposition.

\begin{prop} \label{cor:Optm-IdF}\label{cor:Irrm-Id}
Let $\MM$ be a distributive meet-semilattice with top. For every $F\subseteq M$, 
\begin{enumerate}
\item  $F\in \Irrm(\MM)$ if and only if $F^{c}$ is a $\wedge$-prime order ideal,
\item  $F\in\Optm(\MM)$ if and only if $F^{c}$ is a $\wedge$-prime Frink ideal.
\end{enumerate}
\end{prop}

For a given distributive meet-semilattice $\MM=\langle M,\wedge,1\rangle$ with top,  the \textit{representation map}  $\sigma:M\longrightarrow \mathcal{P}^{\uparrow}(\Optm(\MM))$ is defined by setting for every $a\in M$ 
\begin{equation*}
\sigma(a):=\{P\in \Optm(\MM) :a\in P\}.
\end{equation*}
The set $\sigma[M] := \{\sigma(a): a \in M\}$ is closed under the binary operation of intersection of sets and $\sigma(1) = \Optm(\MM)$. Therefore,  $\sigma[\MM]:= \langle \sigma[M], \cap, \sigma(1)\rangle$ is a meet-semilattice with top and $\sigma$  is an isomorphism between $\MM$ and $\sigma[\MM]$. Note that $\MM$ has a bottom   element if and only if $\emptyset \in \sigma[M]$.

A particular construction of the  \emph{distributive envelope} of $\MM$ is given by the distributive lattice $\Lee(\MM)$ that we obtain by closing  $\sigma[M]$ under the binary operation of set-theoretic union. The embedding from  $\MM$ to $\Lee(\MM)$ that shows that this lattice is the distributive envelope of $\MM$ is $\sigma$. A very useful property of $\sigma$ that follows easily from the fact that $\sigma$ is a sup-homomorphism is the following. For every $a, b_0, \ldots, b_n\in M$,
\begin{equation}
\label{eq:distr-envelop-1}
\bigcap_{i \leq n}{\uparrow}b_i \subseteq {\uparrow}a \  \text{  iff } \ \sigma(a) \subseteq \bigcup_{i \leq n} \sigma(b_i). 
\end{equation} 

A particular construction of the  \emph{distributive envelope} of $\MM$ is given by the distributive lattice $\Lee(\MM)$ that we obtain by closing  $\sigma[M]$ under the binary operation of set-theoretic union and the embedding $\sigma$ from  $\MM$ to $\Lee(\MM)$. 
%
In this particular  construction of the distributive envelope,  the relation between the filters of $\MM$ and the filters of $\Lee(\MM)$ is given by the map  $\Fimg{\sigma[.]}$.  When it  is applied to  $\Fim(\MM)$ establishes an isomorphism between  $\Fim(\MM)$ and the lattice of the filters of   $\Lee(\MM)$ whose inverse is easily seen to be  given by the map    $\sigma^{-1}[.]$. Moreover, the lattice of the Frink ideals of $\MM$ and the lattice of the ideals of $\LL(\MM)$ together with the empty set, when $\MM$ has no lower bound, are also isomorphic. The isomorphism is  given by the map that sends a Frink ideal $I$ of $\MM$ to the ideal generated by $\sigma[I]$, when $I$ is non-empty, and to the emptyset otherwise. The inverse of this  isomorphism  is also given by the map $\sigma^{-1}[.]$. It follows that 
 $\Fimg{\sigma[.]}$  establishes  an isomorphism  between $\Optm(\MM) \setminus \{M\}$ and the set $\Pr(\Lee(\MM))$ of the prime filters of $\Lee(\MM)$, when $\MM$ has no bottom element, and between $\Optm(\MM)$  and $\Pr(\Lee(\MM))$, when it has. 
 
 The dual objects of distributive meet-semilattices with top  are called $\star$-generalized Priestley spaces in  \cite{BeJa11}. 
We delete the star in this paper. 

A \emph{generalized Priestley space} is a tuple
$\spa{X}=\langle X,\tau,\leq,\Xb \rangle$ such that
\begin{itemize}
	\labitem{(DS1)}{itm:PrDS1} $\langle X,\tau,\leq\rangle$ is a Priestley space,
	\labitem{(DS2)}{itm:PrDS2} $\Xb $ is a dense subset of $X$,
	\labitem{(DS3)}{itm:PrDS3}$\Xb =\{x\in X: \{U\in X^{*}:x\notin U\} \text{ is non-empty and up-directed}\}$,
	\labitem{(DS4)}{itm:PrDS4}  for all $x,y\in X$, $x\leq y$ iff  $y\in U$ for all  $U\in X^{*}$  such that $x\in U$,
\end{itemize}
where $X^{*}\coloneqq\{U\in\clopu(X):\max(U^{c})\subseteq \Xb \}$. The elements of $X^{*}$ are called the  \emph{$\Xb$-admissible} clopen up-sets of $\spa{X}$. The collection $X^{*}$ is closed under the binary operation of  intersection, 
and  the structure $\spa{X}^{*}\coloneqq\langle X^{*},\cap,X\rangle$ is a distributive meet-semilattice with top. It follows from the denseness of $\Xb$ that $\spa{X}^{*}$ has a lower bound if and only if $\emptyset \in X^{*}$. The meet-semilattice $\spa{X}^{*}$ is \textit{the dual} of $\spa{X}$.

The following facts hold for every generalized Priestley space $\spa{X}=\langle X,\tau,\leq,\Xb \rangle$:\label{prop:closureadmisibles-distr-envelop}

\begin{enumerate}
\item The closure $(X^{*})^{\cup}$ of $X^{*}$ under the  operation $\cup$ of set-theoretic union is the set of all non-empty clopen up-sets, if $\emptyset \not \in X^{*}$, and it is the set of all clopen up-sets, if  $\emptyset  \in X^{*}$.
\item The distributive envelope of  $\spa{X}^{*}$ is (up to isomorphism) the lattice $\Lee(\spa{X}^{*}) = \langle (X^{*})^{\cup}, \cap, \cup, X\rangle$, with embedding the identity map.
\item For every proper optimal filter $F$ of $\spa{X}^{*}$, the filter generated by $F$  in $\Lee(\spa{X}^{*})$ is a prime filter and for every prime filter $G$ of $\Lee(\spa{X}^{*})$, $G \cap X^{*}$ is an optimal filter of $\spa{X}^{*}$.

\end{enumerate}


The dual of    a distributive meet-semilattice $\MM=\langle M,\wedge,1\rangle$ with top  is defined  in  \cite{BeJa11} as follows. Let 
 $\tau_{\MM}$ be the topology on $\Optm(\MM)$ determined by the subbasis  $\{\sigma(a):a\in M\}\cup\{\sigma(b)^{c}:b\in M\}$. Then   $\langle \Optm(\MM),\tau_{\MM},\subseteq\rangle$ is a Priestley space. 
And the structure $\sOptm(\MM):=\langle \Optm(\MM),\tau_{\MM},\subseteq,\Irrm(\MM)\rangle$ turns out to be a generalized Priestley space such that 
$\Optm(\MM)^{*} = \sigma[M]$. This  implies that $\MM$
 is isomorphic to $(\sOptm(\MM))^{*}$  by means of the map $\sigma$. 
 Furthermore, the following facts hold for every  distributive meet-semilattice with top $\MM=\langle M,\wedge,1\rangle$ (for a proof see  \cite[Sec.\ 5]{BeJa11})\label{prop:gps1}:

\begin{itemize}
	\item[(1)] Every non-empty clopen-upset of $\langle \Optm(\MM),\tau_{\MM},\subseteq\rangle$ is the  union of a finite and non-empty subset of $\sigma[M]$.
	\item[(2)] For every clopen up-set $U\in\clopu(\Optm(\MM))$, $U$ is  an $\Irrm(\MM)$-admissible clopen up-set if and only if there exists $a\in M$ such that $U=\sigma(a)$. 
\end{itemize}

Any \DSp space $\spa{X}=\langle X,\tau,\leq,\Xb \rangle$ is order homeomorphic to 
 $\langle \Optm(\spa{X}^{*}),\tau_{\spa{X}^{*}},\subseteq, \Irrm(\spa{X}^{*})\rangle$ by means of the map 
 $\eps:X\longrightarrow \Optm(\spa{X}^{*})$, given by:
 \begin{equation*}
 \eps(x)\coloneqq\{U\in X^{*}:x\in U\}.
 \end{equation*} 
 Moreover, $\eps[\Xb]=\Irrm(\spa{X}^{*})$.

Regarding morphisms, we just need to recall how the dual of an algebraic homomorphism, which is a relation, is defined in \cite{BeJa11}.  

For  generalized Priestley spaces $\spa{X}_{1}$ and $\spa{X}_{2}$, a relation $R\subseteq X_{1}\times X_{2}$ is a  \emph{generalized Priestley morphism}
  (\cite[Def.\ 6.2]{BeJa11}) when:
 \begin{itemize}
	\labitem{(DSR1)}{itm:PrDSR1} $\Box_{R}(U)\in X_{1}^{*}$ for all $U\in X_{2}^{*}$,
	\labitem{(DSR2)}{itm:PrDSR2} if $(x,y)\notin R$, then there exists $U\in X_{2}^{*}$ such that $y\notin U$ and $R(x)\subseteq U$.
\end{itemize}
A generalized Priestley morphism $R$ is \emph{functional} when in addition: 
 \begin{itemize}
	\labitem{(DSR3)}{itm:PrDSR3} for every $x \in X_1$ there exists $y \in X_2$ such that $R(x) = {\uparrow}y$.
	\end{itemize}

For a given   \DSp morphism $R\subseteq X_{1}\times X_{2}$,  the map $\Box_{R}:\pws(X_{2})\longrightarrow \pws(X_{1})$ is an algebraic homomorphism  between the distributive meet-semilattices with top   $\spa{X}_{2}^{*}$ and $\spa{X}_{1}^{*}$. If $R$ is functional, then $\Box_{R}$ is a sup-homomorphism.
Moreover, for all $x\in X_{1}$ and all $y\in X_{2}$, $(x,y)\in R$ if and only if $(\eps(x),\eps(y))\in R_{\Box_{R}}$. 
For every homomorphism $h:\MM_{1}\longrightarrow \MM_{2}$ between distributive meet-semilattices with top, the relation $R_{h}\subseteq \Optm(\MM_{2})\times \Optm(\MM_{1})$ defined by:
\begin{equation*}
(P,Q)\in R_{h}  \IFF  h^{-1}[P]\subseteq Q,
\end{equation*} 
is a generalized Priestley morphism between the dual generalized Priestley spaces $\sOptm(\MM_{2})$ and $\sOptm(\MM_{1})$; $ R_{h}$  is functional when $h$ is a sup-homo\-mor\-phism.
Moreover,  $\sigma_{2}(h(a))=\Box_{R_{h}}(\sigma_{1}(a))$ for all $a\in M_{1}$. 

 The notion of generalized Priestley morphism has a slight drawback: 
the usual composition of relations  does not always produce a generalized Priestley morphism when applied to two of them; hence it can not be taken as a  category-theoretic operation of composition.
Instead, we have that for any generalized Priestley spaces 
 $\spa{X}_{1}$, $\spa{X}_{2}$ and $\spa{X}_{3}$ and any generalized Priestley morphisms $R\subseteq X_{1}\times X_{2}$ and $S\subseteq X_{2}\times X_{3}$, the operation of composition  of  $R$ and $S$ that we need in order to obtain a category is  the relation  $S\star R \subseteq X_{1} \times X_{3}$ defined by:
\begin{equation*}
(x,z)\in (S\star R)\IFF \forall U\in X_{3}^{*}\big(x \in (\Box_{R}\circ\Box_{S}) (U) \Rightarrow  z\in U   \big).
\end{equation*}

\subsection{Congruential logics}

Congruential logic is a concept studied in  abstract algebraic logic\footnote{%
We follow here the terminology used in \cite{GeJaPa10}.  
Congruential logics were previously called \emph{strongly selfextensional} \cite{FJa09} and \emph{fully selfextensional} \cite{JaPa06}. The last terminology is currently widely used in abstract algebraic logic.}. To be able to present it we need to go throughout some of the basic concepts of this field. 
We  follow the survey \cite{FoJaPi03}. 

 A  logical language, or algebraic type, $\mathscr{L}$ is a set of connectives, a.k.a.\  operation symbols, possibly of arity $0$.  
A \emph{logic}  
   is a pair $\mathcal{S}\coloneqq\langle \Fm_{\mathscr{L}_{\mathcal{S}}},\vdashs\rangle$, where $\mathscr{L}_{\mathcal{S}}$ is a logical language, $\Fm_{\mathscr{L}_{\mathcal{S}}}$ is   the absolutely free algebra of terms of type $\mathscr{L}_{\mathcal{S}}$ over a countably infinite set of generators, called the variables,   and ${\vdashs}\subseteq {\pws(Fm_{\mathscr{L}_{\mathcal{S}}})\times Fm_{\mathscr{L}_{\mathcal{S}}}}$ is a substitution-invariant consequence relation
     on the universe $Fm_{\mathscr{L}_{\mathcal{S}}}$ of $\Fm_{\mathscr{L}_{\mathcal{S}}}$, called the set of formulas (or terms) of $\Sm$.\footnote{
This means that the following conditions are satisfied: 
\begin{itemize}
	\item[(C1)] if $\gamma\in\Gamma$, then $\Gamma \vdashs \gamma$,
	\item[(C2)] if $\Delta\vdashs\gamma$ for all $\gamma\in\Gamma$ and $\Gamma\vdashs\delta$, then $\Delta\vdashs\delta$, 
	\item[(C3)] if $\Gamma\vdashs\delta$, then $\sigma[\Gamma]\vdashs \sigma(\delta)$ for all substitutions $\sigma\in\Hom(\Fm_{\mathscr{L}},\Fm_{\mathscr{L}})$.
\end{itemize}
}
A logic $\Sm$ is \emph{finitary} when for all $\Gamma\cup\{\delta\}\subseteq Fm_{\mathscr{L}}$, 
 if $\Gamma\vdashs\delta$, 
 then there exists a finite $\Gamma_{0}\subseteq\Gamma$ such that $\Gamma_{0}\vdashs\delta$. 
A logic $\Sm$ \emph{has theorems} when there is at least one formula $\delta\in Fm$ such that $\emptyset\vdashs\delta$. 

Let $\Sm$ be a logic and $\Al$  an algebra of type $\mathscr{L}_{\Sm}$, that for short we say it is an $\mathscr{L}_{\Sm}$-algebra.

\begin{defn} \label{defn:FiS} 
 A set $F\subseteq A$ is an \emph{$\Sm$-filter}   of $\Al$ when for every $h\in\Hom(\Fm_{\mathscr{L}_{\Sm}},\Al)$ and  every set of formulas 
$\Gamma\cup\{\delta\}\subseteq Fm_{\mathscr{L}_{\Sm}}$:   
\begin{equation*} 
\text{ if } \Gamma\vdash_{\Sm} \delta \text{ and }h[\Gamma] \subseteq F\text{, then }h(\delta)\in F.
\end{equation*}
\end{defn}
We denote by $\FiSA$  the collection of all $\Sm$-filters of $\Al$. This collection  is always a closure system, \ie it is closed under arbitrary intersections, and therefore is a complete lattice under the order of  inclusion.  If $\Sm$ is a finitary logic, it is then a finitary closure system. We call  \emph{irreducible} $\Sm$-filters the
meet-irreducible elements of the lattice of $\Sm$-filters  and  denote by  $\IrrSA$ the collection of all of them.  
 The closure operator associated with $\FiSA$ will be denoted by $\Consa$.  
 Thus, for every  $U\subseteq A$, $\Consa(U)$ is the least $\Sm$-filter of $\Al$ that contains $U$.  We abbreviate   $\Consa(\{a\})$ by $\Consa(a)$. When  $\FiSA$ is finitary, $\Consa$  is a finitary closure operator. 
 
 A basic property of $\Sm$-filters whose use is ubiquitous in abstract algebraic logic is the following. Let  $\Al_1, \Al_2$ be $\mathscr{L}_{\Sm}$-algebras. For every $h \in \Hom(\Al_1, \Al_2)$ and every $\Sm$-filter $F$ of $\Al_2$, $h^{-1}[F]$ is an $\Sm$-filter of $\Al_1$. 

A logic $\Sm$ is \emph{filter-distributive} when for every $\mathscr{L}_{\Sm}$-algebra $\Al$, the lattice of $\Sm$-filters of $\Al$ is distributive.

 %
The closure operator $\Consa$ allows us to define the \emph{specialization quasiorder} $\leqsa$ on $A$ by saying that for all $a,b\in A$: 
\begin{equation*}
a\leqsa b\IFF\Consa(b)\subseteq\Consa(a). 
\end{equation*}
We denote by $\equivsa$ the equivalence relation associated with $\leqsa$, \ie ${\equivsa}\coloneqq {\leqsa\cap \geqsa}$. 
We use this relation to introduce the following concepts: 

\begin{defn}  \label{defn:congruential} 
A logic $\Sm$ is  \emph{congruential}
when for every $\mathscr{L}_{\Sm}$-algebra $\Al$,  the relation $\equivsa$ is a congruence of $\Al$. 
\end{defn}
This definition is equivalent to the more usual one of the concept as shown   in \cite[Prop.\ 2.42]{FJa09}. The next definition is also equivalent to the more usual one given in for example \cite{FJa09}.

\begin{defn} 
An $\mathscr{L}_{\Sm}$-algebra $\Al$ is an \emph{$\Sm$-algebra} when for every congruence $\theta$ of $\Al$, if $ \theta\subseteq {\equivsa}$, then $\theta$ is the identity.
We denote by $\Alg\Sm$ the collection of all $\Sm$-algebras. This class of algebras is known as \emph{the algebraic counterpart} of $\Sm$.
\end{defn}

We remark that  the trivial algebras, namely the algebras with a single element, are $\Sm$-algebras.
 
Many  well-known logics, including classical and intuitionistic propositional logics,  are congruential. 
The next theorem provides a useful characterization  of congruentiality. 

\begin{thm}[{\cite[Theorem 2.2]{GeJaPa10}}] \label{thm:GeJaPa10}A logic $\Sm$ is congruential if and only if for every $\mathscr{L}_{\Sm}$-algebra $\Al$,  
$\Al\in\Alg\Sm$ if and only if  $\langle A ,\leqsa\rangle$ is a poset.
\end{thm}

Notice that when $\Sm$ is congruential, for every $\Sm$-algebra $\Al$ we have the collection $\Fi(\Al)$ of the order filters of the poset $\langle A, \leqsa \rangle$  and the collection $\Id(\Al)$ of its order ideals.  All $\Sm$-filters of $\Al$ are up-sets with respect to $\leqsa$, but not necessarily order filters of $\langle A, \leqsa \rangle$ because they may not be down-directed. Note also that  for every $a\in A$, $\Consa(a)= {\uparrow}a$.  
If a congruential logic $\Sm$ has theorems, then all $\Sm$-filters of $\Al$ are non-empty and the poset $\langle A,\leqsa\rangle$ has a top element, that we denote by $1^{\Al}$.  Note that then $\{1^{\Al}\}$ is the least $\Sm$-filter of $\Al$.
Furthermore, from the previous theorem we infer that for every  congruential logic $\Sm$, 
\begin{equation*}
\Alg\Sm=\{{\Al}: \Al \text{ is an $\mathscr{L}_{\Sm}$-algebra and } {{\equivsa}\text{ is the identity}}\}.
\end{equation*}

 We recall now the definition of \emph{$\Sm$-ideal}  given in \cite{GeJaPa10}.
\begin{defn}\label{defn:Sideal}
A subset $I\subseteq A$ is an \emph{$\Sm$-ideal} of an $\mathscr{L}_{\Sm}$-algebra $\Al$ provided that for any finite $I'\subseteq I$ and any $a\in A$, 
if $\bigcap\{\Consa(b):b\in I'\}\subseteq \Consa(a)$, then $a\in I$. 
\end{defn}
We denote by $\IdSA$ the collection of all $\Sm$-ideals of $\Al$. Notice that when $\Al$ is an $\Sm$-algebra of a congruential logic, then the $\Sm$-ideals of $\Al$ are exactly the Frink ideals of the poset  $\langle A ,\leqsa\rangle$. It holds  that 
 $\emptyset\in\IdSA$ if and only if the poset $\langle A,\leqsa\rangle$ has no bottom element. 

We are interested in a certain type of $\Sm$-ideals that will help to define the notion of  optimal $\Sm$-filter we need. 

\begin{defn}\label{defn:sSideal}
An $\Sm$-ideal $I$ of $\Al$ is \emph{strong} when for every finite $I'\subseteq I$ and every non-empty and finite $B\subseteq A$, 
if $\bigcap\{\Consa(b):b\in I'\}\subseteq\Consa(B)$, then   $\Consa(B)\cap I\neq\emptyset$. 
\end{defn}

Note that the definition implies that $A$ is a strong $\Sm$-ideal. The next lemma gives a characterization of the strong $\Sm$-ideals for congruential logics.

\begin{lemma}
\label{lem:useful-proof-strongideal}
Let $\Sm$ be a congruential logic and $\Al$ an $\Sm$-algebra. A set $I \subseteq A$ is a strong $\Sm$-ideal if and only if it is a down-set w.r.t.\ $\leqsa$ and satisfies the condition of Definition \ref{defn:sSideal}. 
\end{lemma}
\begin{proof}
We only need to proof the implication from right to left. Assume that $I$ is a down-set w.r.t. $\leqsa$ and satisfies the condition of Definition \ref{defn:sSideal}. We have to see that it is an $\Sm$-ideal. Suppose that $I' \subseteq I$ is finite and $\bigcap\{\Consa(b):b\in I'\}\subseteq \Consa(a)$ for $a \in A$. By the condition on Definition \ref{defn:sSideal}, $I \cap \Consa(a) \not = \emptyset$. Thus there is $c \in I$ such that $a \leqsa c$. Hence $a \in I$.
\end{proof}

\begin{rem}
If $\Sm$ has theorems, then in the definition of strong $\Sm$-ideal we can delete the condition that $B$ is non-empty and we obtain an equivalent definition. If we take $B$ possibly empty in the case that $\Sm$ does not have theorems, then $A$ may not be a strong $\Sm$-ideal in case that there are $a, b \in A$ such that $\Consa(a) \cap \Consa(b) = \emptyset$.
\end{rem}

We denote by $\IdsS(\Al)$ the collection of all strong $\Sm$-ideals. 
It is easy to check that when $\Al$ is an $\Sm$-algebra all order ideals of $\langle A,\leqsa\rangle$ are strong $\Sm$-ideals, in particular for every $a \in A$, ${\downarrow}_{\leqsa}a$ is a 
strong $\Sm$-ideal.

The next notion helps to characterize when the emptyset is a strong $\Sm$-ideal.
%
Let $\Sm$ be a congruential logic and $\Al$ an $\Sm$-algebra. We say that a non-empty finite set $U\subseteq  A$  is a \emph{bottom-family} of $\Al$ when $\Consa(U)=A$. 
%
Note that if $\leqsa$ has a bottom element, then its singleton is a bottom-family. In particular, in the trivial $\Sm$-algebras the domain of the algebra  is a bottom-family. Moreover, it is easy to see that if it  has a bottom family,  then it has a bottom-family of incomparable elements with respect to $\leqsa$.
  
  A straightforward argument shows that  $\emptyset\in\IdsS(\Al)$ if and only if $\Al$ has no bottom-family.

We are now in a position to define  the notion of optimal $\Sm$-filter. Note the similarity with the definition of optimal filter of a meet-semilattice. 

\begin{defn}\label{defn:optimalfilter}
An $\Sm$-filter $F\in\FiSA$ is \emph{optimal} when there is a strong $\Sm$-ideal $I\in\IdsS(\Al)$ such that 
 $F$ is a maximal element of the collection $\{G\in\FiSA:G\cap I=\emptyset\}$ and 
$I$ is a maximal element of the collection $\{J\in\IdsS(\Al):F\cap  J=\emptyset\}$. 
\end{defn}
We denote by $\OptSA$ the collection of all optimal $\Sm$-filters of $\Al$. 

\begin{rem}
\label{rem:optimal-bottom-family}
From the definition it follows that $\emptyset\in\IdsS(\Al)$  if and only if $A\in\OptSA$. Therefore, $A\in\OptSA$ if and only if $\Al$ has no bottom-family. Hence, in the trivial $\Sm$-algebras there is no non-empty optimal $\Sm$-filter.
\end{rem}

\begin{rem}
If $\Sm$ does not have theorems, then since $\emptyset$ is an $\Sm$-filter and $A$ a strong $\Sm$-ideal, $\emptyset$ is an optimal $\Sm$-filter.
\end{rem}

For every finitary congruential logic, we have the following two separation lemmas that we gather in one proposition. They relay on Zorn's lemma  and  the fact that for a finitary logic $\Sm$, $\FiSA$ is closed under unions of non-empty chains and for every logic $\Sm$ both $\IdsS(\Al)$   and $\Id(\Al)$ are closed under unions of non-empty chains.

\begin{prop}[Optimal and irreducible $\Sm$-filter lemmas] \label{prop:ofls}
Let $\Sm$ be a finitary congruential logic, $\Al$  an $\Sm$-algebra and $F\in\FiSA$,
\begin{enumerate}
\item  if $I\in\IdsSA$ is such that $F\cap I=\emptyset$, then  there exists $Q\in \OptSA$ such that $F\subseteq Q$ and $Q\cap I=\emptyset$. 
\item if $I\in\Id(\Al)$ is such that $F\cap I=\emptyset$, then there exists
 $Q\in \IrrSA$ such that $F\subseteq Q$ and $Q\cap I=\emptyset$. 
\end{enumerate}
\end{prop}
\begin{proof}
(1) First of all note  that $\IdsSA$ is closed under unions of non-empty chains.  Let $F\in\FiSA$ and $I\in\IdsSA$ be such that  $F\cap I=\emptyset$ and consider the set 
$\FF := \{G \in \FiSA: F \subseteq G \text{ and } G \cap I = \emptyset\}$, which is non-empty and closed under unions of non-empty chains because the closure operator of $\Sm$-filter generation is finitary. By Zorn's lemma there exists  a maximal element $Q$ of $\FF$. Consider now the set 
$\II:= \{H \in  \IdsSA: I \subseteq H \text{ and } H \cap Q = \emptyset\}$. Then $I \in \II$ and $\II$ is closed 
under unions of non-empty chains. Let, by Zorn's lemma, $H \in \II$ be maximal. Then $H$ is $Q$-maximal and it is easy to see that $Q$ is $H$-maximal. Therefore, $Q$ is optimal.

(2)  Let  $F\in\FiSA$ and  $I\in\Id(\Al)$ be such that $F\cap I=\emptyset$. Consider the set $\FF':= \{G \in \FiSA: F \subseteq G \text{ and } G \cap I = \emptyset\}$. Then $F \in \FF'$ and $\FF'$ is closed under unions of non-empty chains. By Zorn's lemma we take a
 maximal elment $Q$ of $\FF'$.  Clearly, since $I \not = \emptyset$, then $Q$ is proper. To show that $Q$  is irreducible suppose that $F_1, F_2 \in \FiSA$ are such that $F_1 \cap F_2 = Q$ and in search of a contradiction that $F_1 \not = Q$ and $F_2 \not = Q$.
 Let then $a \in F_1 \setminus Q$ and $b \in F_2 \setminus Q$. By the maximality of $Q$ in $\FF'$, there are $a' \in \Consa(Q \cup\{a\}) \cap I$ and $b' \in \Consa(Q \cup\{b\}) \cap I$. Since $I$ is up-directed let $c \in I$ such that $a', b' \leqsa c$. Then $c \in \Consa(Q \cup\{a\}) \cap \Consa(Q \cup\{b\}) \subseteq F_1 \cap F_2 = Q$, contradicting the fact that $Q \cap I = \emptyset$.
 \end{proof}

When we restrict ourselves to finitary filter-distributive  logics,  we have good characterizations of the optimal and of the irreducible $\Sm$-filters. 


\begin{thm}\label{thm:OptS-IdsS}\label{thm:IrrS-Id}
Let $\Sm$ be a  filter-distributive, finitary, and congruential logic  and let $\Al$ be an $\Sm$-algebra. Then for every $F\in \FiSA$, 
\begin{enumerate}
\item 
$F\in\OptSA$ if and only if $F^{c}\in\IdsSA$,
\item  $F\in\IrrSA$ if and only if $F^{c}\in\Id(\Al)$.
\end{enumerate}
\end{thm}
\begin{proof}
(1) If $F^{c} \in \IdsSA$, then   obvioulsy $F$ is  $F^{c}$-maximal  and $F^{c}$ is $F$-maximal; hence $F$ is optimal. To prove the converse, suppose that $F$ is optimal. Let then $I \in \IdsSA$ be such that $F$ is $I$-maximal and $I$ is $F$-maximal. We prove that $F^{c} \in \IdsSA$. We reason by cases. If $F = A$, then $F^{c} = \emptyset$ and therefore $I = \emptyset = F^{c}$; hence $F^{c} \in \IdsSA$. Suppose then that $F \not = A$.  
Let $B$ be a finite and non-empty subset of $A$ and $H$ a finite subset of $F^{c}$. Suppose that $\bigcap\{\Consa(b): b \in H\} \subseteq \Consa(B)$.  If $H = \emptyset$, then $\Consa(B) = A$, and since $F^{c} \not = \emptyset$, $\Consa(B) \cap F^{c} \not = \emptyset$. If $H \not = \emptyset$,  since $F$ is $I$-maximal, for every $b \in H$, let $a_b \in \Consa(F, b) \cap I$, and by the finitarity of the closure operator $\Consa$, let $F_b \subseteq F$ be  finite and such that $a_b \in \Consa(F_b, b)$. Consider then the  finite set $G = \bigcup_{b \in H}F_b$. Then for every $b \in H$, $a_b \in \Consa(G, b) = \Consa(G) \sqcup \Consa(b)$.
Using the distributivity of $\FiSA$ we have
\begin{align*}
\bigcap_{b \in H} \Consa(a_b) & \subseteq \bigcap_{b \in H} \Consa(G) \sqcup \Consa(b)\\
& = \Consa(G) \sqcup \bigcap_{b \in H} \Consa(b)\\
&  \subseteq \Consa(G) \sqcup \Consa(B).
\end{align*}
 From the fact that  $I$ is a strong $\Sm$-ideal that includes $\{a_b: b \in H\}$, it follows that $\Consa(G \cup B) \cap I \not = \emptyset$. Suppose towards a contradiction that $\Consa(B) \cap F^{c} = \emptyset$. Then $\Consa(B) \subseteq F$ and therefore $\Consa(G \cup B) \subseteq F$. Since $G \subseteq F$ we get that $F \cap I \not = \emptyset$, a contradiction.

(2) Suppose that $F\in\IrrSA$. By assumption $F$ is a proper  up-set w.r.t.\ $\leqsa$ and hence $F^{c}$ is a non-empty down-set. To show that it is up-directed, let $a, b \in F^{c}$, so that $\Consa(a), \Consa(b) \not \subseteq F$. Since $\FiSA$ is distributive and $F$ is meet-irreducible, $F$ is meet-prime. Therefore,  $\Consa(a) \cap \Consa(b) \not \subseteq F$. Let $c \in \Consa(a) \cap \Consa(b)$ be such that $c \not \in F$. Then 
$a, b \leqsa c$ and $c \in F^{c}$. Thus $F^{c}$ is up-directed. To prove the converse, assume that $F^{c} \in  \Id(\Al)$. We show that $F$ is a meet-prime element of $\FiSA$. Since $F^{c}$ is non-empty, $F$ is proper. Suppose that $F_1, F_2 \in \FiSA$ are such that $F_1 \cap F_2 \subseteq F$,  $F_1 \not \subseteq F$ and $F_2 \not \subseteq F$. Let then $b_1 \in F_1 \setminus F$ and $b_2 \in F_2 \setminus F$. Since $F^{c}$ is up-directed, let $c \in F^{c}$ be such that $b_1, b_2 \leqsa c$. Then $c \in  \Consa(b_1) \cap \Consa(b_2) \subseteq F_1 \cap F_2$ and hence $c \in F$, a contradiction. 
\end{proof}
As a consequence, we have that when $\Sm$ is a filter-distributive, finitary,  and congruential logic   all the optimal $\Sm$-filters are proper if and only if $\Al$ has a bottom-family.

\begin{rem}
We can separate the proof if we like in two cases depending on whether  $\Sm$ has theorems or not. In the second case, condition (1) holds for $\emptyset$ even if $\Sm$ is not filter-distributive, so the real proof handles the non-empty case. But even in this case, to carry on the proof we need the given definition of strong $\Sm$-ideal (that requires $B$ to be non-empty).
\end{rem}

We conclude this section by introducing a  concept of primness  that will be used later on. 


\begin{defn}\label{defn:Sprime}
A subset $X\subseteq A$ is  \emph{$\Sm$-prime} when it is a proper subset such that for all non-empty and  finite $B\subseteq A$, if $\Consa(B)\cap X\neq \emptyset$, then $B\cap X\neq \emptyset$. 
\end{defn}

\begin{lemma}[{\cite[Lemma 4.4.11]{Es13}}] \label{lemma:FiS-Sprime} Let $\Sm$ be a finitary logic, $\Al$ an $\Sm$-algebra and   $X\subseteq A$ 
non-empty. Then $X\in\FiSA$ if and only if $X^{c}$ is $\Sm$-prime.
\end{lemma}

The requirement that $X$ is non-empty is important when $\Sm$ has no theorems. Otherwise, the lemma would fail for $\emptyset$. The definition of $\Sm$-prime implies that $A$ is not $\Sm$-prime. 

\section{A categorical duality}\label{sec:catduality}

In this section we obtain a dual category  for the category of $\Sm$-algebras and their homomorphisms for finitary, congruential and filter-distributive logics with theorems. Nevertheless we state the results that require less assumptions to hold with the minimal ones they need.

\subsection{Representation theorem for $\Sm$-algebras}\label{sec:repr}

From now on we focus on congruential logics. Let us fix an arbitrary congruential logic  $\Sm$ and  an $\Sm$-algebra $\Al$. We are interested in  closure bases for the closure system  $\FiSA$. 
Recall that $\FF\subseteq\FiSA$ is a \emph{closure base}
for $\FiSA$ provided that every $\Sm$-filter is an intersection of elements in $\FF $. From Proposition \ref{prop:ofls} follows that if $\Sm$ is finitary and congruential, then $\OptSA$ is a closure base for $\FiSA$. These bases will play a central role in our duality.


\label{defn:map:varphi}We fix an arbitrary closure base $\FF  $ for $\FiSA$.   
The \textit{representation map} $\varphif:A\longrightarrow \ups(\FF )$ is defined by setting for every $a\in A$:
\begin{equation*}
\varphif(a):=\{P\in \FF :a\in P\}.
\end{equation*}
Since $\Sm$ is congruential and $\Al$ is an $\Sm$-algebra, it is easily seen that the  map $\varphif:A\longrightarrow \ups(\FF )$ is injective. We lift  $\varphif$ to a map $\varphifm: \mathcal{P}(A) \to  \ups(\FF )$ by defining for every $B\subseteq A$: 
$$\varphifm(B)\coloneqq \bigcap\{\varphif(b):b\in B\}.$$
Therefore $\varphifm(B) =\{P\in \FF :B\subseteq P\}$. 
Notice that for any $B,B'\subseteq A$, we have   $\varphifm(B)\cap\varphifm(B')=\varphifm(B\cup B')$. Notice also that $\varphifm(\emptyset) = \FF$ and that if $\Sm$ has theorems, then $\varphif(1^{\Al}) = \FF = \varphifm(\emptyset)$.
It is not difficult to see that for any non-empty and finite $B\subseteq A$, $\varphifm(B)\neq\{A\}$. In particular, $\varphif(b)  \neq \{A\}$, for every $b \in A$.

%

We denote by $\varphifA$ the algebra whose universe is $\varphif[A]$ and  that
for every $n$-ary connective $f$ of the language $\mathscr{L}_{\Sm}$, 
the operation $f^{\varphifA}$ on $\varphif[A]$ is defined by setting for all elements $a_{1},\dots,a_{n}\in A$, 
\begin{equation*}
{f^{\varphifA} (\varphif(a_{1}),\dots,\varphif(a_{n}))}\coloneqq{\varphif(f^{\Al}(a_{1},\dots,a_{n}))}.
\end{equation*}
In particular, if $c$ is a $0$-ary connective, $c^{\varphifA} = \varphif(c^{\Al})$.
These operations are well defined since the map $\varphif$ is injective.
Thus, $\varphifA $ is  a well defined algebra and moreover $\varphif\in\Hom(\Al,\varphifA )$.

\begin{thm} \label{thm:srepresentation} 
The map $\varphif:A\longrightarrow \mathcal{P}^{\uparrow}(\FF )$ 
is an isomorphism  between  $\Al$ and $\varphifA $, and an isomorphism between the posets $\langle A,\leqsa\rangle$ and $\langle \varphif[A],\subseteq\rangle$.
\end{thm}
\begin{proof}
By definition and the injectivity of $\varphif$ it follows that this map is an isomorphism from $\Al$ onto $\varphifA$.
Moreover,  $a\leqsa b$ implies $\varphif(a)\subseteq\varphif(b)$, and from this and the injectivity of $\varphif$,  we get that $\varphif$ is an order embedding.
\end{proof}

Theorem \ref{thm:srepresentation} is the representation theorem for $\Sm$-algebras we are interested in. It 
was already addressed by Czelakowski in  \cite[Chapter 6]{Cz01}.  
From it we obtain that  the algebra $\varphifA$ is an $\Sm$-algebra, and  
therefore we may consider the closure operator $\Conspa$ associated with the closure system $\FiS(\varphifA)$. 
It is no difficult to see that  the following theorem holds (\cf   \cite[Lem.\ 4.3.6 and Cor.\ 4.3.7]{Es13} for a proof): 

\begin{thm}\label{thm:repr}
Then following facts hold: 
\begin{enumerate}
\item $\{\varphif[P]:P\in\FF \}$ is a closure base for $\FiS(\varphifA)$, 
\item  for any $a\in A$ and any $B\subseteq A$: 
\begin{equation*}
a\in\Consa(B) \IFF 
\varphifm(B)\subseteq \varphif(a) \IFF 
\varphif(a)\in \Conspa (\varphif[B]).
\end{equation*}
\item for all  $B\subseteq A$, $\Conspa (\varphif[B]) = \varphif[\Consa(B)]$.
\item for all $B, B' \subseteq A$,
$\Consa(B) = \Consa(B') \IFF \varphifm(B) = \varphifm(B').$
\end{enumerate}
\end{thm}

In the case of finitary, filter-distributive and congruential logics
 we have the following characterization of the irreducible $\Sm$-filters inside the optimal ones.
  
\begin{prop}\label{prop:irrinsideopt}
Let $\Sm$ be a finitary, filter-distributive and congruential logic,  $\Al$ an $\Sm$-algebra and  $\FF $  a closure base  for $\FiSA$. 
Then, for every $F  \in \OptSA$,
$F \in \IrrSA$ if and only if  $\varphi_{\FF}[F^c] = \{\varphi_{\FF}(a):a \not \in F\}$ is non-empty and up-directed in $\langle \varphif[A],\subseteq\rangle$.
\end{prop}
\begin{proof}
Let $F \in \OptSA$. Note that from Theorem \ref{thm:srepresentation} we get that  $F^{c}$ is an order ideal of $\langle \Al, \leqsa\rangle$ if and only if $\{\varphi_{\FF}(a):a \not \in F\}$ is non-empty and up-directed in $\langle \varphi_{\FF}[A],\subseteq\rangle$. From (2) of Theorem \ref{thm:OptS-IdsS} we have that $F \in \IrrSA$   if and only if $F^{c}$ is an $\leqsa$-order ideal of $\Al$. Thus we obtain the proposition.  
\end{proof}

\subsection{The $\Sm$-semilattice of an $\Sm$-algebra $\Al$}

The $\Sm$-semilattice of an $\Sm$-algebra $\Al$ of a finitary congruential logic $\Sm$  was introduced  in \cite{GeJaPa10}. 
It can be described as the dual of the join-semilattice of the finitely generated $\Sm$-filters of $\Al$. 
Alternatively, to obtain a specific definition of it  we can work with certain closure bases $\FF$  for $\FiSA$, and obtain the $\Sm$-semilattice of $\Al$ as the closure under finite intersections of the image of $A$ under the representation map $\varphif$.

Let $\Sm$ be a congruential logic, $\Al$ an $\Sm$-algebras  and $\FF$  a closure base for the closure system $\FiSA$. We denote by $\MefA$ the closure of $\varphi_{\FF}[A]$ under intersections of finite subsets.  In other words, $\MefA = \{ \varphifm(B): B \subsw A \}$. Note that $\FF \in \MefA$ because $\FF = \varphifm(\emptyset)$. 

\begin{prop}
\label{prop:iso}
If $\FF$ and $\FF'$ are closure bases for $\FiSA$, then $ \langle \MefA,\cap, \FF\rangle$ and $\langle \Mee_{\FF'}(\Al), \cap, \FF'\rangle$ are isomorphic semilattices
\end{prop}
\begin{proof}
The map $H: \MefA \to  \Mee_{\FF'}(\Al)$ defined by $H(\varphifm(B)) = \myhat{\varphi}_{\FF'}(B)$ for every$B \subsw A$ is an isomorphism between $ \langle \MefA,\cap, \FF\rangle$ and $\langle \Mee_{\FF'}(\Al), \cap, \FF'\rangle$.  It is obvious that it is onto. Using using (4) of Theorem \ref{thm:repr}  it easily follows that it is one-to-one. Moreover, 
for any $B,B'\subseteq A$, since $\varphifm(B)\cap\varphifm(B')= \varphifm(B\cup B')$ and $\myhat{\varphi}_{\FF'}(B)\cap \myhat{\varphi}_{\FF'}(B') = \myhat{\varphi}_{\FF'}(B \cup B')$, $H(\varphifm(B)\cap\varphifm(B')) = \myhat{\varphi}_{\FF'}(B)\cap \myhat{\varphi}_{\FF'}(B')$, and since   $\FF = \varphifm(\emptyset)$ and $\FF' = \myhat{\varphi}_{\FF'}(\emptyset)$, $H(\FF) = \FF'$. 
\end{proof}

According to the proposition, we can dispense with the subscript $\FF$ in  $\varphif$ and  $\MefA$. This entitles us to give the next definition.

\begin{defn} \label{defn:Ssemilattice} 
For every  congruential logic $\Sm$  and every $\Sm$-algebra $\Al$, 
the algebra $\MeeA\coloneqq \langle \MeeA,\cap,\FF\rangle$ 
 is called the \emph{{$\Sm$-semilattice} of $\Al$}. 
  \end{defn}

By definition, $\MeeA$ is a meet-semilattice with top element $\FF$. 
 If $\Sm$ has theorems, then $\varphi(1^{\Al})=\FF$. In this case we can describe  $\MeeA$ as the set  $\{ \varphimh(B): B \subsw A \text{ and } B \not = \emptyset\}$.
The following lemma concerns when a  bottom element exists in  $\MeeA$.

\begin{lemma} \label{lemma:stechnical} 
A finite and non-empty $U \subseteq A$ is   a bottom-family  of $\Al$  if and only if $\varphimh(U)$ is a bottom element of $\MeeA$. \end{lemma}
\begin{proof} 
Assume that $U \subseteq A$ is finite and non-empty. Suppose that $\varphimh(U)$ is  a bottom-family of $\Al$. Then $\Consa(U) = A$. From Theorem \ref{thm:repr} follows that for every $a \in A$, $\varphimh(U)\subseteq \varphi(a)$ and this implies that $\varphimh(U)$ is a bottom element of $\MeeA$. To prove the converse 
assume that   $\varphimh(U)$ is a bottom element of $\MeeA$. Then $\varphimh(U) \subseteq \varphi(a)$ for every $a \in A$. Hence, by Theorem \ref{thm:repr} it follows that $\Consa(U)=A$ and hence $U$ is a bottom-family of $\Al$. 
%
\end{proof}

\begin{rem}
Under the assumptions of  Lemma \ref{lemma:stechnical},   if $\Al$ is not trivial and $\MeeA$ has a bottom element, then it is $\varphimh(U)$ where $U$ is any bottom-family of $\Al$. Indeed if $\varphifm(B)$ is a bottom element and $B = \emptyset$, then    
$\varphimh(B) = \FF$ and this implies that $\Al$ is trivial.
\end{rem}

We proceed to study  the relations between the different families of filters and ideals in $\Al$ and $\MeeA$ we have considered so far. 
Let us begin with the $\Sm$-filters of $\Al$ and the  filters of $\MeeA$. The operation of meet filter generation $\Fimg{. }$ is taken in $\MeeA$.
Moreover, we fix a congruential logic $\Sm$, and $\Sm$-algebra and a closure base $\FF$ for $\FiSA$.

\begin{lemma}\label{lemma:me2}
  For every  finite $B \subseteq A$ and every finite family $\{B_i: i \in K\}$ of finite subsets of $A$: \begin{equation*}
\bigcap_{i\in K}\Consa(B_{i})\subseteq \Consa(B)   
\IFF  \bigcap_{i\in K} \Fimg{\varphim(B_{i})}\subseteq \Fimg{\varphim(B)}.
 \end{equation*}
\end{lemma}
\begin{proof} We distinguish two cases depending on wether  $K$ is empty or not. Assume first that $K\not = \emptyset$.
Suppose  that $\bigcap\{\Consa(B_{i}):i\in K\}\subseteq \Consa(B) $ and let $D\subseteq A$ be finite and such that $\varphim(D)\in  \bigcap\{ \Fimg{\varphim(B_{i})}:i\leq n\}$, \ie $\varphim(B_{i})\subseteq \varphim(D)$ for all $i\leq n$. 
Then, by Theorem \ref{thm:repr},  $D\subseteq \Consa(B_{i})$ for all $i\leq n$. 
Thus from the hypothesis follows that $D\subseteq \Consa(B)$ 
and by Theorem \ref{thm:repr} again we get $\varphim(B)\subseteq \varphim(D)$; hence $\varphim(D)\in \Fimg{\varphim(B)}$. 
For the converse, assume that $\bigcap\{ \Fimg{\varphim(B_{i})}:i\in K\}\subseteq \Fimg{\varphim(B)}$ and let $a\in\bigcap \{\Consa(B_{i}):i\in K\}$. 
Then for each $i\leq n$, $a\in\Consa(B_{i})$, and so, by Theorem \ref{thm:repr}, $\varphim(B_{i})\subseteq \varphi(a)$. 
This implies that $\varphi(a)\in \bigcap\{ \Fimg{\varphim(B_{i})}:i\in K\}$, and so, from the hypothesis follows that $\varphi(a)\in\Fimg{\varphim(B)}$, 
\ie  $\varphim(B)\subseteq \varphi(a)$. 
Then  by Theorem \ref{thm:repr} again, we get $a\in\Consa(B)$. 

Now we assume that $K = \emptyset$. If $\bigcap\{\Consa(B_{i}):i\in K\} \subseteq \Consa(B)$, then $\Consa(B) = A$ and hence $B$ is a bottom-family and $\varphim(B)$ is a bottom element of $\MeeA$. Therefore, $\Fimg{\varphim(B)} = \MeeA$ and thus $\bigcap_{i\in K} \Fimg{\varphim(B_{i})}\subseteq \Fimg{\varphim(B)}$. Conversely,  if  $\bigcap_{i\in K} \Fimg{\varphim(B_{i})}\subseteq \Fimg{\varphim(B)}$, then $\Fimg{\varphim(B)} = \MeeA$; thus $\varphim(B)$ is the bottom element of $\MeeA$.  Using Theorem  \ref{thm:repr} it follows that $\Consa(B) = A$.
Hence, $\bigcap\{\Consa(B_{i}):i\in K\} \subseteq \Consa(B)$.
\end{proof}

The next proposition (cf. \cite[Lemmas 4.5 and 4.8]{Es13}) shows the  relation between the $\Sm$-filters of $\Al$ and the filters of $\MeeA$ provided $\Sm$ is finitary.

\begin{prop} \label{prop:mefilter0} 
 Let $\Sm$ be a finitary congruential logic, $\Al$  an $\Sm$-algebra and $\FF$ a closure base for  $\FiSA$.  
\begin{enumerate}
	\item If $F\in \FiSA $, then  $\Fimg{\varphi[F]} \in \Fim(\MeeA)$  and 
 $\varphi^{-1}\big[\Fimg{\varphi[F]}\big]=F$.

	\item If $F \in \Fim(\MeeA)$, then 
	 $\varphi^{-1}[F] \in \FiSA$   and 
	 $\Fimg{F\cap\varphi[A]}=F$.		
\end{enumerate}
\end{prop}
\begin{proof}
(1) Let $F \in \FiSA$. 
By definition $\Fimg{\varphi[F]} \in \Fim(\MeeA)$ and it is clear that $F \subseteq \varphi^{-1}\big[\Fimg{\varphi[F]}\big]$. To prove  the other inclusion let 
$a \in \varphi^{-1}\big[\Fimg{\varphi[F]}\big]$, \ie $\varphi(a) \in  \Fimg{\varphi[F]}$. 
Then there is a finite $B \subseteq F$ such that $\varphim(B) \subseteq \varphi(a)$. Hence, by Theorem \ref{thm:repr} we have $a \in \Consa(B) \subseteq F$. 

 (2) Let $F \in \Fim(\MeeA)$. 
 By definition $\varphi^{-1}[F]\subseteq\Consa(\varphi^{-1}[F])$. Let now
 $a \in \Consa(\varphi^{-1}[F])$. 
 By finitarity of $\Sm$,  let  $B \subseteq  \varphi^{-1}[F]$ be finite with $a \in \Consa(B)$. 
 Then by Theorem \ref{thm:repr}  we have $\varphim(B) \subseteq \varphi(a)$. 
 Since $\varphi[B] \subseteq F$,  $\varphim(B) = \bigcap  \varphi[B] \in F$, 
 and since $F$ is an up-set,  $\varphi(a) \in F$. Hence $a \in \varphi^{-1}[F]$. We conclude that 
 $\varphi^{-1}[F] = \Consa(\varphi^{-1}[F])$ and thus that $\Fimg{\varphi[F]} \in \Fim(\MeeA)$.
  
 For the remaining part of the statement, note that $\Fimg{\varphi[\varphi^{-1}[F]]} =  \Fimg{F \cap \varphi[A]}$, having then $\Fimg{F \cap \varphi[A]} \subseteq F$. 
 For the other inclusion  let $B \subseteq A$ be such that $\varphim(B) \in F$. If $B = \emptyset$, then $\varphim(B) = \FF \in \Fimg{F \cap \varphi[A]}$. If $B \not = \emptyset$, then  $\varphi[B] \subseteq F \cap \varphi[A]$ and hence $\varphim(B) = \bigcap \varphi[B]   \in \Fimg{F \cap \varphi[A]}$.  
\end{proof}

\begin{cor}
 \label{cor:mefilter0}
For every finitary  congruential logic $\Sm$  and every $\Sm$-algebra $\Al$, 
there is an isomorphism between the lattice $\FiSA$ and the lattice $\Fim(\MeeA)$ given by the following maps, each one inverse of the other:
\begin{align*}
\Fimg{\varphi[\phantom{.}]}:\FiSA\,\,\cong\,\, \Fim(\MeeA):\varphi^{-1}[\phantom{.}].
\end{align*}
\end{cor}

As a consequence, we obtain that for every  finitary congruential logic $\Sm$ and every $\Sm$-algebra $\Al$, the lattice $\FiSA$ is distributive if and  only if the lattice $\Fim(\MeeA)$ is distributive. 
Therefore if $\Sm$ is in addition a filter-distributive logic, then $\MeeA$ is  a distributive meet-semilattice. 
Moreover, the previous isomorphism maps irreducible $\Sm$-filters of $\Al$ to irreducible filters of $\MeeA$. 
But even  if this happens,  it  may not be an isomorphism between the optimal $\Sm$-filters of $\Al$ and the optimal filters of $\MeeA$. However, under the assumption of filter-distributivity of the logic, it is  an isomorphism.
In order to show it, 
since the behaviour of optimal $\Sm$-filters and optimal  filters depends on the behaviour of the strong $\Sm$-ideals of $\Al$ and  the Frink ideals of $\MeeA$,  
we need to make a detour and study first  the relation between strong $\Sm$-ideals of $\Al$ and Frink ideals of $\MeeA$. 
Throughout the next proofs
we use ${\downarrow}$ instead of ${\downarrow}_{\MeeA}$. The operation $\Idfg{.}$  of Frink ideal generation is taken in $\MeeA$.

\begin{prop} \label{prop:meideals} 
Let $\Sm$ be a  congruential logic,  $\Al$  an $\Sm$-algebra and   $\FF$ a closure base for  $\FiSA$:
\begin{enumerate}
\item For every  strong $\Sm$-ideal $I$ of $\Al$, $\Idfg{\varphi[I]}={\downarrow}_{\MeeA}\varphi[I]$. 
	\item If $I\in\IdsS(\Al)$, then  $\Idfg{\varphi[I]}\in \IdF(\MeeA)$,   
 $\varphi^{-1}[\Idfg{\varphi[I]}]=I$ and if $I$ is $\wedge$-prime, then $\Idfg{\varphi[I]}$ is $\wedge$-prime.
	\item If $I\in\IdF(\MeeA)$ is $\wedge$-prime, then  $\varphi^{-1}[I]\in \IdsS(\Al)$, ${\Idfg{\varphi[\varphi^{-1}[I]]}=I}$ and  it is $\Sm$-prime.	
\end{enumerate}
\end{prop}
\begin{proof}
(1)\,\, Let $I\in\IdsSA$. If $I = \emptyset$, then $\Al$ has no bottom-family and hence $\MeeA$ has no bottom element, which implies that $\emptyset$ is a Frink ideal of $\MeeA$. Thus  $\Idfg{\varphi[I]} = \emptyset$ and we are done. If  $I$ is non-empty, 
it is enough to show that ${\downarrow} \varphi[I]$ is a Frink ideal. 
So let  $K$ be a finite set of indexes and $\{B_{i}: i \in K\}$ a family of finite subsets of $A$  such that $\varphim(B_{i})\in{\downarrow} \varphi[I]$ for all $i\in K$. Then for every $i \in K$ there exists $a_{i}\in I$  such that $\varphim(B_{i})\subseteq \varphi(a_{i})$. 
Moreover, let  $B\subseteq A$ be  finite and such that $\bigcap\{\Fimg{\varphim(B_{i})}:i\in K\}\subseteq \Fimg{\varphim(B)}$. 
If $K=\emptyset$, then  $\Fimg{\varphim(B)}=\Meea$, and so $\varphim(B)$ is the bottom element of $\MeeA$ that belongs to ${\downarrow} \varphi[I]$ since $I$ is non-empty. 
If $K\not =\emptyset$, then  $\bigcap\{\Fimg{\varphi(a_{i})}:i\in K\}\subseteq \Fimg{\varphim(B)}$. 
From Lemma \ref{lemma:me2} it follows $\bigcap\{\Consa(a_{i}):i\in K\}\subseteq \Consa(B)$, 
and then since $\{a_{i}:i\in K\}\subseteq I$ and $I$ is a strong $\Sm$-ideal, 
there exists $c\in \Consa(B)\cap I\neq \emptyset$. 
Then by Theorem \ref{thm:repr}  $\varphim(B)\subseteq\varphi(c)\in\varphi[I]$, as required.

(2)\,\, Let $I\in\IdsSA$. 
By definition $\Idfg{\varphi[I]}\in\IdF(\MeeA)$, and  clearly $I\subseteq \varphi^{-1}[\Idfg{\varphi[I]}]$. 
Let us prove the other inclusion.  
If $I=\emptyset$ we get as in the previous proof that $\Idfg{\varphi[I]} = \emptyset$ and we are done. Assume $I\neq\emptyset$, and let $a\in \varphi^{-1}[\Idfg{\varphi[I]}]$, \ie $\varphi(a)\in \Idfg{\varphi[I]}$. 
By the characterization of the   generated Frink ideal, there is a finite $I'\subseteq  I$ such that $\bigcap\{\Fimg{\varphi(b)}:b\in I'\}\subseteq \Fimg{\varphi(a)}$. 
As $I\neq\emptyset$, we can assume, without loss of generality, that $I'\neq\emptyset$. 
Then by Lemma \ref{lemma:me2}, $\bigcap\{\Consa(b):b\in I'\}\subseteq \Consa(a)$. 
And since $I$ is an $\Sm$-ideal, we get $a\in I$, as required. 

Assume now that $I$ is $\Sm$-prime. Thus $I$ is proper. Using that  $\varphi^{-1}[\Idfg{\varphi[I]}]=I$, $\Idfg{\varphi[I]}$ should be also proper. To prove that $\Idfg{\varphi[I]}$ is $\wedge$-prime,   let  $B_{1},B_{2}\subseteq  A$ be  finite  and such that $\varphim(B_{1})\cap \varphim(B_{2})\in \Idfg{\varphi[I]}$. 
Using that $\varphim(B_{1})\cap \varphim(B_{2})=\varphim(B_{1}\cup B_{2})$ and $\Idfg{\varphi[I]}={\downarrow} \varphi[I]$, 
we get that there is $c\in I$ such that $\varphim(B_{1}\cup B_{2})\subseteq \varphi(c)$. This implies that $B_1 \cup B_2$ is non-empty, otherwise $c$ is a bottom element and $I$ is not proper.
Then by Theorem \ref{thm:repr}, $c\in \Consa(B_{1}\cup B_{2})$, so $\Consa(B_{1}\cup B_{2})\cap I\neq \emptyset$. 
Moreover, since $I$ is $\Sm$-prime, we get $(B_{1}\cup B_{2})\cap I\neq \emptyset$, so $B_{1}\cap I\neq\emptyset$ or $B_{2}\cap I\neq\emptyset$. 
This implies, by Theorem \ref{thm:repr} again, that either $\varphim(B_{1})\in {\downarrow}{\varphi[I]}$ or $\varphim(B_{2})\in {\downarrow}{\varphi[I]}$. We conclude that  $\Idfg{\varphi[I]}$ is $\wedge$-prime.

(3)\,\, Let now $I\in\IdF(\MeeA)$ be $\wedge$-prime. Thus $I$ is proper. 
We show that the condition in the definition of strong $\Sm$-filter holds for $\varphi^{-1}[I]$. By Lemma \ref{lem:useful-proof-strongideal} this implies that $\varphi^{-1}[I]$ is an $\Sm$-ideal because $\varphi^{-1}[I]$ is easily seen to be a downset. Let $I'\subseteq \varphi^{-1}[I]$ be finite and let $C\subseteq A$ be finite, non-empty, and such that $\bigcap\{\Consa(b):b\in I'\}\subseteq \Consa(C)$. By Lemma \ref{lemma:me2} we have  $\bigcap\{\Fimg{\varphim(b)}: b \in I'\}\subseteq \Fimg{\varphim(C)}$. This implies that $\varphim(C) \in I$. Now, by $\wedge$-primeness of $I$ we get that  $\varphi(c)\in I$ for some $c\in C$ and so $\Consa(C)\cap \varphi^{-1}[I]\neq\emptyset$.

We proceed to  show that $\Idfg{\varphi[\varphi^{-1}[I]]}= I$. 
Clearly the inclusion from left to right holds, so we just have to prove the other inclusion. Since $I$ is proper, $\varphim(\emptyset) \not \in I$. 
Let $B\subseteq  A$ be non-empty, finite and such that $\varphim(B)\in I$. 
Then, as $I$ is $\wedge$-prime, there is $b\in B$, such that $\varphi(b)\in I$. 
So $\varphi(b)\in\varphi[\varphi^{-1}[I]]$ and as $\varphim(B)\subseteq \varphi(b)$ and Frink ideals are down-sets, we have $\varphim(B)\in \Idfg{\varphi[\varphi^{-1}[I]]}$.

Now we prove that  $\varphi^{-1}[I]$ is proper. In view of a contradiction, suppose that it is not. Then $\varphi[\varphi^{-1}[I]] = \varphi[A]$ and since, as we already proved, ${\Idfg{\varphi[\varphi^{-1}[I]]}=I}$, we have, ${\Idfg{\varphi[A]}=I}$. We show that $\Idfg{\varphi[A]} = \MeeA$. Let $B$ be a finite subset of $A$. Observe that $$\bigcap_{a \in B} {\uparrow}_{\MeeA}\varphi(a) \subseteq {\uparrow}_{\MeeA}\varphim(B).$$ Since $\{\varphi(a): a \in B\}$ is a finite subset of $\varphi[A]$ it follows that $\varphim(B) \in \Idfg{\varphi[A]}$. Hence, $\Idfg{\varphi[A]} = \MeeA$ and therefore $I$ is not proper.  To conclude, we show that  $\varphi^{-1}[I]$ is  $\Sm$-prime. 
Let  $B\subseteq  A$ be non-empty, finite and such that $\Consa(B)\cap \varphi^{-1}[I]\neq\emptyset$, and let $c\in \Consa(B)\cap \varphi^{-1}[I]$. 
As $c\in\Consa(B)$, then by Theorem \ref{thm:repr} $\varphim(B)\subseteq \varphi(c)$. 
Moreover, since $\varphi(c)\in I$ and $I$ is a down-set, we get $\varphim(B)\in I$. 
Now, as $I$ is $\wedge$-prime, there is $b\in B$ such that $\varphi(b)\in I$, so $B\cap  \varphi^{-1}[I]\neq\emptyset$, as required. 
\end{proof}

\begin{cor}
The map $\Idfg{\varphi[.]}$  suitably restricted establishes an isomorphism between the poset of the $\Sm$-prime and strong $\Sm$-ideals of $\Al$ and the poset of the $\wedge$-prime Frink ideals of $\MeeA$; and the inverse map  is given by  $\varphi^{-1}$.
\end{cor}
\begin{proof}
It follows from Proposition \ref{prop:meideals}.
\end{proof}

Up to this point, all the results in the present section are valid in general for any  congruential logic, except Proposition \ref{prop:mefilter0} and Corollary \ref{cor:mefilter0} where we required the logic to be in addition finitary. If we assume besides that the logic is filter-distributive, then we get further results. 


\begin{prop} \label{prop:mefilter} 
Let $\Sm$ be a filter-distributive and finitary congruential logic,  $\Al$ an $\Sm$-algebra, and $\FF$ a closure base for  $\FiSA$:
\begin{enumerate}
	\item If $F\in\OptSA$ is non-empty, then $\Fimg{\varphi[F]}\in\Optm(\MeeA)$.
	\item If $F\in \Optm(\MeeA)$, then $\varphi^{-1}[F]\in\OptSA$.
\end{enumerate}
\end{prop}

\begin{proof}
(1)\,\, Let  $F\in\FiSA$ be optimal and non-empty. 
Then by Theorem \ref{thm:OptS-IdsS} and Lemma \ref{lemma:FiS-Sprime}, $F^{c}$ is an $\Sm$-prime strong $\Sm$-ideal of $\Al$, and so by Proposition \ref{prop:meideals} $\Idfg{\varphi[F^{c}]}$ is a $\wedge$-prime Frink ideal of $\MeeA$, hence proper, and moreover by Proposition \ref{cor:Optm-IdF} $\Idfg{\varphi[F^{c}]}^{c}$ is an optimal  filter of $\MeeA$. 
Therefore, it is enough to show that $\Fimg{\varphi[F]}=\Idfg{\varphi[F^{c}]}^{c}$. 
To this end, we prove first the inclusion from right to left. 
Let  $B\subseteq  A$ be finite and such that $\varphim(B)\in \Idfg{\varphi[F^{c}]}^{c}$. If $B = \emptyset$, then  $\varphim(B) = \FF$ and hence $\varphim(B)\in  \Fimg{\varphi[F]}$. If $B$ is non-empty, 
then for all $b\in B$,  $\varphi(b)\notin \Idfg{\varphi[F^{c}]}$, and 
thus by  Proposition \ref{prop:meideals} we get that 
$b\notin \varphi^{-1}[\Idfg{\varphi[F^{c}]}]=F^{c}$. 
Therefore  $\varphi(b)\in \varphi[F]$ for all $b\in B$, and thus 
$\bigcap\{\varphi(b):b\in B\}=\varphim(B)\in \Fimg{\varphi[F]}$. 
For the other inclusion,  let $B\subseteq  A$ be  finite  and such that $\varphim(B)\in \Fimg{\varphi[F]}$. 
Then either $\varphim(B)=\FF$ or there is a non-empty and finite $B'\subseteq  F$ such that $\varphim(B')\subseteq \varphim(B)$. 
In the first case, since $\Idfg{\varphi[F^{c}]}$ is proper, $\varphim(B) \notin \Idfg{\varphi[F^{c}]}$, and we are done. In the second case,  since $\varphim(B')\subseteq \varphim(B)$ implies that $\Fimg{\varphim(B)}\subseteq \Fimg{\varphim(B')}$,  by Lemma \ref{lemma:me2}, $\Consa(B)\subseteq\Consa(B')\subseteq F$, so $B\subseteq F$.  
Therefore for all $b\in B$, $b\notin F^{c}=\varphi^{-1}[\Idfg{\varphi[F^{c}]}]$. 
Hence $\varphi(b) \notin \Idfg{\varphi[F^{c}]}$ for all $b\in B$. 
Moreover, since  $\Idfg{\varphi[F^{c}]}$ is a $\wedge$-prime Frink ideal,  $\varphim(B)\notin \Idfg{\varphi[F^{c}]}$, \ie $\varphim(B)\in \Idfg{\varphi[F^{c}]}^{c}$, as required.

(2)\,\, For $F\in\Optm(\MeeA)$,  
by Corollary \ref{cor:Optm-IdF} $F^{c}$ is a $\wedge$-prime Frink ideal of $\MeeA$, 
and so by Proposition \ref{prop:meideals}, 
$\varphi^{-1}[F^{c}]$ is an $\Sm$-prime strong $\Sm$-ideal of $\Al$, 
and moreover by Lemma \ref{lemma:FiS-Sprime} and Theorem \ref{thm:OptS-IdsS},
 $(\varphi^{-1}[F^{c}])^{c}$ is an optimal $\Sm$-filter of $\Al$. 
 Notice that $(\varphi^{-1}[F^{c}])^{c}=\varphi^{-1}[F]$. 
Therefore $\varphi^{-1}[F]$ is an optimal $\Sm$-filter of $\Al$.
\end{proof}

\begin{cor} \label{cor:mefilter}
For every filter-distributive and finitary  congruential logic $\Sm$  and every $\Sm$-algebra $\Al$, 
there is an order isomorphism between the poset $\langle\OptSA \setminus \{\emptyset\},\subseteq\rangle$ and the poset $\langle\Optm(\MeeA),\subseteq\rangle$ given by the following maps, one inverse to the other:
\begin{align*}
\Fimg{\varphi[\cdot]}:\langle\OptSA\setminus \{\emptyset\},\subseteq\rangle\,\,\cong\,\, \langle\Optm(\MeeA),\subseteq\rangle:\varphi^{-1}[\cdot].
\end{align*}
\end{cor}

Note that if $\Sm$ has theorems, then for every algebra $\Al$ all the $\Sm$-filters are non-empty, in particular so are the optimal filters of the $\Sm$-algebras. Hence in this case the corollary establishes  an isomorphism between $\langle\OptSA,\subseteq\rangle$ and  $\langle\Optm(\MeeA),\subseteq\rangle$.

Also note that the previous isomorphism provides a new characterization of the non-empty optimal $\Sm$-filters of $\Al$ 
as the images of optimal filters of $\MeeA$ by the map $\varphi^{-1}[\cdot]$. 
This is the keystone why we build the 
 dual Priestley space of $\Al$ from the dual Priestley space of $\MeeA$, as it is explained in next section.

\subsection{\Sp spaces}\label{sec:Obj}

We proceed to present the correspondence between $\Sm$-algebras and a certain class of Priestley-style spaces for the logics $\Sm$ that are finitary, congruential, filter-distributive and with theorems. This spaces will be named \Sp spaces.
In order to characterize them, we  use the concept of  \emph{referential algebra}, that goes back to W\'ojcicki \cite{Wo88} (see also \cite{JaPa06}). 
%
%

Given a logical language  $\mathscr{L}$, 
an \emph{$\mathscr{L}$-referential algebra}   is a structure $\ralg{X}=\langle X,\BB\rangle$ where $X$ is a  set and $\BB$ is an $\mathscr{L}$-algebra whose elements are subsets of $X$.\footnote{ We admit $X$ to be empty to cover the case of the trivial algebras. In this case the domain of $\BB$ is $\{\emptyset\}$ and $\BB$ is trivial.}

For any $\mathscr{L}$-referential algebra $\ralg{X}=\langle X,\BB\rangle$,  the relation ${\preceq_{\ralg{X}}}\subseteq X\times X$ defined by setting for every $u, v \in X$:
\begin{equation*}
u \preceq_{\ralg{X}} v \IFF (\forall U\in B)\big(u\in U \Rightarrow  v \in U\big)
\end{equation*}
 is a quasiorder on $X$. Whenever $\preceq_{\ralg{X}}$ is a partial order, the \mbox{$\mathscr{L}$-referential} algebra $\ralg{X}$ is said to be \emph{reduced}. 
In this case, we  denote $\preceq_{\ralg{X}}$ by $\leq_{\ralg{X}}$, or even by $\leq$ when the context is clear.\footnote{ If $X$ is empty, we consider the empty relation as a partial order and the referential algebra $\langle X,\BB\rangle$  is reduced.} 

Referential algebras  can be used to define logics in the following way. 
 For any $\mathscr{L}$-referential algebra $\ralg{X}=\langle X,\BB\rangle$ we  define the relation ${\vdash_{\ralg{X}}}\subseteq {\pws(Fm_{\mathscr{L}})\times Fm_{\mathscr{L}}}$ such that for all $\Gamma\cup\{\delta\}\subseteq Fm_{\mathscr{L}}$:
\[
\Gamma\vdash_{\ralg{X}}\delta \IFF (\forall h\in\Hom(\Fm_{\mathscr{L}},\BB)) \bigcap_{\gamma\in\Gamma}h(\gamma)\subseteq h(\delta).
\]
This relation is such that  
$\langle Fm_{\mathscr{L}}, \vdash_{\ralg{X}}\rangle$ is a logic.  

Given a logic $\Sm$ and an $\mathscr{L}_{\Sm}$-referential algebra $\ralg{X}$, we say that $\ralg{X}$ is an \emph{$\Sm$-referential algebra}
 provided that $\vdash_{\Sm}\,\subseteq\, \vdash_{\ralg{X}}$. 
Moreover, we say that \emph{$\Sm$ admits a (complete local) referential semantics} 
 if there is a class of referential algebras $\clas{X}$ such that 
${\vdash_{\Sm}}={\bigcap\{\vdash_{\ralg{X}}:\ralg{X}\in\clas{X}\}}$. 
 
 It is easy to see that if  $\ralg{X}=\langle X,\BB\rangle$ is an $\Sm$-referential algebra, then for every $x \in X$, the set $\{U \in B: x \in U\}$ is and $\Sm$-filter. Moreover, if $\langle X,\BB\rangle$ is a  reduced $\Sm$-referential algebra, then $\BB\in\Alg\Sm$ (see  \cite[Remark 5.2]{JaPa06}).

We return now to consider $\Sm$-algebras   and closure bases  for their closure systems of $\Sm$-filters,
as they can be seen as reduced $\Sm$-referential algebras when $\Sm$ is congruential.

\begin{thm} \label{thm:ralg} 
Let $\Sm$ be a congruential logic, $\Al$  an $\Sm$-algebra, and   $\FF$  a closure base  for $\FiSA$.  
Then $\langle \FF,\varphifA\rangle $ is a reduced $\Sm$-referential algebra and the associated order is given by the inclusion relation. 
\end{thm}
\begin{proof} By definition, $\langle \FF ,\varphifA\rangle$ is a referential algebra.
We show first that it is reduced. 
Consider the quasiorder ${\preceq}\subseteq \FF \times \FF $ of this referential algebra, and note that for every $P, Q \in \FF$, $P\preceq Q$ if and only if for every $a\in A$ such that $P\in \varphif(a)$,  $Q\in\varphif(a).$
 It follows that  $\preceq$ is the inclusion relation on $\FF$. Therefore the referential algebra is reduced. 
Let us show now that $\langle \FF ,\varphifA\rangle$ is an $\Sm$-referential algebra.
Let $\Gamma\cup\{\delta\}\subseteq Fm_{\mathscr{L}}$ be such that $\Gamma\vdash_{\Sm}\delta$, and 
let $h\in\Hom(\Fm,\varphifA)$.
Since $\varphif\in\Hom(\Al,\varphifA)$ is an isomorphism, there is $h'\in\Hom(\Fm,\Al)$ such that $\varphif\circ h'=h$. 
To show that $\bigcap\{h(\gamma):\gamma\in\Gamma\}\subseteq h(\delta)$,  suppose that $P\in \FF $ is such that $P\in \bigcap\{h(\gamma):\gamma\in\Gamma\}=\bigcap\{\varphif(h'(\gamma)):\gamma\in\Gamma\}$. 
Then $h'[\Gamma]\subseteq P$. 
And since $P\in \FiS(\Al)$, $h'\in\Hom(\Fm,\Al)$ and $\Gamma\vdash_{\Sm}\delta$ we obtain $h'(\delta)\in P$, so $P\in \varphif(h'(\delta))=h(\delta)$, as required.
\end{proof}

In   \cite[Section 5.6.7]{Wo88} the referential algebra $\langle \FF,\varphif[\Al]\rangle$ is called the \emph{canonical referential algebra for $\Consa$ determined by $\FF$}. 
Notice that Theorem \ref{thm:ralg} and Theorem \ref{thm:repr} together imply that for every congruential logic $\Sm$, 
there is a  correspondence between the reduced $\Sm$-referential algebras and the structures of the form $\langle \Al,\FF\rangle$, 
where $\Al$ is an $\Sm$-algebra and $\FF$ is a closure base for $\FiSA$. 
This correspondence between objects, first addressed by Czelakowski in \cite{Cz01}, 
was formulated as a full-fledged duality in \cite{JaPa06}, for the case when the collection $\FiSA$ is taken itself as the closure base. 
But this closure base is not the closure base that properly generalizes the representation theorems on which  the Stone/Priestley dualities that we find in the literature  are based. 
For instance, the algebraic counterpart of intuitionistic logic is the  variety of Heyting algebras and the intuitionistic logical filters  of the Heyting algebras are the lattice filters. Yet the representation theorem on which  the Esakia duality for Heyting algebras is based relays  on the prime  filters and not on all the lattice filters. 
Therefore, for our purposes  we should not  work with the whole collection of $\Sm$-filters, but rather we should identify a closure base that provides us with a direct generalization of the mentioned representation theorems.  The base we need is  the collection of all optimal $\Sm$-filters.

From now on, we let 
\textit{$\Sm$ to be a filter-distributive finitary congruential logic with theorems} and $\Al$ an $\Sm$-algebra, if we do not say otherwise. From Proposition \ref{prop:ofls} follows that $\OptSA$ is a closure base for $\FiSA$ and we 
 further assume that $\OptSA$ is the  closure base from which the representation map  $\varphi$ is defined, \ie for every $a\in A$, $\varphi(a)=\{P\in\OptSA:a\in P\}$. Since $\Sm$ has theorems, $\Al$ has a top element that we denote by $1^{\Al}$, having then  that  $\varphi(1^{\Al}) = \OptSA$.
We define on $\OptSA$ the topology $\tau_{\Al}$ obtained by taking as subbasis the collection:
\begin{equation*}
\{\varphi(a):a\in A\}\cup\{\varphi(b)^{c}:b\in A\}.
\end{equation*}

\begin{prop}\label{prop:isom1}
The isomorphism $\Fimg{\varphi[.]}:\FiSA \to  \Fim(\MeeA)$ suitably restricted establishes an order homeomorphism  between  the ordered topological  spaces 
$\langle\OptSA,\tau_{\Al},\subseteq\rangle$ and ${\langle \Optm(\Mee(\Al)),\tau_{\Mee(\Al)},\subseteq\rangle}$, whose inverse is $\varphi^{-1}[\cdot]$.
\end{prop}
\begin{proof}
By Corollary \ref{cor:mefilter} we already know  that $\Fimg{\varphi[\cdot]}$ establishes an order isomorphism between $\langle\OptSA,\subseteq\rangle$ and ${\langle \Optm(\Mee(\Al)),\subseteq\rangle}$, whose inverse is $\varphi^{-1}[\cdot]$. 
Therefore, to prove the proposition, using that inverse maps preserve intersections, 
we just need to show that $\varphi^{-1}[\cdot]$ sends subbasic opens of the space $\langle \Optm(\Mee(\Al)),\tau_{\Mee(\Al)},\subseteq\rangle$ to opens of $\langle\OptSA,\tau_{\Al},\subseteq\rangle$.  
Recall that $\{\sigma(U):U\in \Mee(A)\}\cup\{\sigma(V)^{c}:V\in \Mee(A)\}$ is a subbasis for $\tau_{\Mee(\Al)}$.  Using 
Definition \ref{defn:Ssemilattice} this subasis can be described as the union of  
$\{\sigma(\varphim(B)):B\subseteq A \text{ finite  }\}$ and $\{\sigma(\varphim(B))^{c}:B\subseteq A \text{  finite  }\}.$
So  let $B\subseteq A$ be  finite. If $B = \emptyset$, then $\varphim(B) = \OptSA$ and so $\sigma(\varphim(B)) = \Optm(\Mee(\Al))$. Therefore, $\varphi^{-1}[\sigma(\varphim(B))] = \OptSA$, which is open, and $\varphi^{-1}[\sigma(\varphim(B))^c] = \emptyset$, which is also open. In the case that  $B$ is non-empty, we prove that $\varphi^{-1}[\sigma(\varphim(B))] =\bigcap\{\varphi(b):b\in B\}$, which implies that  $\varphi^{-1}[\sigma(\varphim(B))]$ is a basic open subset of $\langle\OptSA,\tau_{\Al},\subseteq\rangle$ and $\varphi^{-1}[\sigma(\varphim(B))^{c}]=\bigcup\{\varphi(b)^{c}:b\in B\}$ an open subset of  $\langle\OptSA,\tau_{\Al},\subseteq\rangle$. First note that if  $F\in\OptSA$, then, since $F = \varphi^{-1}[\Fimg{\varphi[F]}]$, $F\in \varphi^{-1}[\sigma(\varphim(B))] $ if and only if $\Fimg{\varphi[F]}\in \sigma(\varphim(B))$. But,
$\Fimg{\varphi[F]}\in \sigma(\varphim(B))$ if and only if $\varphim(B)\in \Fimg{\varphi[F]}$ and this is equivalent to say that $\varphi(b)\in\Fimg{\varphi[F]}$ for every $b \in B$, which is easily seen to be equivalent to say that $F \in \bigcap\{\varphi(b):b\in B\}$.
\end{proof}

\begin{cor}\label{cor:isom1}
Let $\Al$ be an $\Sm$-algebra. Then
\begin{enumerate}
\item  the space $\langle\OptSA,\tau_{\Al},\subseteq\rangle$ is a Priestley space, 
\item the set $\IrrSA$ is dense in $\langle\OptSA,\tau_{\Al},\subseteq\rangle$, 
\item $\IrrSA = \{P \in \OptSA: \{\varphi(a): a \not \in P\} \text{ is non-empty and}$

\hspace{13mm}$\text{ up-directed}\}$.
\end{enumerate}
\end{cor}
\begin{proof}
(1)  follows from the fact that  ${\langle \Optm(\Mee(\Al)),\tau_{\Mee(\Al)},\subseteq\rangle}$ is a Priestley space,  (2) follows from the fact that  the set $\Irrm(\Mee(\Al))$ is  dense  in that space and  that $\IrrSA = \{\varphi^{-1}[F]: F \in \Irrm(\Mee(\Al))\}$, and (3) is a restatement of  Proposition \ref{prop:irrinsideopt}.
\end{proof}

\begin{rem}
Note that for every $P \in \OptSA$,  $\{\varphi(a): a \not \in P\}$ is non-empty and up-directed if and only if $\{\varphim(B) \in  \Mee(\Al):  P \not \in \varphim(B)\}$ is non-empty and up-directed.
\end{rem}

\begin{rem}
In the case of a trivial $\Sm$-algebra $\Al$, since $\Sm$ has theorems the set $\OptSA$ is empty. This forces us to consider the Priestley space with an empty  set of points. 
\end{rem}

We say that a clopen-upset $U$ of $\langle\OptSA,\tau_{\Al},\subseteq\rangle$ is \emph{$\IrrSA $-admissible} whenever  $\max(U^{c})\subseteq \IrrSA$. The next proposition shows that the set of $\IrrSA $-admissible clopen-upsets of $\langle\OptSA,\tau_{\Al},\subseteq\rangle$ is the closure of $\varphi[A]$ under  the binary operation of intersection.

\begin{prop}\label{prop:isom2}
Let $\Al$ be an $\Sm$-algebra and  $U$ be a clopen up-set of the Priestley space
${\langle \OptSA,\tau_{\Al},\subseteq\rangle}$. 
Then $U=\varphim(C)$ for some  finite $C\subseteq A$ if and only if $\max(U^{c})\subseteq  \IrrSA$.
\end{prop}
\begin{proof}
First note that if $U = \OptSA$, then $\varphim(\emptyset) =U$ and $U^{c} = \emptyset$ and therefore the statment holds. Let $U \not = \OptSA$. Suppose first that $U=\varphim(C)$ for some  finite $C\subseteq A$. Then $C \not = \emptyset$ and  there is $P\in\max(\varphim(C)^{c})$, because $U^{c}$ is clopen and non-empty.  
Hence,  $C\nsubseteq P$. Therefore  there is $b\in C\setminus P$. 
Then by the irreducible $\Sm$-filter lemma, 
there is $Q\in\IrrSA$ such that $b\notin Q$ and $P\subseteq Q$. 
This implies $C\nsubseteq Q$, so  $Q\in \varphim(C)^{c}$ and by maximality of $P$  we conclude $P=Q$, \ie $P$ is an irreducible $\Sm$-filter. 
For the converse, let $U$ be a clopen up-set such that $\max(U^{c})\subseteq \IrrSA$. 
Notice that
$$ \Fimg{\varphi[P]}\in\max(\{\Fimg{\varphi[F]}:F\in U\}^{c})\IFF  
P\in\max(U^{c}).$$
This follows from the isomorphism between $\langle \OptSA,\subseteq\rangle$ and $\langle \Optm(\Mee(\Al)),\subseteq\rangle$ given in Proposition \ref{prop:mefilter}.
Therefore, using the homeomorphism given in Proposition \ref{prop:isom1}, from $U$ being an $\IrrSA$-admissible clopen up-set of $\OptSA$ we obtain that  $\{\Fimg{\varphi[F]}:F\in U\}$ is an $\Irrm(\Mee(\Al))$-admissible clopen up-set of $\Optm(\Mee(\Al))$. 
And then by Proposition \ref{prop:gps1}, there is a  finite subset $C\subseteq A$ such that $\sigma_{\Mee(\Al)}(\varphim_{\Al}(C))=\{\Fimg{\varphi[F]}:F\in U\}$, and then we obtain that $U=\varphim_{\Al}(C)$, as required.
\end{proof}

We are ready to introduce the definition of the dual objects of the $\Sm$-algebras.

\begin{defn}\label{defn:Prspace}
A structure $\spa{X}=\langle X,\tau,\BB\rangle$ is an \emph{\Sp space} when: 
\begin{itemize}
	\labitem{(Pr1)}{itm:Pr1}  $\langle X,\BB\rangle$ is a reduced $\Sm$-referential algebra, whose associated order is denoted by $\leq$,
	\labitem{(Pr2)}{itm:Pr2} for any non-empty and finite $\mathcal{V}\subseteq B$ and any $U\in B$, if $\bigcap\mathcal{V}\subseteq U$, then  $U\in \Consb(\mathcal{V})$,
	\labitem{(Pr3)}{itm:Pr3} $\langle X,\tau\rangle$ is a compact space, 
	\labitem{(Pr4)}{itm:Pr4} $B$ is a family of clopen up-sets of $\langle X,\tau,\leq\rangle$ that contains $X$,  
	\labitem{(Pr5)}{itm:Pr5} the set $$\Xb\coloneqq \{x\in X:\{U\in B:x\notin U\} \text{ is non-empty and
	up-directed}\}$$
	 is dense in $\langle X,\tau\rangle$.
\end{itemize}
We say that a clopen up-set $U$ is \emph{$\Xb$-admissible} if $\max(U^{c})\subseteq \Xb$. Thus $X$ is $\Xb$-admissible. We will denote by  $B^{\cap}$  the closure of $B$ under the binary operation of intersection.
\end{defn}

\begin{rem}\label{rem:up-directedfamiles}
Note that if $\spa{X}=\langle X,\tau,\BB\rangle$ is an \Sp space, then for every $x \in X$, $\{U\in B:x\notin U\}$  is non-empty and up-directed if and only if 
$\{U\in B^{\cap}:x\notin U\}$ is non-empty and up-directed. Therefore $$\Xb= \{x\in X:\{U\in B^\cap:x\notin U\} \text{ is non-empty and up-directed}\}.$$
\end{rem}

\begin{rem}
Note that the empty topological space together with the trivial algebra with domain $\{\emptyset\}$ is an $\Sm$-Priestley space.
\end{rem}

\begin{prop}\label{cor:Salg-Prspa}
For every filter-distributive and finitary congruential logic $\Sm$  with theorems  and every  $\Sm$-algebra $\Al$,  the structure $\sOptS(\Al)\coloneqq \langle \OptSA,\tau_{\Al},\varphi[\Al]\rangle$ is an \mbox{\Sp} space.
\end{prop}
\begin{proof}
Condition \ref{itm:Pr1} was proved in Theorem \ref{thm:ralg}. 
 Condition \ref{itm:Pr2} follows from Theorem \ref{thm:repr}. 
Condition \ref{itm:Pr3} is part of  (1) in  Corollary \ref{cor:isom1}. 
Condition \ref{itm:Pr4} follows from the definition of $\tau_{\Al}$. Condition  \ref{itm:Pr5}  follows from 
 Corollary \ref{cor:isom1}. 
\end{proof}

The dual space of an  $\Sm$-algebra $\Al$ will be the \mbox{\Sp} space $\sOptS(\Al)$ and the dual algebra of an  \mbox{\Sp} space $\spa{X}=\langle X,\tau,\BB\rangle$  will be $\BB$.

An enlightening  description of the open up-sets and the clopen-upsets is given next.

\begin{prop}\label{prop:Pr6-2}
Let $\langle X,\tau,\BB\rangle$ be an \Sp space and let $U\subseteq X$.
\begin{itemize}
	\item[(1)] $U$ is a non-empty open up-set of $\langle X,\tau,\leq\rangle$ if and only if $U$
	is the union of a non-empty set of non-empty sets which are  intersections of non-empty finite subsets  of $B$. 
	\item[(2)] $U$ is a non-empty clopen up-set of $\langle X,\tau,\leq\rangle$ if and only if $U$ is  the union of a non-empty finite set   of non-empty sets which are intersections of non-empty finite subsets  of $B$. 
	\end{itemize}
\end{prop}
\begin{proof}
Let $U$ be a non-empty open up-set of $\langle X,\tau,\leq\rangle$. 
When $U=X$ there is nothing to prove, 
so assume that $U\neq X$ and that 
$x\in U$. 
Because $U$ is an up-set, we have that for all $y\notin U$, $x\nleq y$. 
Then by \ref{itm:Pr1}, since the $\Sm$-referential algebra $\langle X,\BB\rangle$ is reduced, 
for all $y\notin U$ there is $V_{y}^{x}\in B$ such that $x\in V_{y}^{x}$ and $y\notin V_{y}^{x}$. 
Then we have a closed set $U^{c}$ and open sets $\{(V_{y}^{x})^{c}:y\notin U\}$ such that 
$U^{c}\subseteq \bigcup\{(V_{y}^{x})^{c}:y\notin U\}$. 
Now by the compactness of the space given by \ref{itm:Pr3}, there are $y_{0},\dots,y_{n_x}\notin U$ such that $U^{c}\subseteq (V_{y_{0}}^{x})^{c}\cup\dots\cup(V_{y_{n_x}}^{x})^{c}$. 
Hence $V_{y_{0}}^{x}\cap\dots\cap V_{y_{n_x}}^{x}\subseteq U$. 
Notice that  $x\in V_{y_{0}}^{x}\cap\dots\cap V_{y_{n_x}}^{x}$. Therefore 
$U\subseteq \bigcup_{x\in U} (V_{y_{0}}^{x}\cap\dots\cap V_{y_{n_x}}^{x})\subseteq U. $
Thus, as $U \not = \emptyset$,  $U$ is the union  of  a non-empty set of  non-empty sets which are     intersections  of finite and non-empty subsets of  $B$, and (1) has been proven, since the other direction is clear. 
(2) follows easily from the compactness of the space.
\end{proof}

Let $\spa{X}=\langle X,\tau,\BB\rangle$ be an \Sp space. 
The map $\xis:X\longrightarrow \ups(B)$ is defined as follows: 
\begin{equation*}
\xis(x)\coloneqq\{U\in B:x\in U\}. 
\end{equation*}
We will show later that  $\xi$ establishes a homeomorphism between    the $\Sm$-Priestsley space $\spa{X}$ and    the dual space  of its dual algebra. In this way we will have   the natural transformation we need to establish the  categorical duality.   

\begin{lemma}\label{lemma:xis1-1}\label{prop:xisonto-0}
Then map $\xis:X\longrightarrow \ups(B)$ is injective and for every $x \in X$,  $\xis(x)\in\FiSB$.
\end{lemma}
\begin{proof}
The injectivity of $\xis$  follows easily from \ref{itm:Pr1} because in reduced referential algebras the elements of the algebra separate points.  Moreover,  since   $\langle X, \BB\rangle$ is a reduced $\Sm$-referential algebra, for every $x \in X$ the set $\{U \in B: x \in U\}$ is an $\Sm$-filter.
\end{proof}

\begin{rem}\label{rem:Pr2}
Notice that Lemma \ref{lemma:xis1-1}  implies that the converse of condition \ref{itm:Pr2} holds for any \Sp space $\langle X,\tau,\BB\rangle$, 
\ie for all non-empty and finite  $\mathcal{V}\subseteq B$ and all $U\in B$, if $U\in\Consb(\mathcal{V})$, then $\bigcap \mathcal{V}\subseteq U$. 
Assume $U\in\Consb(\mathcal{V})$ and let $x\in\bigcap \mathcal{V}$, so $\mathcal{V}\subseteq\xis(x)$. 
Since $\xis(x)$ is an $\Sm$-filter, $\Consb(\mathcal{V})\subseteq \xis(x)$, and therefore, by  the assumption, $U\in \xis(x)$, \ie $x\in U$. 
\end{rem}

\begin{rem}
\label{rem:Pr2-2}
From Remark \ref{rem:Pr2} and condition \ref{itm:Pr2}   we obtain that in every \Sp space $\langle X,\tau,\BB\rangle$ for
 every finite  $\mathcal{U}, \mathcal{V}\subseteq B$, if $\mathcal{U}$ is non-empty, then
\begin{equation*}
\bigcap\mathcal{U} \subseteq \bigcap\mathcal{V} \IFF \Consb(\mathcal{V} ) \subseteq \Consb(\mathcal{U}). 
\end{equation*}In particular, for any $U,V\in B$:
\begin{equation*}
V\subseteq U \IFF U\in\Consb(V). 
\end{equation*}
Therefore, the specialization order $\leqsb$ on $B$ coincides with the inclusion relation on $B$. 
We will repeatedly use this fact  as well as the next generalization.
\end{rem}

\begin{lemma}\label{lemma:pinclu}
Let $\langle X,\tau,\BB\rangle$ be an \Sp space. 
For every non-empty and finite $\mathcal{U}_0, \ldots, \mathcal{U}_n, \mathcal{V}\subseteq B$:
\begin{equation*}
\bigcap \mathcal{V}\subseteq \bigcup_{i \leq n} \bigcap\mathcal{U}_i \IFF \bigcap_{i \leq n}\Consb(\mathcal{U}_i )\subseteq \Consb(\mathcal{V}) \IFF 
\varphim_{\BB}(\mathcal{V}) \subseteq \bigcup_{i \leq n} \varphim_{\BB}(\mathcal{U}_i)
\end{equation*}
where $\varphi_{\BB}$ is the representation map that corresponds to the optimal $\Sm$-filters of $\BB$.
\end{lemma}
\begin{proof}
We proof the first equivalence.  Assume first that $\bigcap \mathcal{V}\subseteq \bigcup_{i \leq n} \bigcap\mathcal{U}_i$ and assume that  $U\in \bigcap_{i \leq n}\Consb(\mathcal{U}_i )$. 
By Remark \ref{rem:Pr2} we get $ \bigcup_{i \leq n} \bigcap\mathcal{U}_i\subseteq U$; therefore $\bigcap \mathcal{V}\subseteq U$. 
It follows from  \ref{itm:Pr2} that $U\in\Consb(\mathcal{V})$. 
Assume now that $\bigcap_{i \leq n}\Consb(\mathcal{U}_i )\subseteq \Consb(\mathcal{V})$. 
We show that ${\bigcap \mathcal{V}\cap\Xb}\subseteq \bigcup_{i \leq n} \bigcap\mathcal{U}_i$, and then the claim follows from the denseness of $\Xb$ and from $\bigcup_{i \leq n} \bigcap\mathcal{U}_i$ being clopen. 
Let ${x\in {{\bigcap \mathcal{V}}\cap{\Xb}}}$ and suppose, towards a contradiction, that $x\notin \bigcup_{i \leq n} \bigcap\mathcal{U}_i$. Then for every $i \leq n$ there exists $U_i \in \mathcal{U}_i$ such that $x \not \in U_i$. 
Then using condition \ref{itm:Pr5}, there is $U\in B$ such that $\bigcup_{i \leq n}U_i\subseteq U$ and $x\notin U$. Thus for every $i \leq n$, $\bigcap \mathcal{U}_i \subseteq U_i \subseteq U$. Therefore, by \ref{itm:Pr2}, $U \in \bigcap_{i \leq n}\Consb(\mathcal{U}_i )$. Hence $U \in  \Consb(\mathcal{V})$ and this by  Remark \ref{rem:Pr2} implies $\bigcap \mathcal{V} \subseteq U$. As $x\in\bigcap \mathcal{V}$, we get $x\in U$, a contradiction.    

To prove the second equivalence, assume first that  $\bigcap\{\Consb(\mathcal{U}_{i}):i\leq n\}\subseteq \Consb(\mathcal{V}) $ and let $G\in\varphim_{\BB}(\mathcal{V})$. 
Then we have $\mathcal{V}\subseteq G$, and so $\Consb(\mathcal{V})\subseteq G$. 
Suppose, towards a contradiction, that $G\notin\bigcup\{\varphim_{\BB}(\mathcal{U}_{i}):i\leq n\}$. 
Then for each $i\leq n$ there is $U_{i}\in \mathcal{U}_{i}$ such that $U_{i}\notin G$. 
Notice that then $\bigcap\{\Consb(U_{i}):i\leq n\}\subseteq \bigcap\{\Consb(\mathcal{U}_{i}):i\leq n\} \subseteq  \Consb(\mathcal{V})$. As $G$ is an optimal $\Sm$-filter,
by Theorem \ref{thm:OptS-IdsS} we know that $G^{c}$ is a strong $\Sm$-ideal such that $\{U_{i}:i\leq n\}\subseteq G^{c}$, 
 and then we obtain $\Consb(\mathcal{V})\cap G^{c}\neq \emptyset$, a contradiction. 
For the converse, assume  $ \varphim_{\BB}(\mathcal{V})\subseteq \bigcup\{\varphim_{\BB}(\mathcal{U}_{i}):i\leq n\} $ and let $U\in \bigcap\{\Consb(\mathcal{U}_{i}):i\leq n\}$. 
Thus $\varphim_{\BB}(\mathcal{U}_{i})\subseteq \varphi_{\BB}(U)$ for all $i\leq n$, 
and then by  the assumption and Theorem \ref{thm:repr} we obtain $U\in\Consb(\mathcal{V})$.  
\end{proof}

%

We can give characterizations of when $\BB$ has a bottom element and of when it has a bottom-family.

\begin{lemma} \label{lemma:bottomb} 
Let $\langle X,\tau,\BB\rangle$ be an \Sp space. 
\begin{itemize}
	\item[(1)]  $\BB$ has a bottom element if and only if  $\emptyset\in B$. 
	\item[(2)] $\BB$ has a bottom-family if and only if there is a  finite $\mathcal{D}\subseteq B$  such that $\bigcap \mathcal{D}=\emptyset$. 
\end{itemize}
 \end{lemma}
 \begin{proof}
(1)  If $\emptyset\in B$, then $\emptyset$ is the bottom element of $\BB$. 
For the converse, assume that $U$ is the bottom element of $\BB$
and suppose, towards a contradiction, that there is $x\in U\cap\Xb$. 
Then by condition \ref{itm:Pr5}, there is $V\in B$ such that $x\notin V$, 
but since $U$ is the bottom element, then $U\subseteq V$. 
This implies $x\in V$, a contradiction. 
Hence, ${U\cap \Xb}={\emptyset}$, and then from the denseness of $\Xb$, $U=\emptyset$.  

For (2), note first that if $X = \emptyset$, then $B = \{\emptyset\}$ and so $\BB$ has a bottom-family and moreover $\bigcap \{\emptyset\} = \emptyset$. Assume then that $X$ is non-empty.  Suppose first that $\mathcal{D}\subseteq B$ is   finite  and $\bigcap \mathcal{D}=\emptyset$. Note that $\mathcal{D}$ must be non-empty, otherwise $\bigcap \mathcal{D} = X$. 
Then for every $U\in B$, $\bigcap \mathcal{D}=\emptyset\subseteq U$, and by \ref{itm:Pr2} it follows that $U\in\Consb(\mathcal{D})$, so $\Consb(\mathcal{D})=B$ and $\mathcal{D}$ is a bottom-family.  
For the converse, assume that $\BB$ has a bottom-family $\mathcal{D}\subseteq B$, and suppose, towards a contradiction, that there is  $x\in \bigcap \mathcal{D}$. 
Since $\mathcal{D}$ is finite, by denseness of $\Xb$ we can assume, without loss of generality, that $x\in\Xb$. 
Then by condition \ref{itm:Pr5} there is $V\in B$ such that $x\notin V$.  
But being $\mathcal{D}$ a bottom-family, $V \in \Consb(\mathcal{D})$. Remark \ref{rem:Pr2} then implies $\bigcap\mathcal{D}\subseteq V$; hence, $x \in V$,   a contradiction.
\end{proof}

\begin{rem}
Note that from the proof of Lemma \ref{lemma:bottomb} it follows that if $\mathcal{D}$ is a bottom-family of $\BB$, then $\bigcap\mathcal{D} = \emptyset$.
\end{rem}

\begin{prop} \label{prop:xisonto-a}
Let $\langle X,\tau,\BB\rangle$ be an \Sp space. 
Then for every ${x\in X}$, $\xis(x)\in\OptSB$. 
Moreover, if $x\in \Xb$, then $\xis(x)\in\IrrSB$. 
\end{prop}
\begin{proof}
Notice that if $\xis(x)=B$, then $\BB$ has no bottom-family. 
On the contrary, using Lemma \ref{lemma:bottomb},  
there would be $\mathcal{D}\subseteq B$ finite   such that $\bigcap \mathcal{D}=\emptyset$, but this contradicts  the assumption, which implies $x\in\bigcap \mathcal{D}$. 
Thus if $\xis(x)=B$, then $\xis(x)$ is an optimal $\Sm$-filter of $\BB$. 
 
Suppose now  that $\xis(x)\neq B$. To prove that $\xis(x)^{c} = \{U\in B:x\notin U\}$ is a strong $\Sm$-ideal we  use Lemma \ref{lem:useful-proof-strongideal}. 
Since $\leqsb$ is the inclusion relation, it is clear that  $\xis(x)^{c}$ is a downset w.r.t.\ $\leqsb$. 
Now we show 
that the condition on the definition of strong $\Sm$-ideal is satisfied. Let $\mathcal{V}\subseteq B$ be finite and let $I \subseteq \xis(x)^{c}$ be finite and such that $\bigcap \{\Consb(U): U \in I\}\subseteq \Consb(\mathcal{V})$.  If $I = \emptyset$,  the hypothesis turns into $\Consb(\mathcal{V})=B$. This implies that there is $\mathcal{V}'\subseteq \mathcal{V}$ that is a bottom-family for $\BB$. 
Thus, reasoning as in the proof of Lemma \ref{lemma:bottomb}, $\bigcap \mathcal{V'}=\emptyset$, and then there is $V\in\mathcal{V}$ such that  $x\notin V$, \ie $V\in\xis(x)^{c}$. If $I \not = \emptyset$ and $\mathcal{V} = \emptyset$, then since $\Sm$ is assumed to have theorems, $\Consb(\mathcal{V})=\Consb(\{X\}) = \{X\}$, because $1^{B} = X$. Hence,  by Lemma \ref{lemma:pinclu} we get $X \subseteq \bigcup \{U: U \in I\}$, and therefore $x \in U$ for some $U \in I$, a contradiction. Thus,  $\mathcal{V}$ is non-empty.  By Lemma \ref{lemma:pinclu} we get $\bigcap \mathcal{V}\subseteq \bigcup \{U: U \in I\}$, and then from $x\notin U$ for all $U \in I$ we get that there is $V\in\mathcal{V}$ such that $x\notin V$, \ie $V\in\xis(x)^{c}$. From either case we get that $\mathcal{V}\cap \xis(x)^{c}\neq\emptyset$, and so $\Consb(\mathcal{V})\cap \xis(x)^{c}\neq\emptyset$. 
Thus, we have shown that $\xis(x)^{c}$ is a strong $\Sm$-ideal, and by Theorem \ref{thm:OptS-IdsS} we conclude that $\xis(x)$ is an optimal $\Sm$-filter. 

Finally, if $x\in\Xb$, we know from \ref{itm:Pr5} that $\xis(x)^{c}$ is non-empty and up-directed, \ie an order ideal of $\BB$. 
Hence, by Theorem \ref{thm:IrrS-Id}, $\xis(x)$ is an irreducible $\Sm$-filter. 
\end{proof}

We aim to show that the map $\xis$ is onto $\OptSB$. To this end
 we note that $\BB^{\cap}\coloneqq\langle B^{\cap},\cap,X\rangle$ is a meet-semilattice isomorphic to  the $\Sm$-semilattice $\Mee(\BB)$ of $\BB$ by the map $h$ given by the rule: 
\begin{align*}
\bigcap \mathcal{U}&\longmapsto \bigcap\{\varphi_{\BB}(U):U\in\mathcal{U}\}=\varphim_{\BB}(\mathcal{U}),
\end{align*}
where $\mathcal{U}$ is a finite subset of $B$. From the condition \ref{itm:Pr2} and Lemma \ref{lemma:pinclu} it follows that the map is well defined and one-to-one. It is obviously onto. Finally, it is easy to see that it is a homomorphism preserving the top element. From the distributivity of $\Mee(\BB)$, which follows from the fact that the logic $\Sm$ is distributive, we obtain that $\BB^{\cap}\coloneqq\langle B^{\cap},\cap,X\rangle$ is a distributive meet-semilattice. 

Using the isomorphism $h: \BB^{\cap} \cong \Mee(\BB)$ and the isomorphisms $$\Fimg{\varphi[\cdot]}:\FiSB\,\,\cong\,\, \Fim(\MeeB):\varphi^{-1}[\cdot].$$ that restrict to the optimals as $$\Fimg{\varphi[\cdot]}:\langle\OptSB,\subseteq\rangle\,\,\cong\,\, \langle\Optm(\MeeB),\subseteq\rangle:\varphi^{-1}[\cdot]$$
it is easy to see that for every $\Sm$-filter $F$ of  $\BB$, $F \in \OptSB$ if and only if $\{E \in B^{\cap}: \exists\ U_0, \ldots, U_n \in F \text{ s.t. }  U_0 \cap \ldots \cap U_n \subseteq E\} \in \Optm(\BB^{\cap})$. 

Let $B^{\cap\cup}$ be the closure of $B^{\cap}$ under the binary operation of union, so that $\emptyset \in B^{\cap}$ if and only if $\emptyset \in B^{\cap\cup}$. We prove that the identity embedding from $\BB^{\cap}$ to the distributive lattice with top  $\BB^{\cap\cup}\coloneqq\langle B^{\cap\cup},\cap,\cup,X\rangle$ is a sup-homomorphism. Then, since the set $B^{\cap}$ is obviously join-dense in $\BB^{\cap\cup}$, it follows that $\BB^{\cap\cup}$ is the distributive  envelope of $\BB^{\cap}$. 

\begin{prop}
\label{prop:id-sup-hom-B}
The identity map from $\BB^{\cap}$ to $\BB^{\cap\cup}$ is a sup-homomorphism, that is, for all non-empty and finite subsets $\mathcal{U}_0, \ldots, \mathcal{U}_n, \mathcal{V}$ of $B$
\[
\text{ if } \bigcap_{i \leq n}{\uparrow}_{\BB^{\cap}}\bigcap\mathcal{U}_i\subseteq {\uparrow}_{\BB^{\cap}} \bigcap\mathcal{V}, \text{ then } \bigcap_{i \leq n}{\uparrow}_{\BB^{\cap\cup}} \bigcap\mathcal{U}_i \subseteq {\uparrow}_{\BB^{\cap\cup}} \bigcap\mathcal{V}
\]
\end{prop}
\begin{proof}
Suppose that $\bigcap_{i \leq n}{\uparrow}_{\BB^{\cap}} \bigcap\mathcal{U}_i \subseteq {\uparrow}_{\BB^{\cap}} \bigcap\mathcal{V}$ and let $C \in \bigcap_{i \leq n}{\uparrow}_{\BB^{\cap\cup}} \bigcap\mathcal{U}_i$. $C = \bigcup_{j \leq m}\bigcap\mathcal{W}_j$ for some non-empty finite subsets $\mathcal{W}_0, \ldots, \mathcal{W}_m$ of $B$. 
We prove that 
$\bigcap\mathcal{V} \cap X_B \subseteq C$. Then from denseness it   will follow that $\bigcap\mathcal{V}  \subseteq C$. Suppose that  $x \in \bigcap\mathcal{V} \cap X_B$ and $x \not \in C$. Then $x \not \in \bigcup_{j \leq m}\bigcap\mathcal{W}_j$. Let for every $j \leq m$, $W_j \in \mathcal{W}_j$ such that $x \not \in W_j$. Being $x \in X_B$, the set $\{U \in B: x \not \in U\}$ is non-empty and up-directed. Thus there exists $W \in B$ such that $W_0 \cup \ldots \cup W_m \subseteq W$ and $x \not \in W$. It follows that  $\bigcup_{j \leq m}\bigcap\mathcal{W}_j \subseteq W$, \ie $C \subseteq W$; therefore $W \in \bigcap_{i \leq n}{\uparrow}_{\BB^{\cap\cup}} \bigcap\mathcal{U}_i$. The assumption implies that 
$ \bigcap\mathcal{V} \subseteq W$. As $x \in  \bigcap\mathcal{V}$, $x \in W$, a contradiction.
\end{proof}

From the fact that $\BB^{\cap}$ and  $\Mee(\BB)$ are isomorphic it follows that $\BB^{\cap\cup}$ is (isomorphic to) the distributive envelope $\Lee(\Mee(\BB))$ of $\Mee(\BB)$.  We describe the isomorphism in the next proposition.

\begin{prop}
\label{prop:isomBcapcup-distr-envelopM(B))}
The map $g: B^{\cap\cup} \to \Lee(\Mee(\BB))$  defined by
\begin{align*}
g(\bigcup_{i \leq n}\bigcap \mathcal{U}_{i}) : =  \bigcup\{\sigma_{\Mee(\BB)}(\varphim_{\BB}(\mathcal{U}_{i})):i\leq n\}, 
\end{align*}
where the $\mathcal{U}_i$'s are non-empty and finite subsets of $B$,
is a lattice isomorphism between $\BB^{\cap\cup}\coloneqq\langle B^{\cap\cup},\cap,\cup,X\rangle$  and the distributive envelope $\Lee(\Mee(\BB))$ of $\Mee(\BB)$.\end{prop}
\begin{proof}
First of all we need to see that $g$ is well defined. To this end we prove that if $\{\mathcal{U}_i:i \leq n\}$ and $\{\mathcal{V}_j:j \leq m\}$ are  finite families of non-empty finite subsets of $B$ such that 
$\bigcup_{i \leq n}\bigcap\mathcal{U}_i \subseteq \bigcup_{j \leq m}\bigcap\mathcal{V}_j$, then $\bigcup_{i \leq n}\sigma_{\Mee(\BB)}(\varphim_{\BB}(\mathcal{U}_{i})) \subseteq \bigcup_{j \leq m}\sigma_{\Mee(\BB)}(\varphim_{\BB}(\mathcal{V}_{j}))$. From this fact  also follows that the map is order preserving.  Suppose that   $\bigcup_{i \leq n}\bigcap\mathcal{U}_i \subseteq \bigcup_{j \leq m}\bigcap\mathcal{V}_j$. From Lemma \ref{lemma:pinclu} it follows that 
$\bigcup_{i \leq n} \varphim_{\BB}(\mathcal{U}_{i}) \subseteq \bigcup_{j \leq m} \varphim_{\BB}(\mathcal{V}_{j})$.
Then   $\bigcap_{j \leq m} {\uparrow} \varphim_{\BB}(\mathcal{V}_{j}) \subseteq {\uparrow}  \varphim_{\BB}(\mathcal{U}_{i})$ for every $i \leq n$. Using condition (\ref{eq:distr-envelop-1})
 we obtain  that $\bigcup_{i \leq n}\sigma_{\Mee(\BB)}(\varphim_{\BB}(\mathcal{U}_{i})) \subseteq \bigcup_{j \leq m}\sigma_{\Mee(\BB)}(\varphim_{\BB}(\mathcal{V}_{j}))$. 
 
 The next move is to show that $g$ is  order reflecting. Suppose that   $\{\mathcal{U}_i:i \leq n\}$ and $\{\mathcal{V}_j:j \leq m\}$ are  finite families of non-empty finite subsets of $B$ such that  $\bigcup_{i \leq n}\sigma_{\Mee(\BB)}(\varphim_{\BB}(\mathcal{U}_{i})) \subseteq \bigcup_{j \leq m}\sigma_{\Mee(\BB)}(\varphim_{\BB}(\mathcal{V}_{j}))$. Condition  (\ref{eq:distr-envelop-1}) implies that $\bigcup_{i \leq n} \varphim_{\BB}(\mathcal{U}_{i}) \subseteq \bigcup_{j \leq m} \varphim_{\BB}(\mathcal{V}_{j})$. Applying
 Lemma \ref{lemma:pinclu}  we obtain that $\bigcup_{i \leq n}\bigcap\mathcal{U}_i \subseteq \bigcup_{j \leq m}\bigcap\mathcal{V}_j$. 
 \end{proof}

Notice that the emptyset $\emptyset\in\clopu(X)$ can be trivially described as an (empty) finite union of non-empty finite intersections of elements of $B$. 
Therefore, 
the previous proposition implies  that for every \Sp space $\langle X,\tau,\BB\rangle$ 
the lattice of clopen up-sets $\aclopu(X)$ (which is
 the dual distributive lattice of the Priestley space $\langle X,\tau,\leq\rangle$) is isomorphic to $\Lee(\Mee(\BB))\cup\{\emptyset\}$,
  which is the distributive envelope of $\Mee(\BB)$ augmented with a bottom element whenever $\Mee(\BB)$ has no bottom element. 
 In particular, this implies that the optimal filters of $\Lee(\Mee(\BB))$ are in one-to-one correspondence with the prime filters of $\aclopu(X)$. 
  This fact will be used in the proof of the following proposition. 
%

\begin{prop} \label{prop:xisonto}
Let $\langle X,\tau,\BB\rangle$ be an \Sp space. 
Then: 
\begin{itemize}
	\item[(1)] For every $P\in\OptSB$ there is $x\in X$ such that $\xis(x)=P$.
	\item[(2)] For every  $Q\in\IrrSB$ there is $x\in \Xb$ such that $\xis(x)=Q$.
\end{itemize}
\end{prop}
\begin{proof}
(1)\, \,Let $P$ be an optimal $\Sm$-filter of $\BB$. 
Then we know by Proposition \ref{prop:mefilter} that $\Fimg{\varphi[P]}$ is an optimal filter of $\Mee(\BB)$. Therefore, $\Fimg{\sigma[\Fimg{\varphi[P]}]}$ is a prime  filter of $\Lee(\Mee(\BB))$, if  $\Fimg{\varphi[P]} \not = \Mee(\BB)$ or is $\Lee(\Mee(\BB))$ in case $\Fimg{\varphi[P]} = \Mee(\BB)$.  
By the isomorphism $g$ in Proposition \ref{prop:isomBcapcup-distr-envelopM(B))}, $g^{-1}[\Fimg{\sigma[\Fimg{\varphi[P]}]}$ is a prime filter of $\BB^{\cap\cup}$ or $g^{-1}[\Fimg{\sigma[\Fimg{\varphi[P]}]} = B^{\cap\cup}$. In any case it is a prime filter of $\aclopu(X)$. By Priestley duality there exists $x \in X$ such that 
$g^{-1}[\Fimg{\sigma[\Fimg{\varphi[P]}]} = \{C \in \clopu(X): x \in C\}$. Now it is easy to see that 
$B \cap g^{-1}[\Fimg{\sigma[\Fimg{\varphi[P]}]} = P$. It follows that $P = \{U \in B: x \in U\}$. Hence, $P = \xis(x)$.

(2)\,\, Let $Q\in\IrrSB$. By (1) we know that there is $x\in X$ such that $Q=\{U\in B:x\in U\}$. 
Moreover by Theorem \ref{thm:IrrS-Id} we know that $Q^{c}=\{U\in B:x\notin U\}$ is and order ideal, so it is non-empty and up-directed. Hence by \ref{itm:Pr5} we conclude that $x\in \Xb$. 
\end{proof}

\begin{cor}\label{cor:xisok}
Let  $\langle X,\tau,\BB\rangle$ be an \Sp space. The map $\xis$ is an order homeomorphism between the ordered topological spaces $\langle X,\tau,\leq\rangle$ and $\langle \OptSB,\tau_{\BB},\subseteq\rangle$ such that for every $U\in B$, $\xis^{-1}[\varphib(U)]=U$ and $\ \xis[U]=\varphib(U)$.
\end{cor}
\begin{proof}
Notice that for all $x\in X$ and all $U\in B$ we have: $x\in U$ if and only if $U\in\xis(x)$ if and only if $\xis(x)\in\varphib(U)$. 
Thus, $\xis^{-1}[\varphib(U)]=U$ and moreover:
\begin{equation*}
x\in \xis^{-1}[\varphib(U)^{c}] \IFF \xis(x)\in\varphib(U)^{c} \IFF U\notin \xis(x) \IFF x\in U^{c}.
\end{equation*}
Therefore $\xis^{-1}[\varphib(U)^{c}]=U^{c}$.  
Since inverse maps preserve intersections, 
this implies that the inverse of $\xis$ sends basic opens to basic opens. 
From condition \ref{itm:Pr1} it follows that $\xis$ is order preserving.
As $\xis$ is one-to-one (Lemma \ref{lemma:xis1-1}), onto (Proposition \ref{prop:xisonto}), and its inverse sends basic opens of $\langle \OptSB,\tau_{\BB}\rangle$ to basic opens of $\langle X,\tau\rangle$, 
we conclude that $\xis$ is an homeomorphism, as required.
\end{proof}

\begin{cor}\label{cor:Prspa-Prspa}
Let $\spa{X}=\langle X,\tau,\BB\rangle$ be an \Sp space. Then the structure $\langle \OptSB,\tau_{\BB},\varphi_{\BB}[\BB]\rangle$ 
is an \Sp space such that $\langle X,\tau\rangle$ and $\langle\OptSB,
\tau_{\BB}\rangle$ are homeomorphic topological spaces by means of the map $\xi:X\longrightarrow \OptSB$,
 that moreover is an order isomorphism between $\langle X,\leq\rangle$ and $\langle \OptSB,\subseteq\rangle$.  
 Furthermore $\BB$ and $\varphi_{\BB}[\BB]$ are isomorphic $\Sm$-algebras by means of the map $\varphi_{\BB}:B\longrightarrow \varphi_{\BB}[B]$ whose inverse if the map $\xis^{-1}[\cdot]:\varphi_{\BB}[B] \to B$. 
\end{cor}


The previous corollary establishes that $\Sm$-algebras and \Sp spaces are equivalent objects by means of the maps $\varphi$ and $\xis$. 

Before dealing with morphisms, let us investigate some other properties of \Sp spaces.
Notice that from Proposition \ref{prop:Pr6-2} we know that for any $\Sm$-algebra, 
the collection of clopen up-sets of the Priestley space $\langle \OptSA,\taua,\subseteq\rangle$ is $\varphi[A]^{\cap\cup}\cup\{\emptyset\}$. 
We show in the following proposition that within all clopen up-sets of $\OptSA$, 
those that are $\IrrSA$-admissible admit a simpler description as finite intersections of elements of $\varphi[A]$, \ie as elements of $\Mee(\Al)$.

\begin{prop}\label{prop:a3}
Let $\langle X,\tau,\BB\rangle$ be an \Sp space. 
Then 
\begin{enumerate}
\item the collection of $\Xb$-admissible clopen up-sets of $X$ coincides with $B^{\cap}$,
\item  the tuple $\langle X,\tau,\leq,\Xb\rangle$ is a generalized Priestley space and its dual meet-semilattice, i.e., the meet-semilattice of all $\Xb$-admissible clopen up-sets,  is isomorphic to $\Mee(\BB)$. 
\end{enumerate}
\end{prop}
\begin{proof}
(1) We claim that for any finite subsets $\{\mathcal{V}_{i}\subseteq B:i\leq n\}$ for some $n\in\omega$, we have 
$$x\in\max( (\bigcup_{i\leq n}\bigcap\mathcal{V}_{i})^{c}) \IFF \xis(x)\in\max((\bigcup_{i\leq n}\varphim_{\BB}(\mathcal{V}_{i}))^{c}).$$
This follows easily from propositions \ref{prop:xisonto-a} and \ref{prop:xisonto} using that for any $V\in B$, $x\in V$ if and only if $\xis(x)\in\varphi_{\BB}(V)$.

Let first assume that  $U\in B^{\cap}$, \ie $U=\bigcap\mathcal{V}$ for some non-empty and finite $\mathcal{V}\subseteq B$. 
By Proposition \ref{prop:isom2} we know that $\varphim_{\BB}(\mathcal{V})$ is an $\IrrSB$-admissible clopen up-set of $\langle \OptSB,\tau_{\BB},\subseteq\rangle$.
Then from the claim and using that $\Xb=\xis^{-1}[\IrrSB]$, we obtain that $\bigcap\mathcal{V}$ is an $\Xb$-admissible clopen up-set of $X$. 

Suppose now that $U$ is an $\Xb$-admissible clopen up-set of $X$. 
If $U=\emptyset$, then from the claim we get that  $\emptyset$ is an $\IrrSB$-admissible clopen up-set of $\OptSB$. 
Then by Proposition \ref{prop:isom2}, there is a finite $\mathcal{W}\subseteq B$ such that $\emptyset=\varphim_{\BB}(\mathcal{W})$, and this implies $U=\emptyset=\bigcap\mathcal{W}$. 
Assume now that $U\neq \emptyset$.
 Then by Proposition \ref{prop:Pr6-2} we know that $U=\bigcup_{i\leq n}\bigcap\mathcal{V}_{i}$ for finite 
 $\mathcal{V}_{i}\subseteq B$ for all $i\leq n$ for some $n\in \omega$. 
Then from the assumption and the claim, 
using that $\xis[\Xb]=\IrrSB$, 
we obtain that $\bigcup_{i\leq n}\varphim_{\BB}(\mathcal{V}_{i})$ is an $\IrrSB$-admissible clopen up-set of $\OptSB$. 
By Proposition \ref{prop:isom2}, 
there is a finite $\mathcal{W}\subseteq B$ such that $\bigcup_{i\leq n}\varphim_{\BB}(\mathcal{V}_{i})=\varphim_{\BB}(\mathcal{W})$. 
And this implies $\bigcup_{i\leq n}\bigcap\mathcal{V}_{i}=\bigcap \mathcal{W}$, as required.

(2) The condition \ref{itm:PrDS1}  in the definition of generalized Priestley space, namely that  $\langle X,\tau,\leq\rangle$ is a Priestley space, follows from   \ref{itm:Pr1}, \ref{itm:Pr3} and \ref{itm:Pr4}.  The condition \ref{itm:PrDS2} follows from from \ref{itm:Pr5}.  
Recall that $X^{*}$ denotes the collection of $\Xb$-admissible clopen up-sets of $\langle X,\tau,\leq,\Xb\rangle$. By  (1), $X^{*}$ coincides with  $B^{\cap}$. 
Then from  \ref{itm:Pr5} and Remark \ref{rem:up-directedfamiles}  it follows that 
$\Xb=\{x\in X:\{U\in B^{\cap}:x\notin U\}\text{ is non-empty and up-directed}\}$, 
so condition \ref{itm:PrDS3} also holds.
Moreover, from \ref{itm:Pr1} it follows that 
for all $x,y\in X$, 
$x\leq y$ if and only if for all $U\in B^{\cap}$,
 if $x\in U$ then $y\in U$, so condition \ref{itm:PrDS4} also holds. 
\end{proof}



%
\begin{cor}\label{cor:a2}
Let $\Al$ be an $\Sm$-algebra. 
Then $\langle \OptSA,\taua,\leq,\IrrSA\rangle$ is a generalized Priestley space, whose dual meet-semilattice is isomorphic to $\Mee(\Al)$.
\end{cor}

Since for any Priestley space $X$ the collection $\{U\setminus V:U,V\in\clopu(X)\}$ is a basis for the space, 
from Proposition \ref{prop:Pr6-2} we obtain that for any \Sp space $X$ the collection $B\cup\{U^{c}:U\in B\}$ is a subbasis of the space $X$. 
The next proposition highlights that this issue   is strongly connected with the fact that the $\Sm$-referential algebra $\langle X,\BB\rangle$ is reduced, and
 leads us to an alternative definition of \Sp space.

\begin{prop}\label{prop:simpl1}
For any $\Sm$-referential algebra $\langle X,\BB\rangle$ augmented with a topology $\tau$ and an order $\leq$ on $X$,  
 if $\langle X,\tau,\leq\rangle$ is a Priestley space, $X\in B$ and $\clopu(X)=B^{\cap\cup}\cup\{\emptyset\}$, then $\langle X,\BB\rangle$ is reduced. 
\end{prop}
\begin{proof}
Let $\langle X,\tau,\leq\rangle$ be a Priestley space satisfying the conditions above mentioned. 
We show that $\langle X,\BB\rangle$ is reduced by showing that $\leq$ is the quasiorder $\preceq$ associated with the referential algebra. 

Let first $x,y\in X$ be such that $x\leq y$. 
As the elements of $B$ are up-sets, it follows that for all $V\in B$, if $x\in V$ then $y\in V$. 
Let now $x,y\in X$ be such that $x\nleq y$. 
Then by totally order disconnectedness of the space, there is $U$ a clopen up-set such that $x\in U$ and $y\notin U$. 
Clearly $U\neq\emptyset$, so by assumption $U\in B^{\cap\cup}$.   
Then there are non-empty and finite subsets $\mathcal{U}_{i}\subseteq B$, with $i\leq n$, for some $n\in \nat$, such that $x\in \bigcup\{\bigcap \mathcal{U}_{i}:i\leq n\}=U$ and $y\notin \bigcup\{\bigcap \mathcal{U}_{i}:i\leq n\}$. 
So there is $i\leq n$ such that $x\in \bigcap \mathcal{U}_{i}$ and $y\notin \bigcap\mathcal{U}_{i}$. 
And then there is $U_{i}\in\mathcal{U}_{i}\subseteq B$ such that $x\in U_{i}$ and $y\notin U_{i}$. 

We conclude that for all $x,y\in X$,  $x\leq y$ if and only for all $V\in B$, if $x\in V$ then $y\in V$. 
Hence ${\leq} ={\preceq}$. 
And since $\leq$  is a partial order, it follows that the referential algebra $\langle X,\BB\rangle$ is reduced. 
\end{proof}

\begin{cor}\label{cor:SstarPr}
A structure $\spa{X}=\langle X,\tau,\BB\rangle$ is an \Sp space if and only if the following conditions are satisfied:
\begin{itemize}
	\labitem{(Pr1$'$)}{itm:Pr1p} $\langle X,\BB\rangle$ is an $\Sm$-referential algebra, whose associated quasiorder is denoted by $\leq$,
	\item[\ref{itm:Pr2}]   for any finite $\mathcal{V}\subseteq B$ and any $U\in B$, if $\bigcap\mathcal{V}\subseteq U$, then $U\in \Consb(\mathcal{V})$,
	\labitem{(Pr3$'$)}{itm:Pr3p} $\langle X,\tau,\leq\rangle$ is a Priestley space, and $B\cup\{U^{c}:U\in B\}$ is a subbasis for it,
	\labitem{(Pr4$'$)}{itm:Pr4p} $X\in B$ and $\clopu(X)=B^{\cap\cup}\cup\{\emptyset\}$,
	\item[\ref{itm:Pr5}] the set $\Xb\coloneqq \{x\in X:\{U\in B:x\notin U\} \text{ is non-empty and up-directed}\}$ is dense in $\langle X,\tau\rangle$.
\end{itemize}
\end{cor}

\subsection{Morphisms between \Sp spaces}\label{sec:Morf}

The approach for this section is similar to that of  
\cite{BeJa11}. 
From now on, let $\Sm$ be a filter-distributive and  finitary congruential logic with theorems.

 Let $\Al_{1}$ and $\Al_{2}$ be  $\Sm$-algebras  
and let $h\in\Hom(\Al_{1},\Al_{2})$. The dual relation of $h$ is the relation   
 $R_{h}\subseteq \OptS(\Al_{2})\times \OptS(\Al_{1})$ defined by: 
\begin{equation*}
(P,Q)\in R_{h} \IFF h^{-1}[P]\subseteq Q.
\end{equation*}

\begin{prop}\label{prop:sRrepr1}
Let $\Al_{1}$, $\Al_{2}$ be  $\Sm$-algebras  and  $h\in\Hom(\Al_{1},\Al_{2})$. 
For all $a\in A_{1}$:
\begin{itemize}
	\item[(1)] $R_{h}^{-1}(\varphi_{1}(a)^{c})=\varphi_{2}(h(a))^{c}$.
	\item[(2)] $\Box_{R_{h}}(\varphi_{1}(a))=\varphi_{2}(h(a))$.
	\item[(3)] $\Box_{R_{h}}\in \Hom(\varphi_{1}[\Al_{1}],\varphi_{2}[\Al_{2}])$. 
\end{itemize}
\end{prop}
\begin{proof}
(1) First we show that $R_{h}^{-1}(\varphi_{1}(a)^{c})\subseteq\varphi_{2}(h(a))^{c}$. Let
 ${P\in \OptS(\Al_{2})}$ such that $P\in R_{h}^{-1}(\varphi_{1}(a)^{c})$, \ie $h^{-1}[P]\subseteq Q$ for some $Q\notin \varphi_{1}(a)$. 
Then from $a\notin Q$ we get $a\notin h^{-1}[P]$, hence $h(a)\notin P$ and so $P\in\varphi_{2}(h(a))^{c}$. 
For the converse, let $P\in \varphi_{2}(h(a))^{c}$, \ie $a\notin h^{-1}[P]$. 
As $P$ is an $\Sm$-filter of $\Al_{2}$, $h^{-1}[P]$ is an $\Sm$-filter of $\Al_{1}$. 
By the optimal $\Sm$-filter lemma, 
since ${\downarrow}a$ is a strong $\Sm$-ideal and ${\downarrow}a\cap h^{-1}[P]=\emptyset$,
 there is $Q \in \OptS(\Al_{1})$  such that $a\notin Q\supseteq h^{-1}[P]$.  
So, $Q\in \varphi_{1}(a)^{c}$ and $Q\in R_{h}(P)$. Hence $P\in R_{h}^{-1}(\varphi_{1}(a)^{c})$. 

 (2) First we show that $\Box_{R_{h}}(\varphi_{1}(a))\subseteq \varphi_{2}(h(a))$, so let $P\in \Box_{R_{h}}(\varphi_{1}(a))$, \ie $R_{h}(P)\subseteq \varphi_{1}(a)$. 
Suppose, towards a contradiction, that $P\notin \varphi_{2}(h(a))$. 
Then by item (1) we have $P\in R_{h}^{-1}(\varphi_{1}(a)^{c})$, so there is $Q\in R_{h}(P)$ such that $Q\notin \varphi_{1}(a)$, a contradiction. 
For the converse, let $P\in \varphi_{2}(h(a))$, so $a\in h^{-1}[P]$. 
Then for any $Q\in R_{h}(P)$, from $h^{-1}[P]\subseteq Q$  we get $a\in Q$, \ie $Q\in\varphi_{1}(a)$. 
This implies that $R_{h}(P)\subseteq \varphi_{1}(a)$, \ie $P\in \Box_{R_{h}}(\varphi_{1}(a))$, as required.

(3) Let $f$ be an $n$-ary connective of the language and let $a_{i}\in A_{1}$ for each $i\leq n$. 
Using the definition of $\varphi_{1}[\Al_{1}]$ and  $\varphi_{2}[\Al_{2}]$, 
item (2), 
and the fact that $h \in \Hom(\Al_{1}, \Al_{2})$, we get: 
\begin{align*} \Box_{R_{h}}(f^{\varphi_{1}[\Al_{1}]}(\varphi_{1}(a_{1}),\dots,\varphi_{1}(a_{n})))
&=\Box_{R_{h}}(\varphi_{1}(f^{\Al_{1}}(a_{1},\dots,a_{n})))\\
&=\varphi_{2}(h(f^{\Al_{1}}(a_{1},\dots,a_{n})))\\
&=\varphi_{2}(f^{\Al_{2}}(h(a_{1}),\dots,h(a_{n})))\\
&=f^{\varphi_{2}[\Al_{2}]}(\varphi_{2}(h(a_{1})),\dots,\varphi_{2}(h(a_{n})))\\
&=f^{\varphi_{2}[\Al_{2}]}(\Box_{R_{h}}(\varphi_{1}(a_{1})),\dots,\Box_{R_{h}}(\varphi_{1}(a_{n}))).
\end{align*}
\end{proof}

\begin{prop}\label{prop:sRPr}
Let $\Al_{1}$, $\Al_{2}$ be  $\Sm$-algebras  and $h\in\Hom(\Al_{1},\Al_{2})$. 
For any $P\in\OptS(\Al_{2})$ and $Q\in\OptS(\Al_{1})$ such that $(P,Q)\notin R_{h}$, there is $a\in A_{1}$ such that $Q\notin\varphi(a)$ and $R_{h}\subseteq \varphi(a)$.
\end{prop}
\begin{proof}
From $(P,Q)\notin R_{h}$ we get $h^{-1}[P]\nsubseteq Q$, so there is $a\in A$ such that $a\in h^{-1}[P]$ and $a\notin Q$. 
This implies that $Q\notin \varphi(a)$ and  for all $Q'\in\Opt_{\Sm}(\Al_{1})$ such that $(P,Q')\in R_{h}$, $a\in Q'$. 
Therefore $R_{h}(P)\subseteq \varphi(a)$ and we are done.
\end{proof}

Notice that the previous propositions hold in general for any finitary congruential logic with theorems, not necessarily a filter-distributive one. 
They lead us to the definition of the Priestley-dual morphisms   between $\Sm$-algebras. 

\begin{defn}\label{defn:Prmor}
Let $\spa{X}_{1}=\langle X_{1},\tau_{1},\BB_{1}\rangle$ and $\spa{X}_{2}=\langle X_{2},\tau_{2},\BB_{2}\rangle$ be two \Sp spaces. 
A relation $R\subseteq X_{1}\times X_{2}$ is an \emph{\Sp morphism} when: \begin{itemize}
	\labitem{(PrR1)}{itm:PrR1} $\Box_{R}\in\Hom(\BB_{2},\BB_{1})$,
	\labitem{(PrR2)}{itm:PrR2} if $(x,y)\notin R$, then there  is $U\in B_{2}$ such that $y\notin U$ and $R(x)\subseteq U$.
\end{itemize}
\end{defn}

\begin{prop}\label{cor:Shom-Prmor}  
Let $\Al_{1}$,  $\Al_{2}$ be $\Sm$-algebras and 
 $h\in\Hom(\Al_{1},\Al_{2})$. 
Then ${R}_{h}$ is an \Sp morphism between \Sp spaces $\sOptS(\Al_{2})$ and $\sOptS(\Al_{1})$. 
\end{prop}
\begin{proof}
\ref{itm:PrR1} follows from Proposition \ref{prop:sRrepr1} and \ref{itm:PrR2} follows from Proposition \ref{prop:sRPr}. 
\end{proof}

Recall that Proposition \ref{prop:a3} tells us  that for any   \Sp space $\langle X,\tau,\BB\rangle$, 
the structure $\langle X,\tau,\leq,\Xb \rangle$ is a generalized Priestley space. 
Analogously, in the next theorem we show how \Sp morphisms and generalized Priestley morphisms are related.

\begin{thm}\label{thm:gpm}
Let $R\subseteq X_{1}\times X_{2}$ be an \Sp morphism between \Sp spaces $\spa{X}_{1}$ and $\spa{X}_{2}$. 
Then $R$ is a generalized Priestley morphism between generalized Priestley spaces $\langle X_{1},\tau_{1},\leq_{1},\Xbone\rangle$ and $\langle X_{2},\tau_{2},\leq_{2},\Xbtwo\rangle$.
\end{thm}
\begin{proof}
We just need to check that for any $\Xbtwo$-admissible clopen up-set of $X_{2}$, we have  that $\Box_{R}(U)$ is an $\Xbone$-admissible clopen up-set of $X_{1}$.   
So let $U\in\clopu(X_{2})$ be such that $\max(U^{c})\subseteq \Xbtwo$. 
By Proposition \ref{prop:a3} there are $U_{0},\dots,U_{n}\in B_{2}$ such that $U=U_{0}\cap\dots\cap U_{n}$. 
Then we have that $\Box_{R}(U)=\{x\in X:R(x)\subseteq U_{0}\cap\dots\cap U_{n}\}=\Box_{R}(U_{0})\cap\dots\cap \Box_{R}(U_{n})$. 
And then by \ref{itm:PrR1} and Proposition \ref{prop:a3} again, 
$\max((\Box_{R}(U))^{c})\subseteq\Xbone$,  as required. 
\end{proof}

The order associated with the $\Sm$-referential algebra plays a prominent role in the duality. The next propositions show that it is an \Sp morphism and that its relational composition with any \Sp morphism $R$ is included in $R$.

\begin{prop}\label{prop:porder}
For any \Sp space $\spa{X}=\langle X,\tau,\BB\rangle$, 
the order associated with the $\Sm$-referential algebra $\langle X,\BB\rangle$ is an \Sp morphism.
\end{prop}
\begin{proof}
Recall that we denote the order associated with the $\Sm$-referential algebra $\langle X,\BB\rangle$ by $\leq$. 
As the referential algebra is reduced, for any $x,y\in X$ such that $x\nleq y$, there is $U\in B$ such that $x\in U$ and $y\notin U$. 
Moreover, as $B$ is a family of clopen up-sets, for every $z\in {\uparrow} x$ we get $z\in U$. 
Therefore ${\uparrow}x\subseteq U$, hence condition \ref{itm:PrR2} is satisfied by ${\leq}$. 
Notice also that $\Box_{\leq}(Y)=\{x\in X:{\uparrow}x\subseteq Y\}$. 
Since the elements of $B$ are up-sets with respect to $\leq$,  for all $U\in B$ we have
$\Box_{\leq}(U)=U$. 
Therefore $\Box_{\leq}$ is the identity map  from $B$ to $B$, 
and so $\Box_{\leq}\in\Hom(\BB,\BB)$ and condition \ref{itm:PrR1} is also satisfied by ${\leq}$. 
Hence the relation ${\leq} \subseteq {X\times X}$ is an \Sp morphism. 
\end{proof}

\begin{prop}
Let $\spa{X}_{1}=\langle X_{1},\tau_{1},\BB_{1}\rangle$ and $\spa{X}_{2}=\langle X_{2},\tau_{2},\BB_{2}\rangle$ be two \Sp spaces and 
 $R\subseteq X_{1}\times X_{2}$ an \emph{\Sp morphism}. Then ${\leq_1} \circ R \subseteq R$ and $R \circ {\leq_2} \subseteq R$. 
\end{prop}
\begin{proof}
If $x \leq_1 y$,  $(y, z) \in R$ and $(x, z) \not \in  R$, let  $U\in B_{2}$ such that $z\notin U$ and $R(x)\subseteq U$. Thus, $x \in \Box_R(U)$ and since $\Box_R(U) \in B_1$ it is an up-set; therefore $y \in \Box_R(U)$ and $R(y) \subseteq U$. Hence $z \in U$, a contradiction. This proves that ${\leq_1} \circ R \subseteq R$.  A similar reasoning gives that $R \circ {\leq_2}  \subseteq R$.
\end{proof}

\subsection{The functors}\label{sec:Cat}

We conclude the presentation of the duality for the algebraic counterpart $\Alg\Sm$ of a fitler-distributive, finitary,  and congruential logic $\Sm$ with theorems  by showing  the functors and the natural transformations involved in it. 
The \mbox{$\Sm$-algebras} together with homomorphisms between them form a category, that we denote by $\cAlg\Sm$. Before proving the categorical duality for $\cAlg\Sm$, 
we need to   show that \Sp spaces and \Sp morphisms form a category as well.

Similarly to the case of distributive meet-semilattices,  the set-theoretic relational composition of two composable $\Sm$-Priestley morphisms may not be an $\Sm$-Priestley morphisms. Hence we can not use this operation  to obtain a category. 
The operation that  works is, as for distributive meet-semilattices, the following one.  
If   $\spa{X}_{1},\spa{X}_{2}$ and $\spa{X}_{3}$ are \Sp spaces and  $R\subseteq X_{1}\times X_{2}$ and $S\subseteq X_{2}\times X_{3}$ are \Sp morphisms, the composition $(S\star R)\subseteq X_{1}\times X_{3}$ is the relation  defined by: 
\begin{align*}
		(x,z)\in (S\star R) &\IFF \forall U\in B_{3}\big(x\in \Box_{R}\circ\Box_{S} (U) \Rightarrow  z\in U   \big)\\
		&\IFF \forall U\in B_{3}\big( (S\circ R)(x)\subseteq U \Rightarrow  z\in U   \big).
\end{align*}

\begin{thm}\label{thm:cPr}  Let $\langle X_{1},\tau_{1},\BB_{1}\rangle$, $\langle X_{2},\tau_{2},\BB_{2}\rangle$ and $\langle X_{3},\tau_{3},\BB_{3}\rangle$ be  \Sp spaces and let $R\subseteq X_{1}\times X_{2}$ and $S\subseteq X_{2}\times X_{3}$ be \Sp morphisms. Then:
\begin{enumerate}
	\item The \Sp morphism ${\leq_{2}}\subseteq X_{2}\times X_{2}$ satisfies:
	\begin{enumerate}
	\item ${\leq_{2}}\circ R=R\text{ and }S\circ{ \leq_{2}}=S,$
	\item ${\leq_{2}}\circ R = {\leq_{2}}\star R$ and $S\circ{ \leq_{2}} = S\star{ \leq_{2}}$,
	\end{enumerate}
	\item $(S\star R)\subseteq X_{1}\times X_{3}$ is an \Sp morphism.
\end{enumerate}
\end{thm}

\begin{proof} To prove (2) note that  Conditions \ref{itm:PrR1} and \ref{itm:PrR2} follow easily from the definition of $\star$. 
We proceed to prove (1.a).  First we show that ${\leq_{2}}\circ R=R$. 
Let $y\in R(x)$ and $y\leq_{2} z$, and suppose, towards a contradiction, that $z\notin R(x)$. 
By \ref{itm:PrR2} there is $U\in B_{2}$ such that $R(x)\subseteq U$ and $z\notin U$. 
Then by assumption $y\in U$, and since $U$ is an up-set, we get $z\in U$, a contradiction. 
Hence we have ${\leq_{2}}\circ R\subseteq R$. 
The other inclusion is immediate. 
Now we show that $S\circ{ \leq_{2}}=S$. 
Let $x\leq_{2} y$ and $z\in S(y)$, and suppose, towards a contradiction, that $z\notin S(x)$. 
By \ref{itm:PrR2} again, there is $U\in B_{3}$ such that $S(x)\subseteq U$ and $z\notin U$. 
Then we have $x\in\Box_{S} (U)$ and by \ref{itm:PrR1} we get $\Box_{S} (U)\in B_{2}$. 
In particular $\Box_{S}(U)$ is an up-set, thus $y\in \Box_{S}(U)$. 
Then $S(y)\subseteq U$, and therefore $z\in U$, a contradiction. 
Hence we have $S\circ{ \leq_{2}}=S$. 
The other inclusion is immediate.  Finally we prove (1.b). 
The inclusion from  left to right follows by definition. 
For the other inclusion, let ${(x,z)}\in {({\leq_{2}}\star {R})}$ and suppose, towards a contradiction, that $(x,z)\notin {\leq_{2}}\circ {R}$. 
By item (1) we know that ${\leq_{2}}\circ {R}=R$, and then from the hypothesis and \ref{itm:Pr2}, there is $U\in B_{2}$ such that $R(x)\subseteq U$ and $z\notin U$. But since $({\leq_{2}}\circ {R})(x)=R(x)$, we conclude $(x,z)\notin ({\leq_{2}}\star {R})$, a contradiction.  A similar proof shows that $S\circ{ \leq_{2}} = S\star{ \leq_{2}}$.
\end{proof}

Proposition  \ref{prop:porder} and Theorem \ref{thm:cPr} imply the next corollary.

\begin{cor} \label{cor:cPr} The \Sp spaces and the  \Sp morphisms with composition  $\star$ form a category.
\end{cor}

Let us denote by $\cPr\Sm$ the category of \Sp spaces and \Sp morphisms. 
On the one hand, we consider the functor $\sOptS:\cAlg\Sm\longrightarrow \cPr\Sm$ such that for  any $\Sm$-algebras $\alg{A}$, $\alg{A}_{1}$ and $\alg{A}_{2}$ and  any homomorphism $h\in\Hom(\alg{A}_{1},\alg{A}_{2})$:
\begin{align*}
\sOptS(\Al)&\coloneqq \langle \OptS(\Al),\tau_{{\Al}},\varphi[\Al]\rangle,\\
\sOptS(h)&\coloneqq R_{h}\subseteq \OptS(\Al_{2})\times \OptS(\Al_{1}).
\end{align*}
Clearly, for  the identity morphism $\funid_{\Al}$ for $\Al \in \cAlg\Sm$, 
we get $R_{\funid_{\Al}}={\subseteq}$,
and  this is  the identity morphism  for $\sOptS(\Al)$ in $\cPr\Sm$. 
Moreover, it easily follows  from the definition of the dual relation of a homomorphism between $\Sm$-algebras  that for $\Sm$-algebras $\Al_{1},\Al_{2}$ and $\Al_{3}$ and  homomorphisms $f\in\Hom(\Al_{1},\Al_{2})$ and $g\in\Hom(\Al_{2},\Al_{3})$, 
$R_{g\circ f}=R_{f}\star R_{g}$. 
Therefore, using propositions \ref{cor:Salg-Prspa} and \ref{cor:Shom-Prmor}, we conclude that the functor $\sOptS$ is well defined. 
%

\label{defn:pfunctor-sa} 
On the other hand, we consider the functor $(\phantom{.})^{\bullet}:\cPr\Sm\longrightarrow \cAlg\Sm$ 
such that for  any \Sp spaces $\spa{X}$, $\spa{X}_{1}$, $\spa{X}_{2}$ and  any \Sp morphism  $R\subseteq  X_{1}\times  X_{2}$:
\begin{align*}
 \spa{X}^{\bullet}&\coloneqq \BB,\\
R^{\bullet}&\coloneqq \Box_{R}:B_{2}\longrightarrow B_{1}.
\end{align*}
%

In order to complete the duality, we need to define two natural isomorphisms, 
one between the identity functor on $\cAlg\Sm$ and $(\sOptS(\phantom{.}))^{\bullet}$, 
and the other between the identity functor on $\cPr\Sm$ and $\sOptS((\phantom{.})^{\bullet})$. 
Consider first the family of morphisms in $\cAlg\Sm$:
$$\nvarphi_{\Sm}:(\varphi_{\Al}:A\longrightarrow \varphi_{\Al}[A])_{\Al\in\cAlg\Sm}$$

\begin{thm} \label{thm:pnat-a} 
$\nvarphi_{\Sm}$ is a natural isomorphism between the identity functor on $\cAlg\Sm$ and $(\sOptS(\phantom{.}))^{\bullet}$.
\end{thm}
\begin{proof}
Let $\Al_{1},\Al_{2}\in\Alg\Sm$ and let $h\in\Hom(\Al_{1},\Al_{2})$. 
It is enough to show  that $\Box_{R_{h}}\circ \varphi_{1}=\varphi_{2}\circ h$.
For $a\in A_{1}$ and $P\in \Box_{R_{h}}(\varphi_{1}(a))$, we have $R_{h}(P)\subseteq \varphi_{1}(a)$. 
It follows that $h(a)\in P$, so $P\in\varphi_{2}(h(a))$. 
For $P'\in \varphi_{2}(h(a))$, we have $h(a)\in P'$. 
It follows that $R_{h}(P')\subseteq \varphi_{1}(a)$, so $P'\in \Box_{R_{h}}(\varphi_{1}(a))$.

From this we have get that $\nvarphi_{\Sm}$ is a natural transformation, 
and since by Theorem \ref{thm:srepresentation} we know that $\varphi_{1}$ is an isomorphism from $\Al_{1}$ to $\varphi_{1}[\Al_{1}]$, 
we conclude that $\nvarphi_{\Sm}$ is a natural isomorphism.
\end{proof}


Before defining the other natural isomorphism, we need to do some work. 
Recall that for any \Sp space $\spa{X}=\langle X,\tau,\BB\rangle$, 
 the map $\xis_{\spa{X}}:X\longrightarrow \OptS(\BB)$ defined on Section \ref{sec:Obj}   is a homeomorphism between the topological spaces $\langle X,\tau\rangle$ and $\langle \OptS(\BB),\tau_{{\BB}}\rangle$. 
This map encodes the natural isomorphism we are looking for, but since morphisms in $\cPr\Sm$ are relations, we need to give a relation associated  with this map.
We define  the relation $\xisa_{\spa{X}}\subseteq X\times \OptS(\BB)$  by:
\begin{equation*}
(x,P)\in\xisa_{\spa{X}} \IFF \xis_{\spa{X}}(x)\subseteq P.
\end{equation*}

\begin{prop} \label{prop:pepsaok}$\xisa_{\spa{X}}$ is an \Sp morphism.
\end{prop}
\begin{proof} 
We have to show that $\Box_{\xisa_{\spa{X}}}\in\Hom(\varphi_{\BB}[\BB],\BB)$. 
Notice that for all $\varphi_{\BB}(b)\in\varphi_{\BB}[B]$, we have:
\begin{align*}\Box_{\xisa_{\spa{X}}}(\varphi_{\BB}(b))&=\{x\in X:\forall y\in X\big( (x,\xis_{\spa{X}}(y))\in\xisa_{\spa{X}}\Rightarrow \xis_{\spa{X}}(y)\in \varphi_{\BB}(b)  \big)\}\\
&=\{x\in X:\forall y\in X\big( \xis_{\spa{X}}(x)\subseteq \xis_{\spa{X}}(y) \Rightarrow b\in\xis_{\spa{X}}(y) \big)\}\\
&=\{x\in X:b\in \xis_{\spa{X}}(x)\}=b.
\end{align*}
Therefore $\Box_{\xisa_{\spa{X}}}=\varphi_{\BB}^{-1}$. 
And since  $\BB$ and $\varphi_{\BB}[\BB]$ are isomorphic $\Sm$-algebras by means of the map $\varphi_{\BB}$, it follows that $\Box_{\xisa_{\spa{X}}}\in\Hom(\varphi_{\BB}[\BB],\BB)$. 
This proves that  condition \ref{itm:PrR1} is satisfied by $\xisa_{\spa{X}}$. 

We prove now that condition \ref{itm:PrR2} is also satisfied by $\xisa_{\spa{X}}$. 
Notice that for each $x\in X$, we have $\xisa_{\spa{X}}(x)={\uparrow}\xis_{\spa{X}}(x)$. 
Let $x\in X$ and $P\in\OptS(\BB)$ be such that $(x,P)\notin \xisa_{\spa{X}}$. 
We have to show that there is $U\in B$ such that $P\notin \varphi_{\BB}(U)$ and $\xisa_{\spa{X}}(x)\subseteq \varphi_{\BB}(U)$. 
By definition of $\xisa_{\spa{X}}$, we have that $\xis_{\spa{X}}(x)\nsubseteq P$, so there is $U\in B$ such that $U\in\xis_{\spa{X}}(x)$ and $U\notin P$. 
Hence $P\notin \varphi_{\BB}(U)$ and $\xis_{\spa{X}}(x)\in\varphi_{\BB}(U)$. 
Now since $\xisa_{\spa{X}}(x)={\uparrow}\xis_{\spa{X}}(x)$, we obtain that $\xisa_{\spa{X}}(x)\subseteq \varphi_{\BB}(U)$, as required.
Finally, by previous argument we conclude that $\xisa_{\spa{X}}$ is an isomorphism in $\cPr\Sm$. 
\end{proof}

Consider now the family of morphisms in $\cPr\Sm$:
$$\nxisa_{\Sm}=\big( \xisa_{\spa{X}}\subseteq X\times \OptS(\BB)\big)_{\spa{X} \in \cPr\Sm}$$

\begin{thm} \label{thm:pnat-s} 
$\nxisa_{\Sm}$ is a natural isomorphism between the identity functor on $\cPr\Sm$ and $\sOptS((\phantom{.})^{\bullet})$.
\end{thm}

\begin{proof}  
Let $\spa{X}_{1}=\langle  X_{1},\tau_{1},\BB_{1}\rangle $ and $\spa{X}_{2}=\langle  X_{2},\tau_{2},\BB_{2}\rangle $ be two \Sp spaces and let $R\subseteq  X_{1}\times  X_{2}$ be an \Sp morphism. 
First we show that:
\begin{equation*}
(x,y)\in R \IFF (\xis_{1}(x),  \xis_{2}(y))\in R_{\Box_{R}}.
\end{equation*}

Let $x\in  X_{1}$ and $y\in  X_{2}$ be such that $(x,y)\in R$, and let $U\in B_{2}$.  
Notice that we have:
\begin{align*} U\in \Box^{-1}_{R}[\xis_{1}(x)] 
&\IFF \Box_{R}(U)\in\xis_{1}(x) 
\IFF x\in \Box_{R}(U)\IFF R(x)\subseteq U.
\end{align*} 
Thus if $U\in \Box^{-1}_{R}[\xis_{1}(x)]$, then $R(x)\subseteq U$, and since $(x,y)\in R$, we obtain $y\in U$, \ie $U\in\xis_{2}(y)$, and therefore $(\xis_{1}(x),  \xis_{2}(y))\in R_{\Box_{R}}$. 
For the converse,  let $x\in  X_{1},y\in  X_{2}$ be such $(\xis_{1}(x) , \xis_{2}(y))\in R_{\Box_{R}}$ and suppose, towards a contradiction, that $y\notin R(x)$. 
Since $R$ is an \Sp morphism, by \ref{itm:PrR1}, there is $U\in B_{2}$ such that $y\notin U$ and $R(x)\subseteq U$. 
From previous equivalences we obtain  $U\in \Box^{-1}[\xis_{1}(x)]$. 
But then from the  hypothesis $U\in\xis_{2}(y)$, so $y\in U$, a contradiction. 

The equivalence that we just proved implies that 
$R_{\Box_{R}}\star \xisa_{\spa{X}_{1}}= \xisa_{\spa{X}_{2}}\star R$. 
Thus $\nxisa_{\Sm}$ is a natural equivalence. 
Moreover, as $\xisa_{\spa{X}}$ is an isomorphism for each \Sp space $\spa{X}$, then $\nxisa_{\Sm}$ is a natural isomorphism.  
\end{proof}

\begin{thm} \label{thm:PrSduality} 
The categories $\cAlg\Sm$ and $\cPr\Sm$ are dually  equivalent by means of the contravariant functors $\sOptS$ and $(\phantom{.})^{\bullet}$ and the natural equivalences $\nvarphi_{\Sm}$ and $\nxisa_{\Sm}$.
\end{thm}

\section{Dual correspondence of some logical properties}

In this  final section we examine how the correspondences between objects of the categories we are considering can be refined depending on the properties of the logic under consideration. 
For information on  the abstract properties of logics we consider in the sequel  we refer the reader to \cite{FJa09}.  
Given the abstract character of our general approach, 
we  carry out this study modularly, 
treating each property independently, obtaining in this way results that can be combined afterwards.

The next proposition contains a fact on congruential logics that will be used in the proofs of the correspondence theorems we present between properties of a  filter-distributive and  finitary congruential logic $\Sm$ with theorems  and properties of its \Sp spaces.

\begin{prop}\label{prop:bilogic-quotient}
Let $\Sm$ be a congruential logic. 
For every algebra $\Al$, the quotient homomorphism  $\pi_{\Al}:\Al\longrightarrow \Al/{\equivsa}$ induces an isomorphism between the lattices $\FiSA$ and $\Fi_{\Sm}(\Al/{\equivsa})$ given by $F \mapsto \pi_{\Al}[F]$ and whose  inverse is given by $G \mapsto \pi_{\Al}^{-1}[G]$. 
\end{prop}

Given a logic $\Sm$ and a set of $\mathcal{L}_{\Sm}$-formulas $\Gamma$, we abbreviate $\Cons^{\Fm}(\Gamma)$ by $\Cons(\Gamma)$.  
It holds that for every formula $\varphi$, $\varphi \in \Cons(\Gamma)$ if and only if $\Gamma \vdash_{\Sm} \varphi$.

\subsection{The Property of Conjunction.}\label{subsec:T-DCL-Corr-PC}   

A logic   $\Sm$ has the \emph{property of conjunction}   (PC) for a formula $x\wedge y$ in two variables if :
$$x\wedge y \vdash_{\Sm}x,\AAAND 
x\wedge y\vdashs y,\AAAND
x,y\vdashs x\wedge y.$$
In other words, if $\Cons(x, y) = \Cons(x \wedge y)$.
By the substitution-invariance of $\vdashs$ we can replace $x$ and $y$ by arbitrary formulas. 
The property of conjunction transfers to every algebra, see  \cite[p.\ 50]{FJa09}, in the sense that if $\Sm$ has (PC) for $x \wedge y$, then for every algebra $\Al$ and every $a,b\in A$ 
 $$\Consa(a\wedge^{\Al}b)=\Consa(a,b),$$
where  $a\wedge^{\Al}b$ is the interpretation of the term $x \wedge y$ in $\Al$ when $a$ is assigned to $x$ and $b$ to $y$.

For the remaining part of the subsection,  let \emph{$\Sm$ be a filter-distributive and  finitary congruential logic with theorems}.
From the transfer of (PC)  to every algebra it easily follows that: 

\begin{lemma}\label{lemma:PCa}
If $\Sm$ satisfies  (PC), then for every $\Sm$-algebra $\Al$ and  all $a,b\in A$, $\varphi(a)\cap \varphi(b)=\varphi(a\wedge^{\Al} b)$.
\end{lemma}

Notice that by the associativity of intersection, 
the previous lemma implies that for any non-empty and finite $B\subseteq A$, $\bigcap \{\varphi(b):b\in B\}=\varphi(\bigwedge^{\Al}B)$. 
Recall that we defined the $\Sm$-semilattice of $\Al$ as the closure of $\varphi[A]$ under  finite intersections. 
Therefore, if $\Sm$ satisfies (PC), then $\langle  A, \wedge^{\Al}, 1^{\Al} \rangle$ and $\Mee(\Al)$ are isomorphic.

\begin{prop}\label{prop:PCb1}
If $\Sm$ satisfies  (PC), then for every $\Sm$-algebra $\Al$, then
\begin{enumerate}
\item  for all $U\subseteq \OptSA$, if $U$ is an $\IrrSA$-admissible clopen up-set of the space $\langle \OptSA,\tau_{\Al},\subseteq\rangle$, 
then there is $a\in A$ such that $U=\varphi(a)$,
\item $\varphi[A]$ is 
 the collection of $\IrrSA$-admissible clopen up-sets of $\langle \OptSA,\tau_{\Al},$ $ \subseteq\rangle$. 
\end{enumerate}
\end{prop}
\begin{proof}
(1) Let $U\subseteq \OptSA$ be a clopen up-set of $\langle \OptSA,\tau_{\Al},\subseteq\rangle$ such that $\max(U^{c})\subseteq \IrrSA$. 
Then by Proposition \ref{prop:isom2}, there is a non-empty finite $B\subseteq A$ such that $U=\varphim(B)=\bigcap\{\varphi(b):b\in B\}$.  
 Lemma \ref{lemma:PCa} implies that   $U=\varphi(\bigwedge^{\Al} B)$, as required. 
\end{proof}


We can conjecture that 
 the property of \Sp spaces that  corresponds  to (PC) is  ``on every \Sp space $\langle X,\tau,\BB\rangle$  the set $B$ is the collection of $\Xb $-admissible clopen up-sets''. This is indeed the case. It can be proved using 
 Proposition \ref{prop:bilogic-quotient}   and the next proposition.
%


\begin{prop}\label{prop:PCb2}
Let $\langle X,\tau,\BB\rangle$ be an \Sp space such that $B$ is the collection of the $\Xb$-admissible clopen up-sets of $X$. 
Then for all  $U,V\in B$, $\Consb(U,V)=\Consb(U\cap V)$. 
\end{prop}
\begin{proof}
First notice that the hypothesis implies that $B$ is closed under finite intersections. 
Now let $U,V\in B$. 
By \ref{itm:Pr2} we get $\bigcap \{U,V\}\subseteq U\cap V$ if and only if $U\cap V\in\Consb(U,V)$. Therefore,  $U\cap V\in\Consb(U,V)$. Then, considering Remark \ref{rem:Pr2-2}, it is easy to see that $\Consb(U,V)=\Consb(U\cap V)$. 
\end{proof}

\begin{prop}\label{prop:PCb3}
If $\Sm$ is  such that for every \Sp space $\langle X,\tau,\BB\rangle$ the set  $B$ is the collection of the $\Xb$-admissible clopen up-sets of $X$,  
then $\Sm$ satisfies (PC). 
\end{prop}
\begin{proof}
Recall that the Lindenbaum-Tarski algebra $\Fm/{\equivsfm}$ is an \mbox{$\Sm$-algebra}. We abbreviate $\equivsfm$ by $\equiv$. 
Let $x,y\in Var$. 
The assumption implies that there is a formula  $\rho \in Fm_{\mathscr{L}}$ such that $\varphi(x/{\equiv})\cap\varphi(y/{\equiv})=\varphi(\rho/{\equiv})$. 
Moreover, by Proposition \ref{prop:PCb2} we have 
\begin{equation*}
\mathrm{Fg}_{\Sm}^{\varphi[\Fm/{\equiv}]}\left(\varphi(x/{\equiv}),\varphi(y/{\equiv})\right)=\mathrm{Fg}_{\Sm}^{\varphi[\Fm/{\equiv}]}\left(\varphi(x/{\equiv})\cap\varphi(y/{\equiv})\right). 
\end{equation*}
Then by Theorem \ref{thm:repr} we obtain $\mathrm{Fg}_{\Sm}^{\Fm/{\equiv}}(x/{\equiv},y/{\equiv})=\mathrm{Fg}_{\Sm}^{\Fm/{\equiv}}(\rho/{\equiv})$. 
And using Proposition \ref{prop:bilogic-quotient} we get
$\Cons(x,y)=\Cons(\rho )$. 
By the substitution-invariance of $\vdash_\Sm$, it follows that
there is a formula $\rho'(x,y)$ in at most the variables $x$ and $y$ such that $\Cons(x,y)=\Cons(\rho'(x,y))$. 
Hence $\Sm$ satisfies (PC)  for  $\rho'(x, y) $. 
\end{proof}

Propositions \ref{prop:PCb1} and  \ref{prop:PCb3} imply  the desired theorem. 

\begin{thm}\label{cor:PC3b}
Let $\Sm$ be a filter-distributive and finitary congruential logic with theorems. Then $\Sm$ satisfies (PC) if and only if for every \Sp space $\langle X,\tau,\BB\rangle$, $B$ is the collection of the $\Xb$-admissible clopen up-sets of $X$. 
\end{thm}

\subsection{The Property of Disjunction.}\label{subsec:T-DCL-Corr-PDI}   

A logic  $\Sm$  satisfies the \textit{property of  disjunction} (PDI) for  a set of formulas in two variables  $\nabla(x, y)$  if for all formulas $\delta,\gamma,\mu$ and set of formulas $\Gamma$ satisfies: 
$$\text{(a) }  \delta\vdash_{\Sm} \nabla(\delta, \gamma), \,\,\,\,\,\,\,\,\,
	\text{(b) }\delta\vdash_{\Sm} \nabla(\gamma,\delta), $$
	$$
	\text{(c) } \text{ if }\Gamma,\delta\vdash_{\Sm}\mu \text{ and } \Gamma,\gamma\vdash_{\Sm}\mu \text{, then } \Gamma, \nabla(\delta, \gamma)\vdash_{\Sm}\mu.$$
	%
%
We say that $\Sm$ satisfies (PDI) if it satisfies (PDI) for some set of formulas. 
%

If $\Sm$ satisfies (PDI) for  $\nabla(x, y)$, then this property transfers to every algebra (\cf Corollary 2.5.4 in \cite{Cz01} or Theorem 2.52 in \cite{FJa09}), that is,   for every algebra $\Al$ and every $\{a,b\}\cup X\subseteq A$: 
$$\Consa(X,\nabla^{\Al}(a, b)) =\Consa(X,a)\cap\Consa(X,b).$$
Moreover, it is  known that if a logic satisfies (PDI) then it is filter-distributive (\cf \cite{Cz84}). 
%
Next lemma (see Proposition 2.5.7 in \cite{Cz01}) will be used in the sequel.

\begin{lemma}\label{lemma:PDa-3} 
If $\Sm$  satisfies (PDI) for a set of formulas $\nabla(x,y)$, then for every   $\Sm$-algebra $\Al$, 
an $\Sm$-filter $F\in \FiSA$ is irreducible if and only if for every $a, b \in A$, 
if $\nabla^{\Al}(a, b) \subseteq F$, then $a \in F$ or $b \in F$. 
\end{lemma}

Another important fact on (PDI)  follows from Theorem 2.5.9 in \cite{Cz01}, as observed in \cite{CiNo13} 
 (taking into account that in \cite{Cz01} the irreducible $\Sm$-filters are called prime):

\begin{thm}
\label{thm:f-dist-PDI-Charact}
A filter-distributive and finitary logic $\Sm$ satisfies (PDI)  if and only if 
for all algebras $\Al, \BB$, 
every homomorphism $h\in\Hom(\Al,\BB)$ 
and every $G\in\IrrSB$, the $\Sm$-filter $h^{-1}[G]$ of $\Al$ is irreducible.
\end{thm}

In our setting of congruential, filter-distributive and finitary logics, 
we can restrict the class of algebras in Theorem \ref{thm:f-dist-PDI-Charact} to $\Alg\Sm$ using   the next proposition.

\begin{prop}
Let $\Sm$ be a congruential logic. 
If for all algebras $\Al, \BB \in \Alg\Sm$, 
every homomorphism $h\in\Hom( \Al, \BB)$ 
and every  $G\in\IrrSB$ the $\Sm$-filter $h^{-1}[G]$ of $\Al$ is irreducible, then the same holds for arbitrary algebras. 
\end{prop}
\begin{proof}
Let $\Al, \BB$ be arbitrary algebras and let $h: \Al \to \BB$ be a homomorphism. 
Consider the quotient algebras $\Al/{\equivsa}$ and $\BB/{\equivsb}$, which belong to $\Alg\Sm$. 
Let  $\pi_\Al: \Al \to  \Al/{\equivsa}$ and $\pi_\BB: \BB \to  \BB/{\equivsb}$ the  quotient   homomorphisms.  
We define the map $h': A/{\equivsa} \to B/{\equivsb}$ by setting for every $a \in A$,
$h'(a/{\equivsa}) = h(a)/{\equivsb}$. 
This map is well defined, as if   $a \equivsa b$ and $F \in \FiSB$ is such that $h(a) \in F$, 
then $a \in h^{-1}[F] \in \FiSA$ and therefore, since $a \equivsa b$, we obtain $b \in h^{-1}[F]$ and so $h(b) \in F$. 
Moreover, $h'$ is a homomorphism and $h' \circ \pi_{\Al} = \pi_{\BB} \circ h$. 
By Proposition \ref{prop:bilogic-quotient}, $\pi_{\Al}[\cdot]$ establishes a lattice isomorphism between $\FiSA$ and $\Fi_{\Sm}(\Al/{\equivsa})$ and $\pi_{\BB}[\cdot]$  a lattice isomorphism between $\FiSB$ and $\Fi_{\Sm}(\BB/{\equivsb})$.
Thus if $G$ is an irreducible $\Sm$-filter of $\BB$, then $\pi_{\BB}[G]$ is an irreducible $\Sm$-filter of $\BB/{\equivsb}$. 
Hence from the assumption follows that $h'^{-1}[\pi_{\BB}[G]]$  is an irreducible $\Sm$-filter of $\Al/{\equivsa}$. 
Therefore, $\pi_{\Al}^{-1}[h'^{-1}[\pi_{\BB}[G]]]$ is an  irreducible $\Sm$-filter of $\Al$. 
But $\pi_{\Al}^{-1}[h'^{-1}[\pi_{\BB}[G]]] = (h' \circ \pi_{\Al})^{-1}[\pi_{\BB}[G]] = (\pi_{\BB} \circ h)^{-1}[\pi_{\BB}[G]] = h^{-1}[\pi_{\BB}^{-1}[\pi_{\BB}[G]]] = h^{-1}[G]$ and therefore $h^{-1}[G]$ is  an irreducible $\Sm$-filter  of $\Al$.
\end{proof}

\begin{cor}
\label{cor:PDI-hom}
A congruential, filter-distributive and finitary logic $\Sm$ satisfies (PDI)  if and only if 
for all algebras $\Al, \BB \in \Alg\Sm$, 
every homomorphism $h\in\Hom(\Al,\BB)$ 
and every $G\in\IrrSB$,
 the $\Sm$-filter $h^{-1}[G]$ of $\Al$ is irreducible.
\end{cor}

We look now for a translation of this property on the inverse images of irreducible $\Sm$-filters under homomorphisms  
into a property of morphisms of the dual \Sp spaces. 
From now on in this subsection \emph{we fix a congruential, filter-distributive and finitary logic $\Sm$}.

\begin{prop}
\label{prop:inverse-irr-irr}
For all algebras $\Al, \BB \in \Alg\Sm$ and every $h\in\Hom(\Al,\BB)$ 
the following conditions are equivalent:
\begin{enumerate} 
\item for every $G \in \IrrSB$, $h^{-1}[G] \in \IrrSA$,
\item for every $G \in \IrrSB$ there exists $F \in \IrrSA$ such that for all $a \in A$, $F \in \varphi(a)$ if and only if $R_{h}(G) \subseteq \varphi(a)$.
\end{enumerate}
\end{prop}
\begin{proof}
Note that condition (2) is equivalent to saying  that $\bigcap R_{h}(G) \in \IrrSA$, for every $G \in \IrrSB$. 
Also note that  from Proposition \ref{prop:ofls} it follows that for any $G\in\FiSB$, $h^{-1}[G]$ equals the intersection of all the irreducible $\Sm$-filters of $\Al$ that include $h^{-1}[G]$. 
Moreover, since every irreducible $\Sm$-filter is optimal, every irreducible $\Sm$-filter that includes $h^{-1}[G]$ belongs to $R_h(G)$. Therefore, since $h^{-1}[G] \subseteq H$ for every $H \in R_{h}(G)$ it follows that $h^{-1}[G]  =  \bigcap R_{h}(G)$. From this last observation follows  that (1) implies (2).
Suppose now (2). Let $G \in \IrrSB$. Then, by (2), $\bigcap R_{h}(G) \in \IrrSA$. Since $h^{-1}[G]  =  \bigcap R_{h}(G)$ we conclude that  $h^{-1}[G] \in \IrrSA$. 
\end{proof}

We say that an \Sp  morphism $R \subseteq X_1 \times X_2$  from an \Sp space $\langle X_1,\tau_1,\BB_1\rangle$ to an \Sp space $\langle X_2,\tau_2,\BB_2\rangle$ is \textit{functional} if 
\begin{equation*}
\tag{E1}\label{eq:PDI-dual-2}
(\forall u \in X_1)(\exists v \in X_2) R(u)= {\uparrow}v.
\end{equation*}
 The name is due to the fact that then from $R$ we can define a function $f_R: X_1 \to X_2$ by setting for every $u \in X_1$ $f_R(u)$ as  $\text{the only $v \in X_2$ such that } R(u) = {\uparrow}v$.

 The condition (\ref{eq:PDI-dual-2}) is easily seen to be equivalent to the following one:
 \begin{equation*}
\tag{E2}\label{eq:PDI-dual-a}
(\forall u \in X_{1})(\exists v \in X_{2})(\forall U \in B_2)  (v \in U \Leftrightarrow R(u) \subseteq U). 
\end{equation*}

We obtain the next characterization of having (PDI) where condition (\ref{eq:PDI-dual-a}) is restricted to the elements of $X_{\BB_1}$ and $X_{\BB_2}$.

\begin{thm}
\label{thm:PDI-dual}
The logic $\Sm$ satisfies (PDI) if and only if  
for  every \Sp  morphism $R \subseteq X_1 \times X_2$  from an \Sp space $\langle X_1,\tau_1,\BB_1\rangle$ to an \Sp space $\langle X_2,\tau_2,\BB_2\rangle$ the following condition holds:
\begin{equation*}
\tag{E3}\label{eq:PDI-dual}
(\forall u \in X_{\BB_1})(\exists v \in X_{\BB_2})(\forall U \in B_2)  (v \in U \Leftrightarrow R(u) \subseteq U). 
\end{equation*}
 \end{thm}
\begin{proof}
If $\Sm$ satisfies (PDI), then using the duality we have developed and Proposition \ref{prop:inverse-irr-irr}  we obtain condition (\ref{eq:PDI-dual}). 
Conversely, if (\ref{eq:PDI-dual}) holds  for  every \Sp morphism $R \subseteq X_1 \times X_2$  from an \Sp space $\langle X_1,\tau_1,\BB_1\rangle$ to an \Sp space $\langle X_2,\tau_2,\BB_2\rangle$, then the duality we have developed, Proposition \ref{prop:inverse-irr-irr} and Corollary \ref{cor:PDI-hom} allow us to conclude that $\Sm$ satisfies (PDI).
\end{proof}

When $\Sm$ satisifes (PDI) for a finite set of formulas we have the next proposition. 
\begin{prop}\label{prop:PDa-4} 
If $\Sm$ satisfies (PDI) for a finite set of formulas $\nabla(x,y)$, then for every $\Sm$-algebra $\Al$, $\OptSA=\IrrSA$.
\end{prop}
\begin{proof}
We just need to show that every optimal $\Sm$-filter of $\Al$ is irreducible,  
so let $P$ be an optimal $\Sm$-filter of $\Al$ and let $a, b \in A$ 
be such that $\nabla^{\Al}(a, b) \subseteq P$. 
Assume, towards a contradiction, that $a, b \not \in P$. 
Since $\Consa(a) \cap \Consa(b) = \Consa(\nabla^{\Al}(a, b))$ and $P^{c}$ is a strong $\Sm$-ideal it follows that $ \Consa(\nabla^{\Al}(a, b)) \cap P^{c} \neq \emptyset$, a contradiction. 
Hence, if $\nabla^{\Al}(a, b) \subseteq P$, then $a \in P$ or $b \in P$. Thus, Lemma \ref{lemma:PDa-3} implies that $P$ is irreducible. 
\end{proof}

Taking into account Proposition \ref{prop:PDa-4} we have the following result:

\begin{prop}
\label{prop:PDI-dual-2}
If  $\Sm$ satisfies (PDI) for some finite  set of formulas $\nabla(x,y)$, then   for every \Sp space $\langle X ,\tau ,\BB\rangle$ we have $X = \Xb$  and  every \Sp  morphism between \Sp spaces is functional. 
 \end{prop}
 \begin{proof}
 If $\Sm$ satisfies (PDI)  for a finite  set of formulas $\nabla(x,y)$, then  the duality we have developed together with Proposition \ref{prop:PDa-4}  show that 
 $X = \Xb$ for every \Sp space $\langle X ,\tau ,\BB\rangle$. Then  condition (\ref {eq:PDI-dual}) in Theorem \ref{thm:PDI-dual} implies condition (\ref{eq:PDI-dual-a}) that is equivalent to the condition that defines being functional. 
 \end{proof}


 The consequent of the statement of Proposition \ref{prop:PDI-dual-2} implies that $\Sm$ has (PDI), because it implies condition  (\ref{eq:PDI-dual}) in Theorem \ref{thm:PDI-dual}.

Now we concentrate in the case where (PDI) holds for a one-element set $\nabla = \{x \vee y\}$, where $x \vee y$ is a formula in two variables.  
In this case we say $\Sm$ satisfies (PDI) for  $x \vee y$.

\begin{lemma}\label{lemma:PDa} 
If $\Sm$  satisfies (PDI) for $x \vee y$, then for every   $\Sm$-algebra $\Al$ and all $a,b\in A$, $\varphi(a)\cup \varphi(b)= \varphi(a \vee^{\Al} b)$.
\end{lemma}
\begin{proof}
Notice that since (PDI) transfers to every algebra, 
for all $a,b\in A$ we have $\Consa(a \vee^{\Al} b)=\Consa(a)\cap\Consa(b)$. 
This implies that $a,b\leqsa a \vee^{\Al} b$ and therefore 
for any $P\in\OptSA$, we have that if $a\in P$ or $b\in P$, then $a \vee^{\Al} b\in  P$, because $P$ is an up-set. This proves that $\varphi(a)\cup \varphi(b) \subseteq \varphi(a \vee^{\Al} b)$. To prove the other inclusion assume that $P \in \varphi(a \vee^{\Al} b)$.  Since, by Proposition \ref{prop:PDa-4}, $P$ is irreducible, by Lemma \ref{lemma:PDa-3} follows that $a \in P$ or $b \in P$; hence $P \in \varphi(a)\cup \varphi(b)$.
\end{proof}

\begin{cor}\label{cor:PDb1}
If $\Sm$  satisfies (PDI) for a single formula $x \vee y$, then for every $\Sm$-algebra $\Al$, 
 $\varphi[A]$ is closed under the binary operation of union.  
\end{cor}

We recall that for every logic $\Sm$ and every algebra $\Al$, every $\Sm$-filter $F$ of $\Al$ is equal to the intersection of all the irreducible $\Sm$-filters of $\Al$ that include $F$. In particular this holds for the theories of $\Sm$.

\begin{prop}
\label{prop:PDi-single}
If $\Sm$  satisfies (PDI) and for every $\Sm$-algebra $\Al$ the set   
 $\varphi[A]$ is closed under the binary operation of union, then $\Sm$ satisfies (PDI) for a single formula. 
\end{prop}
\begin{proof}
Let us consider the relation $\equiv_{\Sm}^{\Fm}$ that we abbreviate all along the proof by $\equiv$. Then the quotient algebra $\Fm/{\equiv} \in \Alg\Sm$. Let $p, q$ be variables
and consider the equivalence classes $x/{\equiv}, y/{\equiv}$. By assumption  $\varphi(x/{\equiv}) \cup \varphi(y/{\equiv}) = \varphi(\delta/{\equiv})$ for some formula $\delta$. 

We first prove that $\Cons(x) \cap \Cons(y) = \Cons(\delta)$.  Let $T$ be an irreducible $\Sm$-theory such that $\Cons(x) \cap \Cons(y) \subseteq T$. Then $T = (\Cons(x) \cap \Cons(y)) \sqcup T$. Since $\Sm$ is filter-distributive, $T = (\Cons(x) \sqcup T) \cap (\Cons(y)\sqcup T)$. Now being $T$ irreducible, it follows that $T = \Cons(x) \sqcup T$ or $T = \Cons(y) \sqcup T$. Hence $x \in T$ or $y \in T$. Therefore, $x/{\equiv} \in \pi[T]$ or $y/{\equiv} \in \pi[T]$, where $\pi$ is the quotient homomorphism from $\Fm$ onto $\Fm/{\equiv}$. By Proposition \ref{prop:bilogic-quotient},  the fact that $T$ is irreducible implies that $\pi[T]$ is an irreducible $\Sm$-filter of $\Fm/{\equiv}$; therefore $\pi[T] \in  \varphi(x/{\equiv}) \cup \varphi(y/{\equiv})$. Hence, $\delta/{\equiv} \in \pi[T]$. 
This implies, since by the definition of $\equiv$, $\delta \equiv \delta'$ if and only if $\Cons(\delta) = \Cons(\delta')$, that $\delta \in T$. We conclude that  $\Cons(\delta) \subseteq \Cons(x) \cap \Cons(y)$. To prove the other inclusion, let $T$ be an irreducible $\Sm$ theory such that $\delta \in T$. Then $\delta/{\equiv} \in \pi[T]$ and $\pi[T]$ is an irreducible $\Sm$-filter of $\Fm/{\equiv}$. Therefore, $\pi[T] \in \varphi(\delta/{\equiv})$. Hence $x/{\equiv} \in \pi[T]$ or $y/{\equiv} \in \pi[T]$.
This, by a similar reasoning as before, implies that $x \in T$ or $y \in T$. In both cases $\Cons(x) \cap \Cons(y) \subseteq T$. We conclude that $\Cons(x) \cap \Cons(y) \subseteq \Cons(\delta)$.

Now, since $\Sm$ has (PDI), let us assume that $\Sm$ has (PDI) for  $\nabla(x, y)$. Then $\Cons(x) \cap \Cons(y) = \Cons(\nabla(x,y))$.  Therefore, $\Cons(\nabla(x,y)) = \Cons(\delta)$. Let $\sigma$ be the substitution that maps $p$ to $p$ and all the remanning variables to $q$. Then $\sigma[\nabla(x, y)] = \nabla(x,y)$ and therefore $\Cons(\nabla(x,y)) = \Cons(\sigma(\delta))$ using  invariance under substitutions.
It follows that $\Cons(x) \cap \Cons(y) = \Cons(\sigma(\delta))$ and the variables in $\sigma(\delta)$ are at most $x, y$. Thus $\Sm$ has (WPDI) for $\sigma(\delta)$ and being filter-distributive it also has (PDI) for $\sigma(\delta)$. 
\end{proof}

Combining Proposition \ref{prop:PDI-dual-2} and Proposition \ref{prop:PDi-single} we 
easily obtain that the property of \Sp spaces that corresponds to having  (PDI) for a single formula is:  $X=\Xb$, $B$  is closed under the binary operation of union,  and the \Sp morphisms are functional. We state it in a theorem.

\begin{thm}
The logic $\Sm$  satisfies (PDI) for a single formula if and only if for every \Sp space $\langle X ,\tau ,\BB\rangle$ we have that $X=\Xb$, $B$  is closed under the binary operation of union,  and moreover the \Sp morphisms between \Sp spaces are functional.
\end{thm}

We find another equivalent condition,  this time on morphisms.  
To this end  we consider the next proposition.

  \begin{prop}\label{prop:PDI-duals-hom}
If $\Sm$  satisfies (PDI) for a single formula $x \vee y$, then for every $\Al, \BB \in  \Alg\Sm$ and every $h \in \Hom(\Al, \BB)$, the relation $R_h \subseteq  \OptSB \times  \OptSA$
satisfies that for every $a, b \in A$ and every $P \in \OptSB$, if  $R_h(P) \subseteq \varphi^{\Al}(a) \cup \varphi^{\Al}(b)$, then $R_h(P) \subseteq \varphi^{\Al}(a)$ or $R_h(P) \subseteq \varphi^{\Al}(b)$.
\end{prop}
\begin{proof}
Assume that $R_h(P) \subseteq \varphi^{\Al}(a) \cup \varphi^{\Al}(b)$. Then  $R_h(P) \subseteq \varphi^{\Al}(a \vee^{\Al} b)$. Suppose that $R_h(P) \not \subseteq \varphi^{\BB}(a)$ and  $R_h(P) \not \subseteq \varphi^{\Al}(b)$. Let $Q, Q' \in R_h(P)$ such that $a \not \in Q$ and $b \not \in Q'$. Then $h^{-1}[P] \subseteq Q$ and $h^{-1}[P] \subseteq Q'$. Hence,  $h(a), h(b) \not \in P$ and  therefore, $h(a \vee^{\Al} b) \not \in P$.  
\end{proof}
  
 We  need the following lemma.
 
\begin{lemma}\label{cor:PDb2}
Let $\langle X,\tau,\BB\rangle$ be an \Sp space such that $B$ is closed under the binary operation of union. 
Then for all $\{U,V\}\cup\mathcal{W}\subseteq B$, 
$$\Consb(\mathcal{W},U)\cap \Consb(\mathcal{W},V)=\Consb(\mathcal{W},U\cup V).$$
\end{lemma}
\begin{proof}
We first prove that for all $U,V\in B$, $\Consb(U)\cap \Consb(V)=\Consb(U\cup V)$.  By hypothesis we have for any $W\in B$ that $W\in\Consb(U)\cap\Consb(V)$ if and only if $U\subseteq W$ and $V\subseteq W$ if and only if $W\in\Consb(U\cup V)$, as required. Now using the filter distributivity of the logic we have:
\begin{align*}\Consb(\mathcal{W},U)\cap \Consb(\mathcal{W},V)  &=   \big(\Consb(\mathcal{W})\sqcup\Consb(U)\big)\cap
 \big(\Consb(\mathcal{W})\sqcup\Consb(V)\big)   \\
& =   \Consb(\mathcal{W})\sqcup\big(\Consb(U)\cap\Consb(V)\big)\\ & =  \Consb(\mathcal{W})\sqcup\Consb(U\cup V)  \\
& =   \Consb(\mathcal{W},U\cup V).
\end{align*}
\end{proof}

Lemma \ref{cor:PDb2} implies the next characterization of the irreducible $\Sm$-filters. 

\begin{lemma}
\label{lem:irr-SPri-funions}
Let $\langle X,\tau,\BB\rangle$ be an \Sp space such that $B$ is closed under the binary operation of union and let $F$ be an $\Sm$-filter of $\BB$. 
Then $F$ is irreducible if and only if 
for all $U, V \in B$, if $U \cup V \in F$, then $U \in F$ or $V \in F$.
\end{lemma}
\begin{proof}
Suppose that $F$ is an irreducible $\Sm$-filter of $\BB$ and $U \cup V \in F$ with  $U, V \in B$. 
By Lemma \ref{cor:PDb2}, $\Consb(F,U)\cap \Consb(F,V)=\Consb(F,U\cup V)$. 
Since $U \cup V \in F$, it follows that  $\Consb(F,U)\cap \Consb(F,V)= F$. 
This implies, being $F$ irreducible,  that $\Consb(F,U)= F$ or $\Consb(F,V) = F$; 
hence $U \in F$ or $V \in F$.
Conversely, suppose that for all $U, V \in B$, if $U \cup V \in F$, then $U \in F$ or $V \in F$. 
Assume that $F = H \cap G$ with $F, G \in \FiSA$ and that $F \neq H$ and $F \neq G$. 
Let then $U \in H \setminus F$ and $V \in G \setminus F$. 
Since $B$ is closed under finite unions $U \cup V \in B$. 
Hence, since $H$ and $G$ are upsets, $U \cup V \in H \cap G = F$. 
From the hypothesis follows that $U \in F$ or $V \in F$,  a contradiction. 
\end{proof}

\begin{lemma}
 Let $R \subseteq X_1 \times X_2$ be an \Sp morphism  from an \Sp space $\langle X_1,\tau_1,\BB_1\rangle$ to an \Sp space $\langle X_2,\tau_2,\BB_2\rangle$ such that $\BB_1$ and $\BB_2$ are closed under the binary operation of union and moreover for every $U, V \in B_2$ and every $u \in X_{\BB_1}$, if $R(u) \subseteq U \cup V$, then $R(u) \subseteq U$ or $R(u) \subseteq V$.   Then condition    (\ref{eq:PDI-dual}) holds. 
\end{lemma}
\begin{proof}
Assume the antecedent. Let $u \in X_{\BB_1}$. Then $\varepsilon(u)$ is an irreducible $\Sm$-filter of $\BB_1$. Consider the homomorphism $\Box_{R}: \BB_2 \to \BB_1$. We prove that 
$\Box_{R}^{-1}[\varepsilon(u)]$ is an irreducible $\Sm$-filter of $\BB_2$. Let $U, V \in B_2$ be such that $U \cup V \in \Box_{R}^{-1}[\varepsilon(u)]$. Then $\Box_{R}(U \cup V) \in \varepsilon(u)$. Therefore, $u \in \Box_{R}(U \cup V)$. Thus $R(u) \subseteq U \cup V$. Therefore, by the assumption, $R(u) \subseteq U$ or  $R(u) \subseteq V$. Hence $u \in \Box_{R}(U)$ or $u \in \Box_R(V)$. This implies that $\Box_{R}(U) \in \varepsilon(u)$ or $\Box_{R}(V) \in \varepsilon(u)$; hence $U \in \Box_{R}^{-1}[\varepsilon(u)]$ or $V \in \Box_{R}^{-1}[\varepsilon(u)]$. By Lemma \ref{lem:irr-SPri-funions}
 we obtain that $\Box_{R}^{-1}[\varepsilon(u)]$ an irreducible $\Sm$-filter of $\BB_2$. Then there is $v \in X_{\BB_2}$ such that $\varepsilon(v) = \Box_{R}^{-1}[\varepsilon(u)]$. Therefore, for every $U \in B_2$, $v \in U$ if and only if $u \in \Box_R(U)$ if and only if $R(u) \subseteq U$.
\end{proof}

The next theorem follows using the duality we have developed, Proposition \ref{prop:PDI-duals-hom} and the last lemma.

\begin{thm}
The logic $\Sm$ has (PDI) for a single formula  if and only if for every \Sp space $\langle X,\tau,\BB\rangle$ the set $B$ is closed under the binary operation of union and  for  every \Sp  morphism $R \subseteq X_1 \times X_2$  from an \Sp space $\langle X_1,\tau_1,\BB_1\rangle$ to an \Sp space $\langle X_2,\tau_2,\BB_2\rangle$ it holds that for every $u \in X_{\BB_1}$ and every $U, V \in B_2$, if $R(u)  \subseteq U \cup V$, then $R(u) \subseteq U$ or $R(u) \subseteq V$.
\end{thm}

To conclude, let us  consider the case where $\Sm$ satisfies both (PC) and (PDI) for a single formula. 
Then it is well known that all $\Sm$-algebras have a distributive lattice reduct (see Proposition 2.8 in \cite{GeJaPa10}) and the $\Sm$-filters are the same as the order filters of the spezialization order of the algebras in $\Alg\Sm$. 
In this case, by Propoosition \ref{prop:PCb1}, Corollary \ref{cor:PDb1}, and Proposition \ref{prop:Pr6-2} we know the following: 
if $\Al$ has a bottom element, then $\varphi[A]$ is the collection of clopen up-sets of $\langle \OptSA,\tau_{\Al},\leq\rangle$. 
Since in this case the optimal $\Sm$-filters coincide with the prime filters, 
 we obtain  exactly what Priestley duality for bounded distributive lattices gives us. 
Notice that if no bottom element is assumed, 
we still need to deal with the $\IrrSA$-admissible clopen up-sets for recovering the algebra from the space. 
This collection coincides with all clopen up-sets when the algebra has a bottom element, but  excludes the emptyset when the algebra has no bottom element.

\subsection{Deduction-Detachment Theorem.}\label{subsec:T-DCL-Corr-PDD}   

A logic $\Sm$ has \emph{the deduction-detachment property} (DDT) for a non-empty set of formulas in two variables $\Delta(x,y)$ if 
for all  $\{\delta,\gamma\}\cup\Gamma\subseteq Fm_{\mathscr{L}}$:
	 $$\Gamma,\delta\vdash_{\Sm}  \gamma \IFF \Gamma\vdash_{\Sm} \Delta(\delta,\gamma).$$	 
If $\Sm$ has (DDT) for $\Delta$, then this property transfers to every algebra in the sense that    for every algebra $\Al$, and every $\{a,b\}\cup X\subseteq A$
 \begin{equation*}
 b\in\Consa(X,a)\IFF \Delta^{\Al}(a,b)\subseteq\Consa(X),
 \end{equation*}
 see, for example,  Theorem 2.48 in \cite{FJa09}. 
A logic $\Sm$  \emph{satisfies} (uDDT) for a term $x\to y$ if it satisfies (DDT) for the set $\{x \to y\}$. 
It  is well known that (DDT) implies filter-distributivity of the logic (see \cite{Cz84}). 

The (DDT) implies the much weaker property of being protoalgebraic enjoyed for many logics. A logic $\Sm$ is \textit{protoalgebraic} if there exists  set of formulas in two variables $\Delta(x,y)$ such that
$$\text{(a) } \vdash_{\Sm} \Delta(x, x), \,\,\,\,\,\,\,\,\,
	\text{(b) }x, \Delta(x, y) \vdash_{\Sm} y. $$
	We 	will need this property in our characterization of having (uDDT).

Again we fix for  the remaining part of the subsection   \textit{a filter-distributive, finitary, and  congruential logic $\Sm$ with theorems}.

\begin{lemma}\label{lemma:PDDa} 
If $\Sm$ satisfies (uDDT), then for every  $\Sm$-algebra $\Al$ and  all $a,b\in A$, 
$({\downarrow}(\varphi(a)\cap \varphi(b)^{c}))^{c}=\varphi(a\to^{\Al} b)$.
\end{lemma}
\begin{proof}
Since (uDDT) transfers to every algebra, for any ${\{a,b\}\cup X\subseteq A}$ we have $b\in\Consa(X,a)$ if and only if $a\to^{\Al}b\in\Consa(X)$. 
Let first $P\in\varphi(a\to^{\Al}b)$, 
and suppose, towards a contradiction, that 
$P\notin ({\downarrow}(\varphi(a)\cap \varphi(b)^{c}))^{c}$. 
Then it follows that $P\in {\downarrow}(\varphi(a)\cap \varphi(b)^{c})$, 
and so there is $Q\in\OptSA$ such that $P\subseteq Q$, $Q\in\varphi(a)$ and $Q\notin\varphi(b)$. 
By assumption, from $P\subseteq Q$ we get $a\to^{\Al}b\in Q$, 
and then by (uDDT) we obtain $b\in\Consa(Q,a)$. 
Since $a\in Q$, then $b\in \Consa(Q,a)=\Consa(Q)=Q$,  a contradiction. 
We conclude that $P\in ({\downarrow}(\varphi(a)\cap \varphi(b)^{c}))^{c}$, as required.

Let now $P\in\OptSA$ be such that $P\notin \varphi(a\to^{\Al}b)$, \ie  $a\to^{\Al} b\notin P$. 
By (uDDT) we get that $b\notin \Consa(P,a)$. 
Then by the optimal $\Sm$-filter lemma, 
there is $Q\in \OptSA$ such that $b\notin Q$ and $\Consa(P,a)\subseteq Q$. 
So, we have $a\in Q$, $P\subseteq Q$ and $b\notin Q$, \ie $Q\in\varphi(a)\cap\varphi(b)^{c}$, and so $P\in{\downarrow}(\varphi(a)\cap\varphi(b)^{c})$. 
Therefore $P\notin  ({\downarrow}((\varphi(a)\cap \varphi(b)^{c}))^{c})$, as required. 
\end{proof}

From the previous corollary we conjecture that  for any \Sp space $\langle X,\tau,\BB\rangle$, 
the Priestley-dual property of (uDDT) is the property  
``$B$ is closed under $({\downarrow}((\phantom{.})\cap (\phantom{.})^{c}))^{c}$''. 
Let us check now that this condition is enough for recovering the implication.

\begin{prop}\label{prop:PDDb2} 
Let $\langle X,\tau,\BB\rangle$ be an \Sp space such that for all $U,V\in B$, $({\downarrow}(U\cap V^{c}))^{c}\in B$. 
Then for all $\{U,V\}\cup \mathcal{W}\subseteq B$:
\begin{equation*}
V\in\Consb(\mathcal{W},U)\IFF ({\downarrow}(U\cap V^{c}))^{c}\in\Consb(\mathcal{W}).
\end{equation*} 
\end{prop}
\begin{proof}
Assume first that $({\downarrow}(U\cap V^{c}))^{c}\in\Consb(\mathcal{W})$. 
Then as the logic $\Sm$ is finitary, there is a finite $\mathcal{W}'\subseteq \mathcal{W}$  such that
$({\downarrow}(U\cap V^{c}))^{c}\in\Consb(\mathcal{W'})$. 
Thus by \ref{itm:Pr2}, $\bigcap \mathcal{W}'\subseteq ({\downarrow}(U\cap V^{c}))^{c}$. 
We show that 
$U\cap \bigcap \mathcal{W}'\subseteq V$, so assume that  $u\in U\cap \bigcap \mathcal{W}'$ and suppose in view of a contradiction that $u\notin V$. 
On the one hand $u\in U$. 
Moreover $u\in \bigcap \mathcal{W}'\subseteq ({\downarrow}(U\cap V^{c}))^{c}$, 
 \ie $u\notin{\downarrow}(U\cap V^{c})$. 
But from $u\notin V$ and $x\in U$ we get $u\in U\cap V^{c}\subseteq {\downarrow}(U\cap V^{c})$, a contradiction. 
We conclude that $U\cap \bigcap \mathcal{W}'\subseteq V$, and thus by \ref{itm:Pr2}, $V\in\Consb(\mathcal{W}')\subseteq \Consb(\mathcal{W})$. 

Assume now that $V\in\Consb(\mathcal{W},U)$. 
Then by finitarity again, there is a finite $\mathcal{W}'\subseteq \mathcal{W}$ such that $V\in\Consb(\mathcal{W}',U)$. 
We show that $\bigcap \mathcal{W}'\subseteq ({\downarrow}(U\cap V^{c}))^{c}$. 
Suppose that $u\in \bigcap \mathcal{W}'$ and, 
towards a contradiction, that $u\notin ({\downarrow}(U\cap V^{c}))^{c}$. 
Then there is $v\in U\cap V^{c}$ such that $u\leq v$. 
Let $\xis(v)=\{W\in B: v\in W\}$. This set is an optimal $\Sm$-filter of $\BB$ 
by Proposition \ref{prop:xisonto-a}. Suppose that $W\in \mathcal{W}'$. 
By assumption $u\in W$ and since $W$ is an up-set, $v\in W$, \ie $W\in \xis(v)$. 
Therefore $\mathcal{W}'\subseteq \xis(v)$ and moreover, since $v\in U$, $U\in \xis(v)$. 
Furthermore, as $\xis(v)$ is an $\Sm$-filter ${\Consb(\mathcal{W}',U)\subseteq \xis(v)}$,
 so by hypothesis $V\in \xis(v)$, \ie $v\in V$, a contradiction. 
Thus we conclude that $\bigcap \mathcal{W}'\subseteq ({\downarrow}(U\cap V^{c}))^{c}$, 
and then by \ref{itm:Pr2}, $({\downarrow}(U\cap V^{c}))^{c}\in\Consb(\mathcal{W}')\subseteq \Consb(\mathcal{W})$, as required. 
\end{proof}

Assuming that $\Sm$ is protoalgebraic, we can  find  conditions over the dual space that imply that  the logic has (uDDT).  
 This result is supported in the following theorem due to Czelakowski.

 \begin{thm}[Theorem 2.6.8 in \cite{Cz01}]\label{thm:Cz01} Let $\Sm$ be a protoalgebraic logic. Then 
 $\Sm$  satisfies (DDT) if and only if for every $\Sm$-algebra $\Al$, the lattice of \mbox{$\Sm$-filters} $\lFiSA$ is infinitely meet-distributive over its compact elements, \ie for any finite $B\subseteq A$ and any $\{G_{i}:i\in I\}\subseteq \FiSA$: 
 $$\Consa(B)\sqcup \bigcap_{i\in I}G_{i}=\bigcap_{i\in I}(\Consa(B)\sqcup G_{i}).$$
 \end{thm}

\begin{thm}\label{thm:PDDb3}
Let $\Sm$ be a protoalgebraic  logic. If for every \Sp space $\langle X,\tau,\BB\rangle$, 
 $({\downarrow}(U\cap V^{c}))^c\in B$ for all $U,V\in B$, 
then $\Sm$  satisfies (uDDT). 
\end{thm}
\begin{proof}
Let $\Sm$ be a protoalgebraic logic satisfying the assumption. First we prove that $\Sm$ has (DDT) and then we will see that it satisfies (uDDT). Let $\Al$ be an $\Sm$-algebra.  
By Theorem \ref{thm:Cz01}, it is enough to show that $\FiSA$ is infinitely meet-distributive over its compact elements. 
By the representation theorem for $\Sm$-algebras, 
and Corollary \ref{cor:Salg-Prspa}, 
 we know that for any $\Sm$-algebra $\Al$ there is an \Sp space $\langle X,\tau,\BB\rangle$ such that $\Al$ is isomorphic to $\BB$. 
 Therefore, it is enough to show that for any \Sp space $\langle X,\tau,\BB\rangle$, $\FiSB$ is infinitely meet-distributive over its compact elements. 
 
 So let $\langle X,\tau,\BB\rangle$ be an \Sp space, 
let $\{G_{i}:i\in I\}\subseteq \FiSB$ and  let $U_{1},\dots,U_{n}\subseteq B$ be finite sets. 
We show that 
$$\Consb(\{U_{1},\dots,U_{n}\})\sqcup \bigcap \{G_{i}:i\in I\}=\bigcap \{ \Consb(\{U_{1},\dots,U_{n}\}\cup G_{i}):i\in I\}.$$ 
Notice that the inclusion from left to right is immediate by finitarity of the logic, so we just have to show the other inclusion. 
Let $G:= \bigcap \{G_{i}:i\in I\}$ and suppose that  
$V\in \bigcap \{ \Consb(\{U_{1},\dots,U_{n}\}\cup G_{i}):i\in I\}$. 
Then for each $i\in I$ we have that $V\in \Consb(\{U_{1},\dots,U_{n}\}\cup G_{i})$. 
For any $W_{1},W_{2}\in B$, let us denote $({\downarrow}(W_{1}\cap W_{2}^{c}))^{c}$ by $W_{1}\Rightarrow W_{2}$. 
Then for each $i\in I$, by assumption and Proposition \ref{prop:PDDb2} we get $U_{1}\Rightarrow (\dots(U_{n}\Rightarrow V)\dots)\in G_{i}$. 
Thus $U_{1}\Rightarrow (\dots(U_{n}\Rightarrow V)\dots)\in G$. 
Recall that $\bigcap \{G_{i}:i\in I\}=\Consb(\bigcap \{G_{i}:i\in I\})$ is an $\Sm$-filter of $\BB$.  
So by assumption and Proposition \ref{prop:PDDb2} again we conclude 
\begin{align*}
V\in \Consb(\{U_{1},\dots,U_{n}\}\cup G)&= \Consb\big(\Consb(\{U_{1},\dots,U_{n}\})\cup G))\\
&=\Consb(\{U_{1},\dots,U_{n}\})\sqcup \bigcap \{G_{i}:i\in I\}.
\end{align*} 
We conclude that $\Sm$ has (DDT).

Let us consider now the congruence relation  $\equiv_{\Sm}^{\Fm}$ on $Fm$, that we abbreviate all along the proof by $\equiv$. Let $\pi$ be the quotient homomorphism from $\Fm$ to $\Fm/{\equiv}$. Recall that  $\Fm/{\equiv} \in \Alg\Sm$. Let $x, y$ be two  variables
and consider the equivalence classes $x/{\equiv}$ and $y/{\equiv}$. By the assumption we have $({\downarrow}(\varphi(x/{\equiv}) \cap  \varphi(y/{\equiv})^{c})^c = \varphi(\delta/{\equiv})$ for some formula $\delta$. We prove that for every $\Gamma \subseteq Fm$
\begin{equation*}
\tag{\pseudoddt}
\label{eq:pseudo-ddt}
y \in \Cons(\Gamma, x) \IFF \delta \in \Cons(\Gamma).
\end{equation*}
Suppose that $y \in \Cons(\Gamma, x)$  and $\delta \not \in \Cons(\Gamma)$. Then there is an irreducible $\Sm$-theory $T$ such that $\Cons(\Gamma) \subseteq T$ and $\delta \not \in T$. Then, by Proposition \ref{prop:bilogic-quotient}, $\pi[T]$ is irreducible in  $\Fm/{\equiv}$ and $\pi[T] \not \in \varphi(\delta/{\equiv})$. Thus, $\pi[T] \in {\downarrow}(\varphi(x/{\equiv}) \cap  \varphi(y/{\equiv})^{c})$. Let $Q \in (\varphi(p/{\equiv}) \cap  \varphi(y/{\equiv})^{c})$ such that $\pi[T] \subseteq Q$. Then
$\pi^{-1}[Q]$ is an $\Sm$-theory such that $\Gamma \cup \{x\} \subseteq \pi^{-1}[Q]$. Therefore $q \in \pi^{-1}[Q]$ and this implies that $Q \in \varphi(y/{\equiv})$, a contradiction. To prove the converse, assume that $\delta \in \Cons(\Gamma)$ and 
$y \not \in \Cons(\Gamma, x)$. Let $T$ be an irreducible $\Sm$-theory such that $\Cons(\Gamma, x) \subseteq T$ and $y \not \in T$. Then $\pi[T] \in \varphi(x/{\equiv}) \cap  \varphi(y/{\equiv})^{c}$. Therefore, $\pi[T] \not \in \varphi(\delta/{\equiv})$; hence $\delta \not \in T$. Since $\Gamma \subseteq T$ and $\delta \in \Cons(\Gamma)$ we have a contradiction. 

By the first part of the proof let $\Delta(x, y)$ be a (DDT) set for $\Sm$. We show that $\Cons(\delta) = \Cons(\Delta)$. This easily follows from (\ref{eq:pseudo-ddt}) and the assumption that $\Delta(x, y)$ is a (DDT) set. Indeed, $y \in \Cons(\Delta, x)$ holds. Then by (\ref{eq:pseudo-ddt}) we have $\delta \in \Cons(\Delta)$. On the other hand, since $\delta \in \Cons(\delta)$,  (\ref{eq:pseudo-ddt}) gives $y \in \Cons(\delta, x)$. Therefore, $\Delta \subseteq \Cons(\delta)$.  Now let $\sigma$ be the substitution that maps $x$ to $x$ and all the remaining variables to $y$. The by invariance under substitutions follows that 
$\Cons(\Delta) = \Cons(\sigma(\delta))$. Then it easily follows that  $\Sm$ has (DDT) for the formula $\delta'(x, y) = \sigma(\delta)$. 
\end{proof}

\begin{cor}
$\Sm$ has (uDDT) if and only if $\Sm$ is protoalgebraic and for every \Sp space $\langle X,\tau,\BB\rangle$, 
 $({\downarrow}(U\cap V^{c}))^c\in B$ for all $U,V\in B$.
\end{cor}

\subsection{The Property of an Inconsistent element.} \label{subsec:T-DCL-Corr-PIE}  

A logic $\Sm$
   satisfies the \emph{property of an inconsistent element} (PIE) for a formula $\psi$ if for every
formula $\delta\in Fm_{\mathscr{L}}$:
	$$\psi\vdash_{\Sm} \delta.$$
	Such a formula is known as  an \emph{inconsistent formula}.
	It is immediate that (PIE) transfers to every algebra in the following sense. If $\Sm$ has (PIE) for $\psi$, then for every algebra $\Al$ and any homomorphism $h: \Fm \to \Al$,  
	$$a\in \Consa(h(\psi))$$ 
	for every $a \in A$.  
	If $\Sm$ is congruential, then it easily follows that if $\Sm$ satisfies (PIE) for a formula $\psi$, then for every $\Sm$-algebra $\Al$ and all $h, h' \in \Hom(\Fm, \Al)$,  $h(\psi) = h'(\psi)$, that is, $\psi$ is a constant term on $\Sm$-algebras. Moreover, it also holds that if $\Sm$ satisfies (PIE) for two inconsistent formulas   $\psi$ and $\psi'$, then their interpretations on every $\Sm$-algebra are the same. Thus if $\Sm$ satisfies (PIE), then   in every $\Sm$-algebra $\Al$ there is a unique element that is the unique possible interpretation of all the inconsistent formulas and this element is the bottom element of $\Al$ (w.r.t. $\leq^{\Al}_{\Sm}$). We denote this element by  $\inc^{\Al}$ or $0^{\Al}$ and refer to  it as  the inconsistent element of $\Al$.

For the remaining part of the subsection  let \emph{$\Sm$ be a filter-distributive, finitary, and congruential logic with theorems}.

\begin{lemma}\label{lemma:PIEa}
 If $\Sm$    satisfies (PIE), then for every $\Sm$-algebra $\Al$
and all $a\in A$, $\varphi(0^{\Al})=\emptyset\subseteq \varphi(a)$. 
\end{lemma}
\begin{proof}
Notice that, since (PIE) transfers to every algebra, 
we have that $\Al$ has a bottom element $0^{\Al}$. 
Recall that when $\Al$ has a bottom element, $\emptyset\notin \IdsS(\Al)$ and so the optimal $\Sm$-filters are proper; therefore  $\varphi(0^{\Al})=\emptyset$. 
It follows trivially that $\varphi(0^{\Al})=\emptyset \subseteq \varphi(a)$ for all $a\in A$. 
\end{proof}

\begin{thm}\label{thm:PIEb3}
If  for everyy \Sp space $\langle X,\tau,\BB\rangle$, $\emptyset\in B$,  
then $\Sm$ satisfies (PIE). 
\end{thm}
\begin{proof}
Recall that the Lindenbaum-Tarski algebra $\Fm/{\equivsfm}$ is an \mbox{$\Sm$-algebra}. Let us abbreviate $\equivsfm$ by $\equiv$. By assumption 
$\emptyset \in \varphi[\Fm/{\equivsfm}]$. Therefore there is  $\psi \in Fm_{\mathscr{L}}$ such that $\emptyset = \varphi(\psi/{\equiv})$. 
Let $\delta\in Fm_{\mathscr{L}}$. If $\psi \not \vdash_{\Sm} \delta$, there is an irreducible $\Sm$-theory $T$ such that $\psi \in T$ and $\delta \not \in T$. Then $\psi/{\equiv} \in \pi[T]$, which implies that $\pi[T] \in \varphi(\psi/{\equiv})$, a contradiction. Therefore, $\psi$ is an inconsistent formula.
\end{proof}

\begin{cor}\label{cor:PIE3b}
The logic $\Sm$ satisfies (PIE) if and only if for every \Sp space $\langle X,\tau,\BB\rangle$ it holds that $\emptyset\in B$.
\end{cor}

Observe that  when the logic $\Sm$ satisfies (PIE), 
we have that $\emptyset$ is the inconsistent element in $\BB$, 
\ie the inconsistent element in the referential algebra $\BB$ is represented by the emptyset.

When $\Sm$ satisfies both (PC) and (PIE) we know that in every $\Sm$-Priestley space $\langle X,\tau,\BB\rangle$, 
the emptyset is an  $X_{B}$-admissible clopen up-set of $X$, so $\max(X)\subseteq X_{B}$ or,  in other words, ${\downarrow} X_{B}=X$. 
This property corresponds, in the general case, to the property that the $\Sm$-algebras have a bottom-family.

\bibliographystyle{plain}

\begin{thebibliography}{10}

\bibitem{BeJa11}
G.~Bezhanishvili and R.~Jansana.
\newblock Priestley style duality for distributive meet-semilattices.
\newblock {\em Studia Logica}, 98 (2011) 83--123.

\bibitem{BlPi89}
W.~J. Blok and D.~Pigozzi.
\newblock Algebraizable logics.
\newblock {\em Mem. Amer. Math. Soc.}, 396, January 1989.

\bibitem{Ce03b}
S.~A. Celani.
\newblock Topological representation of distributive semilattices.
\newblock {\em Scientiae Mathematicae Japonicae}, 8 (2003) 41--51.

\bibitem{CeJa12}
S.~A. Celani and R.~Jansana
\newblock On the free implicative semilattice extension of a Hilbert algebra
\newblock {\em Mathematical Logic Quarterly}, 58 (2012) 188--207.

\bibitem{CiNo13}
P.~Cintula and C.~Noguera.
\newblock The proof by cases property and its variants in structural
  consequence relations.
\newblock {\em Studia Logica}, 101 (2013) 713--747.

\bibitem{Cz84}
J.~Czelakowski.
\newblock Filter distributive logics.
\newblock {\em Studia Logica}, 43 (1984) 353--377.

\bibitem{Cz01}
J.~Czelakowski.
\newblock Protoalgebraic Logics.
\newblock (Kluwer, 2001)

\bibitem{Es13}
M.~Esteban.
\newblock Duality Theory and Abstract Algebraic Logic.
\newblock PhD thesis, Universitat de Barcelona, (2013)
  http://www.tdx.cat/handle/10803/125336

\bibitem{FJa09}
J.~M. Font and R.~Jansana.
\newblock A General Algebraic Semantics for Sentential Logics, 2nd ed. volume~7
  of {\em Lectures Notes in Logic}.
\newblock (The Association for Symbolic Logic, Ithaca, N.Y.,
  2009)

\bibitem{FoJaPi03}
J.~M. Font, R.~Jansana, and D.~Pigozzi.
\newblock A survey of abstract algebraic logic.
\newblock {\em Studia Logica}, 74(1/2) (2003)13--97.

\bibitem{Fr54}
O.~Frink.
\newblock Ideals in partially ordered sets. 
\newblock {\em American Mathematical Monthly}, 61(1954) 223--234.



\bibitem{GeJaPa10}
M.~Gehrke, R.~Jansana, and A.~Palmigiano.
\newblock Canonical extensions for congruential logics with the deduction
  theorem.
\newblock {\em Annals of Pure and Applied Logic}, 161 (2010) 1502--1519.

\bibitem{Gr78}
G.~Gr{\"a}tzer.
\newblock {\em General lattice theory}.
\newblock (Academic Press. Birkh{\"a}user Verlag, 1978)

\bibitem{JaPa06}
R.~Jansana and A.~Palmigiano.
\newblock Referential semantics: duality and applications.
\newblock {\em Reports on Mathematical Logic}, 41(2006) 63--93.

\bibitem{Wo88}
R.~W\'ojcicki.
\newblock Theory of Logical Calculi. Basic Theory of Consequence
  Operations.
\newblock (Kluwer, Dordrecht, 1988)

\end{thebibliography}

\end{document}